\documentclass[a4paper,twoside]{amsart}
\usepackage{xypic,amscd,amssymb,latexsym}
\makeatletter
\def\@setcopyright{}
\makeatother
\renewcommand{\-}{\hbox{-}}
\newcommand{\X}{\raisebox{0.43ex}{$\chi$}}
\newcommand{\eins}{\mbox{\rm 1\hspace{-0.24em}l}}
\renewcommand{\mod}{\operatorname{mod}\nolimits}
\newcommand{\Sub}{\operatorname{Sub}\nolimits}
\newcommand{\add}{\operatorname{add}\nolimits}
\newcommand{\Sets}{\operatorname{Sets}\nolimits}
\newcommand{\Hom}{\operatorname{Hom}\nolimits}
\newcommand{\End}{\operatorname{End}\nolimits}
\newcommand{\Aut}{\operatorname{Aut}\nolimits}
\newcommand{\Irr}{\operatorname{Irr}\nolimits}
\newcommand{\Inn}{\operatorname{Inn}\nolimits}
\newcommand{\Out}{\operatorname{Out}\nolimits}
\renewcommand{\Im}{\operatorname{Im}\nolimits}
\newcommand{\Ker}{\operatorname{Ker}\nolimits}
\newcommand{\Coker}{\operatorname{Coker}\nolimits}
\newcommand{\rad}{\operatorname{rad}\nolimits}
\newcommand{\soc}{\operatorname{soc}\nolimits}
\newcommand{\ann}{\operatorname{ann}\nolimits}
\newcommand{\Tor}{\operatorname{Tor}\nolimits}
\newcommand{\Tr}{\operatorname{Tr}\nolimits}
\newcommand{\GL}{\operatorname{GL}\nolimits}
\newcommand{\GF}{\operatorname{GF}\nolimits}
\newcommand{\Mat}{\operatorname{Mat}\nolimits}
\newcommand{\PSL}{\operatorname{PSL}\nolimits}
\newcommand{\SL}{\operatorname{SL}\nolimits}
\newcommand{\Fi}{\operatorname{Fi}\nolimits}
\newcommand{\Co}{\operatorname{Co}\nolimits}
\newcommand{\J}{\operatorname{J}\nolimits}
\newcommand{\kar}{\operatorname{char}\nolimits}
\newcommand{\Ext}{\operatorname{Ext}\nolimits}
\newcommand{\op}{{\operatorname{op}}}
\newcommand{\tr}{{\operatorname{tr}}}
\newcommand{\Ab}{\operatorname{Ab}\nolimits}
\newcommand{\id}{\operatorname{id}\nolimits}
\newcommand{\comp}{\operatorname{\scriptstyle\circ}}
\newcommand{\E}{\operatorname{\mathcal E}\nolimits}
\newcommand{\I}{\operatorname{\mathcal I}\nolimits}
\newcommand{\M}{\operatorname{\mathcal M}\nolimits}
\renewcommand{\P}{\operatorname{\mathcal P}\nolimits}
\newcommand{\R}{\operatorname{\mathfrak{R}}\nolimits}
\newcommand{\Y}{\operatorname{\mathcal Y}\nolimits}
\newcommand{\w}{\operatorname{\mathcal W}\nolimits}
\newcommand{\G}{\Gamma}
\renewcommand{\L}{\Lambda}
\newcommand{\la}{\lambda}
\renewcommand{\r}{\operatorname{\underline{r}}\nolimits}
\renewcommand{\a}{\underline{\mathrm a}}
\renewcommand{\b}{\underline{\mathrm b}}
\newcommand{\m}{\underline{\mathrm m}}
\newcommand{\dwnsimeq}{\Big\downarrow\lefteqn{\wr}}
\newcommand{\upsimeq}{\Big\uparrow\lefteqn{\wr}}
\newcommand{\subneq}{\hspace*{2pt}\mathrel
{{\raisebox{-.3ex}{\hbox{$\scriptscriptstyle\not$}}}\!\!\!\subseteq}}
\newcommand{\semi}{\mathbin{\vcenter{\hbox{$\scriptscriptstyle|$}}\;\!\!\!\times}}
\newcommand{\longequals}{\Relbar\joinrel=}
\newcommand{\A}{\alpha}
\newcommand{\B}{\beta}
\newcommand{\C}{\gamma}
\newcommand{\D}{\delta}
\newcommand{\T}{\tau}
\newcommand{\Z}{\xi}

\newtheorem{lemma}{Lemma}[section]
\newtheorem{proposition}[lemma]{Proposition}
\newtheorem{corollary}[lemma]{Corollary}
\newtheorem{theorem}[lemma]{Theorem}

\newtheorem{cclass}[lemma]{}
\newtheorem{ctab}[lemma]{}

\renewcommand{\arraystretch}{1.5}

\begin{document}

\thispagestyle{empty}
\title{CONSTRUCTION OF $\Co_1$ FROM AN IRREDUCIBLE SUBGROUP $\M_{24}$ of $\GL_{11}(2)$}

\author{Hyun Kyu Kim}
\address{Department of Mathematics, Yale University, New Haven,
  CT. 06511, USA}

\author{Gerhard O. Michler}
\address{Department of Mathematics, Cornell University, Ithaca,
  N.Y. 14853, USA}

\begin{abstract}
In this article we give an self contained existence proof for J. Conway's
sporadic simple group $\Co_1$ \cite{conway1} using
the second author's algorithm \cite{michler1} constructing finite
simple groups from irreducible subgroups of $\GL_n(2)$. Here $n = 11$ and
the irreducible subgroup is the Mathieu group $\M_{24}$. From the split extension $E$ of $\M_{24}$ by a uniquely determined $11$-dimensional $\GF(2)\M_{24}$-module $V$
we construct the centralizer $H = C_G(z)$ of a
$2$-central involution $z$ of $E$ in an unknown target group $G$. Then we prove that all the
conditions of Algorithm 2.5 of \cite{michler1} are satisfied. This
allows us to construct a simple subgroup $G$ of $\GL_{276}(23)$ which we prove to be isomorphic with Conway's original sporadic simple group $\Co_1$ by means of a constructed faithful permutation representation of $G$ and Soicher's presentation \cite{soicher1} of the original Conway group $\Co_1$.
\end{abstract}

\maketitle

\section{Introduction}
In $1969$ J.H. Conway \cite{conway1} discovered $3$ sporadic
simple groups which he defined in terms of the automorphism group
$A = Aut(\Lambda)$ of the $24$-dimensional Leech lattice
$\Lambda$, see also \cite{conway}. The center $Z(A)$ of $A$ has
order $2$, and $\Co_1 = A/Z(A)$ is the largest of these $3$ simple
groups.

The results of this article are part of our joint research project
{\em Simultaneous construction of the sporadic simple groups of
Conway, Fischer and Janko}. Its goal is to provide uniform
existence proofs for the sporadic simple groups discovered by Conway, Fischer and Janko by
means of Algorithm 2.5 of \cite{michler1} constructing finite
simple groups from irreducible subgroups $T$ of $\GL_n(2)$. In
\cite{kim} we constructed Conway's and Fischer's sporadic groups $\Co_2$ and $\Fi_{22}$ simultaneously
from the irreducible subgroup $\M_{22}$ in $\GL_{10}(2)$. In \cite{kim1} the
first author applied the same methods to the
irreducible subgroup $\M_{23}$ of $\GL_{11}(2)$ and realized
$\Fi_{23}$ as an irreducible subgroup of $\GL_{782}(17)$. The authors gave such an existence proof of Fischer's sporadic simple group $\Fi_{24}'$ in \cite{kim2}. Janko's large sporadic group $\J_4$ is dealt with in \cite{michler2}.

In section $2$ we construct the split extension of $\M_{24}$ by its uniquely determined $11$-dimensional irreducible $\GF(2)$-module $V_1$ satisfying $dim_{\GF(2)}[H^2(M_{24},V_1)] = 0$. A presentation of $E$ is stated in Lemma \ref{l.
M24-extensions}. We determine conjugacy classes and character table of $E$. It follows that $E$ has a unique conjugacy class $z^E$ of $2$-central involutions.

In section $3$ we apply Algorithm 2.5 of \cite{michler1} to the
extension group $E$. Therefore we construct $D = C_E(z)$. Using a faithful
permutation representation $PE$ of $E$ of degree $2^{11}$ we find a uniquely determined
non abelian normal subgroup $Q$ of $D$ with center $Z(Q) = Z(D) = \langle z \rangle$ such that $V = Q/Z(Q)$ is
elementary abelian 
normal subgroup of order $2^8$ in $D_1 = D/Z(Q)$. Furthermore, $MD
= D_1/V$ has an elementary abelian normal Fitting subgroup $MY$ of
order $2^6$ which has a complement $MC \cong A_8$ in $MD$.
Furthermore, $MY$ contains a $2$-central involution $Mt$ of $MD$.
Since the alternating group $A_8$ has an irreducible
representation of degree $6$ over $F = \GF(2)$ we apply Algorithm
2.5 of \cite{michler1} again and find a matrix $Mh \in \GL_8(2)$
of order $3$ such that $C_{MO}(Mt) = \langle C_{MD}(Mt), Mh
\rangle$ where $MO = \langle MD, Mh\rangle \cong O^{+}_8(2)$.
Furthermore, $|MO : MD|$ is odd and all conditions of Step 5 of
Algorithm 2.5 of \cite{michler1} are satisfied.

Since $D_1$ does not split over $V$ we have to construct all non split extensions $H^{*}$ of $MO$ by
$V$. Unfortunately, the implementation of Holt's Algorithm \cite{holt5} in MAGMA was not able to calculate the dimension of the second cohomology group $[H^2(MO,V)]$. Therefore we
construct an extension $H_2$ of the centralizer $C_{MO}(Mt)$ by $V$ as a finitely
presented group. Furthermore, we show that $H_2$ has
a central extension $H_1$ by $Z(Q) = \langle z \rangle$ such that its Sylow $2$-subgroup are isomorphic to the Sylow $2$-subgroups of $D$.
Moreover, we found a pair of isomorphic subgroups $U$ and $U_1$ of
$D$ and $H_1$, respectively. The amalgam $D \leftarrow U
\rightarrow H_1$ has then been used to construct a matrix subgroup
$MH$ of $\GL_{128}(23)$ which has an extra-special normal subgroup
$MQ$ of order $2^9$ such that $MH/MQ \cong O^{+}_8(2)$ and $MQ \cong
Q$. A presentation and a faithful permutation representation of
degree $61440$ of this group $H$ are stated in Propositions
\ref{prop. H(Co_1)} and \ref{prop. H(Co_2 in Co_1)}, respectively. Furthermore, all conditions of Step 5 of
Algorithm 2.5 of \cite{michler1} are verified there.

In section 4 we apply Algorithm 7.4.8 of \cite{michler} to this
group $H$ of even order. It returns a simple subgroup $\mathfrak G$ of
$\GL_{276}(23)$ of order $2^{21}\cdot 3^9\cdot 5^4\cdot 7^2\cdot
11\cdot 13\cdot 23$. We construct a faithful permutation representation
of degree $98280$ with a documented stabilizer. Using it and MAGMA we obtain the character table of $\mathfrak G$. It agrees with that of Conway's sporadic group $\Co_1$. All these results are proved in
Theorem \ref{thm. existenceCo_1}. Using then L. Soicher's presentation \cite{soicher1} and Theorem 7.5.1 of \cite{michler} we show in Corollary \ref{cor. identCo_1} that $\mathfrak G \cong \Co_1$. Furthermore, we describe a faithful permutation representation of degree $98280$ of Soicher's finitely presented group $G$ defined in \cite{soicher1}. Its stabilizer $U$ is isomorphic to Conway's sporadic simple group $\Co_2$ constructed in \cite{kim}.

In the appendix we collect the systems of representatives of
conjugacy classes in terms of the given generators of the local
subgroups of $G$ which have been used to construct the matrix
group $G \cong \Co_1$. We also state the character tables of
these subgroups. The four generating matrices of the simple
subgroup $\mathfrak G = \langle \mathfrak x, \mathfrak y,
\mathfrak h, \mathfrak e\rangle$ of $\GL_{276}(23)$ and its faithful permutation representation $P\mathfrak G$ can be
downloaded from the first author's website\\
$\verb"http://www.math.yale.edu/~hk47/Co1/index.html"$.

Concerning our notation and terminology we refer to the books
\cite{carter}, \cite{holt2} and \cite{michler}. The computer algebra system
MAGMA is described in Cannon-Playoust \cite{magma}.

Acknowledgements: The first author kindly acknowledges financial
support by the Hunter R. Rawlings III Cornell Presidential
Research Scholars Program for participation in the research
project {\em Simultaneous construction of the sporadic simple
groups of Conway, Fischer and Janko}. This collaboration has also
been supported by the grant NSF/SCREMS DMS-0532/06.

\newpage

\section{Extensions of Mathieu group $\M_{24}$}

J.A. Todd's presentations and permutation representations of the Mathieu groups
are stated in Definition 8.2.1 and Lemma 8.2.2 of \cite{michler}, respectively.
The irreducible $2$-modular representations of the Mathieu group $\M_{24}$ were determined by
G. James \cite{james}. Therefore all conditions of Holt's Algorithm 7.4.5 of \cite{michler} implemented
in MAGMA are satisfied. It constructs all split and non split extensions of $\M_{24}$ by the two irreducible $11$-dimensional modules $V_1$ and $V_2$ of $\M_{24}$ up to isomorphism. Here we only describe the split extension $E$ of $\M_{24}$ by $V_1$.

\begin{lemma}\label{l. M24-extensions}
Let $\M_{24} = \langle a,b,c,d,t,g,h,i,j,k \rangle$ be the
finitely presented group of Definition 8.2.1 of \cite{michler}. Let ${\mathcal R}(\M_{24})$ be its set of defining relations.
Then the following statements hold:

\begin{enumerate}
\item[\rm(a)] A faithful permutation representation $PM_{24}$ of degree $24$
of $\M_{24}$ is stated in Lemma 8.2.2 of \cite{michler}.

\item[\rm(b)] The first $2$-modular irreducible representation $V_1$ of
$\M_{24}$ is described by the following matrices:
{\renewcommand{\arraystretch}{0.5}
$$
a_1 = \left( \begin{array}{*{11}{c@{\,}}c}
1 & 0 & 0 & 0 & 0 & 0 & 0 & 0 & 0 & 0 & 0\\
0 & 1 & 0 & 0 & 0 & 0 & 0 & 0 & 1 & 0 & 0\\
0 & 0 & 1 & 0 & 0 & 0 & 0 & 0 & 1 & 0 & 0\\
0 & 0 & 0 & 1 & 0 & 0 & 0 & 0 & 1 & 0 & 0\\
0 & 0 & 0 & 0 & 0 & 1 & 0 & 0 & 0 & 0 & 0\\
0 & 0 & 0 & 0 & 1 & 0 & 0 & 0 & 0 & 0 & 0\\
0 & 0 & 0 & 0 & 0 & 0 & 0 & 0 & 0 & 1 & 0\\
0 & 0 & 0 & 0 & 0 & 0 & 0 & 0 & 1 & 0 & 1\\
0 & 0 & 0 & 0 & 0 & 0 & 0 & 0 & 1 & 0 & 0\\
0 & 0 & 0 & 0 & 0 & 0 & 1 & 0 & 0 & 0 & 0\\
0 & 0 & 0 & 0 & 0 & 0 & 0 & 1 & 1 & 0 & 0
\end{array} \right),\quad
b_1 = \left( \begin{array}{*{11}{c@{\,}}c}
1 & 0 & 0 & 0 & 0 & 0 & 0 & 1 & 1 & 0 & 1\\
0 & 1 & 0 & 0 & 0 & 0 & 0 & 1 & 1 & 0 & 1\\
0 & 0 & 1 & 0 & 0 & 0 & 0 & 1 & 1 & 0 & 1\\
0 & 0 & 0 & 1 & 0 & 0 & 0 & 0 & 1 & 0 & 0\\
0 & 0 & 0 & 0 & 0 & 0 & 1 & 1 & 0 & 0 & 0\\
0 & 0 & 0 & 0 & 0 & 0 & 0 & 0 & 1 & 1 & 1\\
0 & 0 & 0 & 0 & 1 & 0 & 0 & 0 & 0 & 0 & 1\\
0 & 0 & 0 & 0 & 0 & 0 & 0 & 0 & 0 & 0 & 1\\
0 & 0 & 0 & 0 & 0 & 0 & 0 & 0 & 1 & 0 & 0\\
0 & 0 & 0 & 0 & 0 & 1 & 0 & 1 & 1 & 0 & 0\\
0 & 0 & 0 & 0 & 0 & 0 & 0 & 1 & 0 & 0 & 0
\end{array} \right),
$$
}

{\renewcommand{\arraystretch}{0.5}
$$
c_1 = \left( \begin{array}{*{11}{c@{\,}}c}
1 & 0 & 0 & 0 & 0 & 0 & 1 & 0 & 0 & 1 & 0\\
0 & 1 & 0 & 0 & 0 & 0 & 1 & 0 & 0 & 1 & 0\\
0 & 0 & 1 & 0 & 0 & 0 & 1 & 0 & 1 & 1 & 0\\
0 & 0 & 0 & 1 & 0 & 0 & 0 & 0 & 1 & 0 & 0\\
0 & 0 & 0 & 0 & 0 & 0 & 1 & 1 & 1 & 0 & 0\\
0 & 0 & 0 & 0 & 0 & 0 & 0 & 0 & 0 & 1 & 1\\
0 & 0 & 0 & 0 & 0 & 0 & 0 & 0 & 1 & 1 & 0\\
0 & 0 & 0 & 0 & 1 & 0 & 0 & 0 & 0 & 1 & 0\\
0 & 0 & 0 & 0 & 0 & 0 & 0 & 0 & 1 & 0 & 0\\
0 & 0 & 0 & 0 & 0 & 0 & 1 & 0 & 1 & 0 & 0\\
0 & 0 & 0 & 0 & 0 & 1 & 1 & 0 & 1 & 0 & 0
\end{array} \right),\quad
d_1 = \left( \begin{array}{*{11}{c@{\,}}c}
1 & 0 & 0 & 0 & 0 & 0 & 1 & 0 & 0 & 1 & 0\\
0 & 1 & 0 & 0 & 0 & 1 & 1 & 0 & 0 & 0 & 1\\
0 & 0 & 1 & 0 & 0 & 1 & 1 & 1 & 0 & 0 & 0\\
0 & 0 & 0 & 1 & 0 & 1 & 1 & 1 & 0 & 0 & 0\\
0 & 0 & 0 & 0 & 0 & 0 & 0 & 1 & 1 & 1 & 0\\
0 & 0 & 0 & 0 & 0 & 0 & 1 & 0 & 0 & 0 & 1\\
0 & 0 & 0 & 0 & 0 & 0 & 0 & 1 & 0 & 1 & 1\\
0 & 0 & 0 & 0 & 0 & 1 & 0 & 0 & 0 & 1 & 0\\
0 & 0 & 0 & 0 & 1 & 1 & 1 & 1 & 0 & 1 & 1\\
0 & 0 & 0 & 0 & 0 & 0 & 1 & 1 & 0 & 0 & 1\\
0 & 0 & 0 & 0 & 0 & 1 & 0 & 1 & 0 & 1 & 1
\end{array} \right),
$$
}

{\renewcommand{\arraystretch}{0.5}
$$
t_1 = \left( \begin{array}{*{11}{c@{\,}}c}
1 & 0 & 0 & 0 & 0 & 0 & 1 & 0 & 0 & 1 & 0 \\
0 & 1 & 0 & 0 & 0 & 1 & 0 & 0 & 1 & 1 & 0 \\
0 & 0 & 1 & 0 & 0 & 0 & 0 & 1 & 0 & 1 & 0 \\
0 & 0 & 0 & 1 & 0 & 0 & 0 & 0 & 1 & 0 & 0 \\
0 & 0 & 0 & 0 & 1 & 0 & 1 & 1 & 1 & 1 & 0 \\
0 & 0 & 0 & 0 & 0 & 1 & 0 & 1 & 1 & 0 & 0 \\
0 & 0 & 0 & 0 & 0 & 1 & 1 & 0 & 1 & 0 & 0 \\
0 & 0 & 0 & 0 & 0 & 0 & 1 & 1 & 1 & 0 & 0 \\
0 & 0 & 0 & 0 & 0 & 1 & 1 & 1 & 1 & 0 & 0 \\
0 & 0 & 0 & 0 & 0 & 1 & 0 & 1 & 0 & 1 & 1 \\
0 & 0 & 0 & 0 & 0 & 1 & 1 & 0 & 0 & 1 & 0
\end{array} \right),\quad
g_1 = \left( \begin{array}{*{11}{c@{\,}}c}
1 & 0 & 0 & 0 & 0 & 0 & 0 & 0 & 0 & 1 & 1 \\
0 & 1 & 0 & 0 & 0 & 0 & 0 & 0 & 1 & 0 & 0 \\
0 & 0 & 1 & 0 & 0 & 0 & 0 & 0 & 0 & 0 & 1 \\
0 & 0 & 0 & 0 & 1 & 0 & 0 & 0 & 1 & 1 & 1 \\
0 & 0 & 0 & 1 & 0 & 0 & 0 & 0 & 1 & 1 & 1 \\
0 & 0 & 0 & 0 & 0 & 1 & 0 & 0 & 1 & 1 & 1 \\
0 & 0 & 0 & 0 & 0 & 0 & 1 & 0 & 1 & 1 & 0 \\
0 & 0 & 0 & 0 & 0 & 0 & 0 & 1 & 1 & 0 & 1 \\
0 & 0 & 0 & 0 & 0 & 0 & 0 & 0 & 1 & 0 & 0 \\
0 & 0 & 0 & 0 & 0 & 0 & 0 & 0 & 0 & 1 & 0 \\
0 & 0 & 0 & 0 & 0 & 0 & 0 & 0 & 0 & 0 & 1
\end{array} \right),
$$
}

{\renewcommand{\arraystretch}{0.5}
$$
h_1 = \left( \begin{array}{*{11}{c@{\,}}c}
1 & 0 & 0 & 0 & 0 & 0 & 1 & 0 & 0 & 1 & 0 \\
0 & 1 & 0 & 0 & 0 & 0 & 0 & 0 & 0 & 1 & 0 \\
0 & 0 & 0 & 1 & 0 & 0 & 0 & 0 & 0 & 1 & 0 \\
0 & 0 & 1 & 0 & 0 & 0 & 0 & 0 & 0 & 1 & 0 \\
0 & 0 & 0 & 0 & 1 & 0 & 1 & 0 & 0 & 0 & 0 \\
0 & 0 & 0 & 0 & 0 & 1 & 0 & 0 & 0 & 1 & 0 \\
0 & 0 & 0 & 0 & 0 & 0 & 1 & 0 & 0 & 0 & 0 \\
0 & 0 & 0 & 0 & 0 & 0 & 1 & 0 & 0 & 0 & 1 \\
0 & 0 & 0 & 0 & 0 & 0 & 1 & 0 & 1 & 1 & 0 \\
0 & 0 & 0 & 0 & 0 & 0 & 0 & 0 & 0 & 1 & 0 \\
0 & 0 & 0 & 0 & 0 & 0 & 1 & 1 & 0 & 0 & 0
\end{array} \right),\quad
i_1 = \left( \begin{array}{*{11}{c@{\,}}c}
1 & 0 & 0 & 0 & 0 & 0 & 0 & 0 & 0 & 0 & 0 \\
0 & 0 & 1 & 0 & 0 & 0 & 1 & 1 & 0 & 0 & 0 \\
0 & 1 & 0 & 0 & 0 & 1 & 1 & 0 & 1 & 0 & 0 \\
0 & 0 & 0 & 1 & 0 & 1 & 1 & 1 & 0 & 0 & 0 \\
0 & 0 & 0 & 0 & 1 & 1 & 0 & 1 & 1 & 0 & 0 \\
0 & 0 & 0 & 0 & 0 & 0 & 0 & 1 & 1 & 0 & 0 \\
0 & 0 & 0 & 0 & 0 & 1 & 0 & 1 & 0 & 0 & 0 \\
0 & 0 & 0 & 0 & 0 & 0 & 1 & 1 & 1 & 0 & 0 \\
0 & 0 & 0 & 0 & 0 & 1 & 1 & 1 & 1 & 0 & 0 \\
0 & 0 & 0 & 0 & 0 & 1 & 0 & 0 & 1 & 0 & 1 \\
0 & 0 & 0 & 0 & 0 & 1 & 1 & 0 & 0 & 1 & 0
\end{array} \right),
$$
} {\renewcommand{\arraystretch}{0.5}
$$
j_1 = \left( \begin{array}{*{11}{c@{\,}}c}
0 & 1 & 0 & 0 & 0 & 1 & 0 & 0 & 0 & 0 & 0 \\
1 & 0 & 0 & 0 & 0 & 1 & 0 & 0 & 0 & 0 & 0 \\
0 & 0 & 1 & 0 & 0 & 1 & 0 & 0 & 1 & 0 & 0 \\
0 & 0 & 0 & 1 & 0 & 0 & 0 & 0 & 1 & 0 & 0 \\
0 & 0 & 0 & 0 & 1 & 0 & 0 & 0 & 1 & 0 & 0 \\
0 & 0 & 0 & 0 & 0 & 1 & 0 & 0 & 0 & 0 & 0 \\
0 & 0 & 0 & 0 & 0 & 1 & 1 & 0 & 1 & 0 & 0 \\
0 & 0 & 0 & 0 & 0 & 1 & 0 & 1 & 0 & 0 & 0 \\
0 & 0 & 0 & 0 & 0 & 0 & 0 & 0 & 1 & 0 & 0 \\
0 & 0 & 0 & 0 & 0 & 0 & 0 & 0 & 1 & 0 & 1 \\
0 & 0 & 0 & 0 & 0 & 0 & 0 & 0 & 1 & 1 & 0
\end{array} \right)\quad \text{and}~\quad
k_1 = \left( \begin{array}{*{11}{c@{\,}}c}
1 & 1 & 1 & 1 & 1 & 1 & 1 & 0 & 1 & 1 & 1 \\
0 & 1 & 0 & 0 & 0 & 1 & 0 & 1 & 0 & 0 & 0 \\
0 & 0 & 1 & 0 & 0 & 1 & 1 & 1 & 0 & 0 & 0 \\
0 & 0 & 0 & 1 & 0 & 1 & 1 & 1 & 0 & 0 & 0 \\
0 & 0 & 0 & 0 & 1 & 0 & 1 & 0 & 0 & 0 & 0 \\
0 & 0 & 0 & 0 & 0 & 0 & 0 & 1 & 0 & 0 & 0 \\
0 & 0 & 0 & 0 & 0 & 0 & 1 & 0 & 0 & 0 & 0 \\
0 & 0 & 0 & 0 & 0 & 1 & 0 & 0 & 0 & 0 & 0 \\
0 & 0 & 0 & 0 & 0 & 1 & 1 & 1 & 1 & 0 & 0 \\
0 & 0 & 0 & 0 & 0 & 0 & 1 & 1 & 0 & 0 & 1 \\
0 & 0 & 0 & 0 & 0 & 1 & 1 & 0 & 0 & 1 & 0
\end{array} \right).
$$
}



\item[\rm(c)] $dim_F[H^2(\M_{24},V_1)] = 0$.


\item[\rm(d)] The split extension $$E = E(Co_1) = \langle
a,b,c,d,t,g,h,i,j,k,v_1,v_2,v_3,v_4,v_5,v_6,v_8,v_8,v_9,v_{10},v_{11}
\rangle$$ of $\M_{24}$ by $V_1$ has a set ${\mathcal R}(E)$ of
defining relations consisting of ${\mathcal R}(\M_{24})$ and the following set
of relations:

\begin{eqnarray*}
&&v_i^2 = 1 \quad \mbox{and} \quad v_kv_j = v_jv_k \quad \mbox{for all} \quad 1 \le i, j, k \le 11, \quad av_1a^{-1}v_1 = 1,\\
&& av_2a^{-1}v_2v_9 = av_3a^{-1}v_3v_9 = av_4a^{-1}v_4v_9 = av_5a^{-1}v_6 = av_6a^{-1}v_5 = 1,\\
&& av_7a^{-1}v_{10} = av_8a^{-1}v_9v_{11} = av_9a^{-1}v_9 = av_{10}a^{-1}v_7 = av_{11}a^{-1}v_8v_9 = 1,\\
&& bv_1b^{-1}v_1v_8v_9v_{11} = bv_2b^{-1}v_2v_8v_9v_{11} = bv_3b^{-1}v_3v_8v_9v_{11} = bv_4b^{-1}v_4v_9 = 1,\\
&& bv_5b^{-1}v_7v_8 = bv_6b^{-1}v_9v_{10}v_{11} = bv_7b^{-1}v_5v_{11} = bv_8b^{-1}v_{11} = 1,\\
&& bv_9b^{-1}v_9 = bv_{10}b^{-1}v_6v_8v_9 = bv_{11}b^{-1}v_8 =1,\\
&& cv_1c^{-1}v_1v_7v_{10} = cv_2c^{-1}v_2v_7v_{10} = cv_3c^{-1}v_3v_7v_9v_{10} = cv_4c^{-1}v_4v_9 = 1,\\
&& cv_5c^{-1}v_7v_8v_9 = cv_6c^{-1}v_{10}v_{11} = cv_7c^{-1}v_9v_{10} = cv_8c^{-1}v_5v_{10} =1,\\
&& cv_9c^{-1}v_9 = cv_{10}c^{-1}v_7v_9 = cv_{11}c^{-1}v_6v_7v_9 = 1,\\
&& dv_1d^{-1}v_1v_7v_{10} = dv_2d^{-1}v_2v_6v_7v_{11} = dv_3d^{-1}v_3v_6v_7v_8 = dv_4d^{-1}v_4v_6v_7v_8 = 1,\\
&& dv_5d^{-1}v_8v_9v_{10} = dv_6d^{-1}v_7v_{11} = dv_7d^{-1}v_8v_{10}v_{11} = dv_8d^{-1}v_6v_{10} =1,\\
&& dv_9d^{-1}v_5v_6v_7v_8v_{10}v_{11} = dv_{10}d^{-1}v_7v_8v_{11} = dv_{11}d^{-1}v_6v_8v_{10}v_{11} = 1,\\
&& tv_1t^{-1}v_1v_8v_9v_{11} = tv_2t^{-1}v_2v_7v_8v_9v_{11} = tv_3t^{-1}v_3v_6v_7v_8v_9v_{11} =1,\\
&& tv_4t^{-1}v_4v_6v_7v_8 = tv_5t^{-1}v_5v_6v_{11} = tv_6t^{-1}v_8v_9 = tv_7t^{-1}v_6v_9 = tv_8t^{-1}v_7v_9 = 1,\\
&& tv_9t^{-1}v_6v_7v_8 = tv_{10}t^{-1}v_6v_8v_{11} = tv_{11}t^{-1}v_6v_7v_{10}v_{11} = gv_1g^{-1}v_1v_{10}v_{11} = 1,\\
&& gv_2g^{-1}v_2v_9 = gv_3g^{-1}v_3v_{11} = gv_4g^{-1}v_5v_9v_{10}v_{11} = gv_5g^{-1}v_4v_9v_{10}v_{11} = 1,\\
&& gv_6g^{-1}v_6v_9v_{10}v_{11} = gv_7g^{-1}v_7v_9v_{10} = gv_8g^{-1}v_8v_9v_{11} = gv_9g^{-1}v_9 = 1,\\
&& gv_{10}g^{-1}v_{10} = gv_{11}g^{-1}v_{11} = 1,\\
&& hv_1h^{-1}v_1v_7v_{10} = hv_2h^{-1}v_2v_{10} = hv_3h^{-1}v_4v_{10} = hv_4h^{-1}v_3v_{10} = 1,\\
\end{eqnarray*}
\begin{eqnarray*}
&& hv_5h^{-1}v_5v_7 = hv_6h^{-1}v_6v_{10} = hv_7h^{-1}v_7 = hv_8h^{-1}v_7v_{11} = 1,\\
&& hv_9h^{-1}v_7v_9v_{10} = hv_{10}h^{-1}v_{10} = hv_{11}h^{-1}v_7v_8 = 1,\\
&& iv_1i^{-1}v_1 = iv_2i^{-1}v_3v_7v_8 = iv_3i^{-1}v_2v_6v_7v_9 =  iv_4i^{-1}v_4v_6v_7v_8 = 1,\\
&& iv_5i^{-1}v_5v_6v_8v_9 = iv_6i^{-1}v_8v_9 = iv_7i^{-1}v_6v_8 = iv_8i^{-1}v_7v_8v_9 = 1,\\
&& iv_9i^{-1}v_6v_7v_8v_9 = iv_{10}i^{-1}v_6v_9v_{11} = iv_{11}i^{-1}v_6v_7v_{10} = jv_6j^{-1}v_6 = 1,\\
&& jv_1j^{-1}v_2v_6 = jv_2j^{-1}v_1v_6 = jv_3j^{-1}v_3v_6v_9 = jv_4j^{-1}v_4v_9 = jv_5j^{-1}v_5v_9  = =1,\\
&& jv_7j^{-1}v_6v_7v_9 = jv_8j^{-1}v_6v_8 = jv_9j^{-1}v_9 = jv_{10}j^{-1}v_9v_{11} = jv_{11}j^{-1}v_9v_{10} = 1,\\
&& kv_1k^{-1}v_1v_2v_3v_4v_5v_6v_7v_9v_{10}v_{11} = kv_2k^{-1}v_2v_6v_8 = kv_3k^{-1}v_3v_6v_7v_8 = 1,\\
&& kv_4k^{-1}v_4v_6v_7v_8 = kv_5k^{-1}v_5v_7 = kv_6k^{-1}v_8 = kv_7k^{-1}v_7 = kv_8k^{-1}v_6 = 1,\\
&& kv_9k^{-1}v_6v_7v_8v_9 = kv_{10}k^{-1}v_7v_8v_{11} = kv_{11}k^{-1}v_6v_7v_{10}=1.\\
\end{eqnarray*}

\item[\rm(e)] The split extension $E$ has a
faithful permutation representation $PE$ of degree
$2048$ with stabilizer $T = \M_{24}$.

\item[\rm(f)] $E$ has $80$ conjugacy classes. The element $z_1 =v_1$
represents the unique conjugacy classes of $2$-central involutions
of $E$.

\item[\rm(g)] $V_1$ is a self centralizing maximal elementary
abelian normal subgroup of order $2^{11}$ of each Sylow
$2$-subgroup $S$ of $E$.
\end{enumerate}
\end{lemma}

\begin{proof}
The $2$ irreducible $F\M_{24}$-modules $V_i$, $i =1,2$, occur as
composition factors with multiplicity $1$ in the permutation
module $(1_{\M_{23}})^{\M_{24}}$ and can easily be constructed
using the faithful permutation representation of $\M_{24}$ stated
in (a) and the Meataxe algorithm implemented in MAGMA. The
corresponding matrices of the generators of $\M_{24}$ with respect
to the first irreducible representation of $\M_{24}$ are stated in
(b).

(c) The cohomological dimension $d_1 = dim_F[H^2(\M_{24},V_1)]$ has been calculated by means of MAGMA using Holt's
Algorithm 7.4.5 of \cite{michler}, the presentation of $\M_{24}$
of Definition 8.2.1 of \cite{michler} and all the data stated in
(a) and (b). It follows that $d_1 = 0$.

(d) Since $d_1 = 0$ there exists only the split extension $E$ of
$\M_{24}$ by $V_1$. The presentation of $E$ has been obtained
automatically by application of Lemma 1.4.7 of \cite{michler2} and
MAGMA.

(e) This assertion is an immediate consequence of (b) and (c).

(f) Using Kratzer's Algorithm 5.3.18 of \cite{michler},
the faithful permutation representation $PE$ and MAGMA we calculated a system of representatives of the $80$ conjugacy classes of $E$. It follows that $z_1 = v_1$ represents the unique conjugacy classes of
$2$-central involutions of $E$.

(g) This statement has been checked computationally by means of the
faithful permutation representation $PE$ and MAGMA.
\end{proof}

\newpage

\section{Centralizer of a $2$-central involution of Conway group $\Co_1$}

In this section we apply Step $5$ of Algorithm 2.5 of
\cite{michler1} to the extension group $E$ to provide a new
existence proof for Conway's largest sporadic simple group $\Co_1$.

Let $E$ be the extension group constructed in Lemma \ref{l.
M24-extensions}. It has a unique conjugacy class $z^{G}$ of $2$-central
involutions. Let $D = C_E(z)$. Using a faithful permutation
representation of $E$ we find a uniquely determined non abelian
normal subgroup $Q$ of $D$ such that $V = Q/Z(Q)$ is elementary
abelian of order $2^8$, where $Z(Q)$ denotes the center of $Q$.
Let $W = D/Q$. Then $W$ has an elementary abelian normal Fitting
subgroup $B$ of order $2^6$ which has a complement $L$ in $W$
isomorphic to the alternating group $A_8$. Let $t$ be a
$2$-central involution of $W$ and $D_2 = C_W(t)$. Applying now
Algorithm 7.4.8 of \cite{michler} we construct a simple subgroup
$K$ of $\GL_8(2)$ such that $|K : W| = 135$.
Furthermore, we determine a presentation of the orthogonal
simple group $K \cong O^{+}_8(2)$. It allows us to build a non
split extension $H_1$ of $K$ by $V$ in terms of generators and
relations. Then we construct all central extensions $H$ of $H_1$
by $Z(Q)$ and check which ones have a Sylow $2$-subgroup $S$ which
is isomorphic to the ones of $D$. We prove that this happens
exactly once. Thus the group $H$ has a center $Z(H) = \langle z
\rangle$. Its presentation is given in Proposition \ref{prop.
DCo_1}. Furthermore, we prove in Proposition \ref{prop. H(Co_1)} that $S$ has
a unique elementary abelian normal subgroup $A$ of order $2^{11}$
such that $N_{H}(A) \cong D$ as it is required by the Algorithm
2.5 of \cite{michler1}.

\begin{proposition}\label{prop. DCo_1} Keep the notation of Lemma \ref{l.
M24-extensions}. Let $$E = \langle a,b,c,d,t,g,h,i,j,k,v_i|
1 \le i \le 11 \rangle$$ be the split extension of $\M_{24}$ by
its simple module $V_1$ of dimension $11$ over $F = \GF(2)$. Then
the following statements hold:

\begin{enumerate}
\item[\rm(a)] $E = \langle x, y, e\rangle$, where $x =
(cgjhi)^7bt$, $y = (bv)^2jkj$ and $e = b$ have orders $7$, $4$ and
$2$, respectively.

Furthermore, $z = v_1 = (xy^3)^{14}$ is $2$-central involution of
$E$ with centralizer $D = C_{E}(z) = \langle x, y \rangle$ of
order $2^{21}\cdot3^2\cdot5\cdot7$.

\item[\rm(b)] A system of representatives $e_i$ of the $80$
conjugacy classes of $E = \langle x,y,e\rangle$ and the
corresponding centralizers orders $|C_E(e_i)|$ are given in Table
\ref{Co_1cc E}.

The character table of $E$ is Table \ref{Co_1ct_E}.

\item[\rm(c)]$D$ has a unique normal subgroup $Q$ of order $2^9$.
It is extra-special and generated by the following $8$ elements:
\begin{eqnarray*}
&&q_1 = y^2, \quad q_2 = (xy)^7,\quad q_3 = (yx)^7,\quad q_4 = (xy^2)^7,\quad q_5 = (yxy)^7, \\
&&q_6 = (x^5yx)^7, \quad q_7 = (x^4yx^2)^7, \quad q_8 = (x^4yxy)^6.\\
\end{eqnarray*}
Furthermore, the elements $q_i$ satisfy the following set
$\mathcal R(Q)$ of relations:
\begin{eqnarray*}
&&q_1^2 = q_2^2 = q_3^2 = q_4^2 = q_5^2 = q_6^4 = q_7^4 = q_8^2 = 1,\\
&&(q_1, q_2) = (q_1, q_3) = z,\quad (q_1, q_4) = (q_1, q_5) = (q_1, q_6) = 1,\\
&&(q_1, q_7) = (q_1, q_8) = z, \quad (q_2, q_3) = z, \quad (q_2, q_4) = (q_2, q_5) = 1,\\
&&(q_2, q_6) = z, \quad(q_2, q_7) = 1, \quad(q_2, q_8) = z,\\
&&(q_3, q_4) = (q_3, q_5) = (q_3, q_6) =1, \quad (q_3, q_7) = (q_3, q_8) =z,\\
&&(q_4, q_5) = (q_4, q_6) = (q_4, q_7) = 1,\quad (q_4, q_8) = z,\\
&&(q_5, q_6) = (q_5, q_7) = z,\quad (q_5, q_8) = 1,\\
&&(q_6, q_7) = (q_6, q_8) = z, \quad (q_7, q_8) = z.\\
\end{eqnarray*}

\item[\rm(d)] Let $\alpha: D \rightarrow D_1 = D/Z(Q)$ be the
canonical epimorphism with kernel $ker(\alpha) = Z(Q) = \langle z
\rangle$. Let $V = \alpha(Q) = \langle v_i = \alpha(q_i) \in D_1
\mid 1 \le i \le 8 \rangle$. Then $V$ is an elementary abelian
normal subgroup of order $2^8$ of $D_1 = \langle \alpha(x),
\alpha(y)\rangle$.

With respect to $\mathcal B = \{v_i| 1 \le i \le 8\}$ the
conjugate actions of $\alpha(x)$ and $\alpha(y)$ on $V$ have the
following matrices:

{\renewcommand{\arraystretch}{0.5} \scriptsize
$$
Mx = \left( \begin{array}{*{8}{c@{\,}}c}
 1 &   0 &   0 &   1 &   1 &   0 &   0 &   0 \\
 0 &   0 &   1 &   0 &   0 &   0 &   0 &   0 \\
 1 &   1 &   0 &   1 &   1 &   0 &   1 &   0 \\
 0 &   0 &   0 &   1 &   0 &   0 &   0 &   0 \\
 0 &   1 &   0 &   0 &   0 &   1 &   1 &   0 \\
 0 &   0 &   0 &   0 &   0 &   0 &   1 &   0 \\
 1 &   1 &   1 &   1 &   0 &   1 &   0 &   0 \\
 1 &   0 &   0 &   1 &   0 &   1 &   1 &   1 \\
\end{array} \right)\quad \text{and}~\quad
My = \left( \begin{array}{*{8}{c@{\,}}c}
 1 &   0 &   0 &   0 &   0 &   0 &   0 &   0 \\
 0 &   0 &   1 &   0 &   0 &   0 &   0 &   0 \\
 0 &   1 &   0 &   0 &   0 &   0 &   0 &   0 \\
 0 &   0 &   0 &   0 &   1 &   0 &   0 &   0 \\
 0 &   0 &   0 &   1 &   0 &   0 &   0 &   0 \\
 0 &   0 &   1 &   0 &   1 &   0 &   0 &   1 \\
 1 &   1 &   0 &   1 &   0 &   1 &   1 &   1 \\
 0 &   1 &   0 &   1 &   0 &   1 &   0 &   0 \\
\end{array} \right)
$$
} in $\GL_8(2)$.

\item[\rm(e)] The map $\varphi : D \rightarrow \GL_8(2)$ defined
by $\varphi(x) = Mx$ and  $\varphi(y) = My$ is a group epimorphism
onto $MD = \langle Mx, My \rangle$ with kernel $Q$. The group $MD$
has an elementary abelian Fitting subgroup $MY$ of order $2^6$
generated by $\varphi(m_j), 1 \le j \le 6 $, where
\begin{eqnarray*}
&&m_1 = (x^4(yx)^2y)^5, \quad m_2 = (x^3(yx)^3)^5,\quad m_3 = (x^2(yx)^2yx^2)^5,\\
&&m_4 = (xyx^4yxy)^5,\quad m_5 = (yx^4yxyx)^5, \quad m_6 = (x^3yx^3yxy)^5.\\
\end{eqnarray*}

Furthermore, $MC = \langle \varphi(x), \varphi(c) \rangle$ and $C
= \langle x, c \rangle$ are complements of $MY$ in $MD$ and of $M
= \langle Q, m_j \mid 1 \le j \le 6 \rangle$ in $D$, respectively,
which both are isomorphic to the alternating group $A_8$. This has been checked with the isomorphism testing program of Cannon and Holt \cite{cannon} implemented in MAGMA. In
particular, $D$ has a faithful permutation representation $PD$ of
degree $2^{15}$ with stabilizer $C$, where $c =
(x^2yx^2yxyxy^2xy)^3$ has order $2$.

\item[\rm(f)] $D_1$ has a faithful permutation representation
$PD_1$ of degree $2^{14}$ with stabilizer $\alpha(C) = \langle
\alpha(x), \alpha(c) \rangle$.

Furthermore, $D_1$ is a non split extension of $D/Q$ by $V$ and
$C_{D_1}(V) = V$.

\item[\rm(g)] $D$ is a finitely presented group $D = \langle Q, x,
y \rangle$ having a set ${\mathcal R}(D)$ of defining relations
consisting of ${\mathcal R}(Q)$ and the following relations:
\begin{eqnarray*}
&&x^7 = y^4 = 1, \quad z = (xy^3)^{14}, \quad z^2 = 1, \quad (x, z) = (y, z) = 1,\\
&&(yx^{-1})^7 = q_1q_2q_4,\quad (xyx^{-1}yx)^4 = q_3q_6q_7,\\
&&(xyx^{-1}yx^{-1}y)^4 = q_4q_6q_7q_8,\\
&&(xyx^{-1}y)^6 = q_1,\quad (xyx^{-1}yx^{-1}yx)^4 = q_5z,\\
&&x^{-1}yx^{-2}yx^2yxyx^{-1}yx^3yx^{-1}yxyx^2yx^{-2}yx^{-1}y = q_2q_7z,\\
&&x^{-2}yx^{-1}yxyx^{-2}yx^{-1}yx^{-3}yx^2yx^2yxyx^{-1}yx^{-1}yx^{-1} = q_2q_4q_5q_8,\\
&&yxyx^{-1}yx^{-2}yx^2yx^{-1}yx^{-2}yxyx^{-2}yx^2yx^{-2}yx^2 = q_5q_6q_8,\\
&&x^{-1}q_1x = q_1q_4q_5q_6^2,\quad x^{-1}q_2x = q_3,\quad x^{-1}q_3x = q_1q_2q_4q_5q_7,\\
&&x^{-1}q_4x = q_4,\quad x^{-1}q_5x = q_2q_6q_7,\quad x^{-1}q_6x = q_7,\\
&&x^{-1}q_7x = q_1q_3q_2q_4q_6,\quad x^{-1}q_8x = q_1q_4q_6q_8q_7,\\
&&q_1^{y} = q_1,\quad y^{-1}q_2yq_3 = (q_1q_2)^2,\quad y^{-1}q_3yq_2 = 1,\\
&&y^{-1}q_4yq_5 = 1,\quad y^{-1}q_5yq_4 = 1,\quad y^{-1}q_6yq_3q_5q_8 = 1,\\
&&y^{-1}q_7yq_1q_2q_4q_6q_7q_8 = 1,\quad y^{-1}q_8yq_2q_4q_6 = (q_1q_2)^2.\\
\end{eqnarray*}

\item[\rm(h)] $MY$ contains a $2$-central involution $Mt =
\varphi[(x^3yx^3yx^2)^3]$ of $MD$ with centralizer $C_{MD}(Mt)$ of
order $2^{12}\cdot 3^2$. Furthermore, $MY$ has a complement $MC_1$
in $C_{MD}(Mt)$ having an elementary abelian Fitting subgroup
$MW_1 = \langle \varphi(f_i) \mid 1\le i \le 4 \rangle$ of order
$2^4$ which again has a complement $MC_2 = \langle \varphi(r_k),
\varphi(t_k) \mid 1 \le k \le 2 \rangle$ of order $2^2\cdot 3^2$,
where
\begin{eqnarray*}
&&f_1 = (yx^3yxyx^2yxyxy^2)^2,\quad f_2 = (x^2yxyx^2yxyxyx^5y)^3,\\
&&f_3 = (xyx^2yx^2yxyx^5yx^2)^2,\quad f_4 = (x^2yx^4yx^3yx^3y^2xy)^3,\\
&& r_1 = (yx^5yx^5yx^2yx^2)^2,\quad r_2 = (x^3yxyx^2yx^4yxyx^2y)^3,\\
&&t_1 = (xyx^3yx^2y^2x^2yxyx^3yx)^5(x^5yx^3yxyx^2yx^3yxy)^3\\
&& \quad \cdot(xyx^3yx^2y^2x^2yxyx^3yx)^{10}(x^5yx^3yxyx^2yx^3yxy)^3,\\
&&t_2 = ((xyx^3yx^2y^2x^2yxyx^3yx)^5(x^5yx^3yxyx^2yx^3yxy)^3)^3,\\
\end{eqnarray*}
and where $r_1$ and $t_1$ have order $3$ and all other elements
are involutions. Furthermore, $r_1^{r_2} = r_1^2$ and $t_1^{t_2} =
t_1^2$.

\item[\rm(i)] The centralizer $MX = C_{\GL_{8}(2)}(Mt)$ of $Mt$ in
$\GL_8(2)$ contains the matrix

{\renewcommand{\arraystretch}{0.5}

\scriptsize

$$
Mh = \left( \begin{array}{*{8}{c@{\,}}c}
 1 &   1 &   0 &   1 &   1 &   1 &   0 &   1 \\
 1 &   0 &   0 &   1 &   1 &   1 &   0 &   0 \\
 0 &   0 &   1 &   0 &   0 &   0 &   0 &   0 \\
 1 &   1 &   0 &   1 &   1 &   1 &   0 &   0 \\
 1 &   0 &   0 &   1 &   1 &   0 &   0 &   1 \\
 1 &   1 &   0 &   0 &   1 &   0 &   0 &   0 \\
 0 &   1 &   0 &   0 &   1 &   1 &   1 &   1 \\
 0 &   1 &   0 &   1 &   0 &   1 &   1 &   1 \\
\end{array} \right).
$$
}

of order $3$ such that the subgroup $MO = \langle MD, Mh \rangle$
of $\GL_8(2)$ has order $2^{12}\cdot3^5\cdot5^2\cdot7$ and the
subgroup
$$MH_3 = C_{MO}(Mt) = \langle C_{MD}(Mt), Mh \rangle = \langle Mm_j, Mf_i, Mr_k,Mt_k,Mh \rangle$$

of $\GL_8(2)$ has order $2^{12}\cdot3^3$, where $1 \le j \le 6, 1
\le i \le 4, 1 \le k \le 2$.

Furthermore, $MO$ is isomorphic to the simple group $O^{+}_8(2)$.

\item[\rm(j)] The Fitting subgroup $MY_3$ of $MH_3$ is
extra-special of order $2^9$ with center $Z(Y_3) = Mt$, and
$$MY_3 = \langle Mf_1, Mf_2, Mf_3, Mf_4, Mf_5, Mf_6, Mf_7, Mf_8\rangle,$$
where $Mf_5 = Mm_3$, $Mf_6 = Mm_4$, $Mf_7 = Mm_6$ and $Mf_8 = (Mm_5Mf_2)^2$.
Furthermore, $MY_3$ is isomorphic to the finitely presented
group $Y_3 = \langle f_i \mid 1 \le i \le 8\rangle$ with set $\mathcal
R(Y_3)$ of defining relations:
\begin{eqnarray*}
&&f_1^2 =  f_2^2 = f_3^2 =  f_4^2 = f_5^2 = f_6^2 =  f_7^2 =  f_8^2 = 1,\\
&&(f_1,  f_2) = (f_1,  f_3) = (f_2,  f_3) = (f_1,  f_4) = (f_2,  f_4) = (f_3,  f_4) = (f_1,  f_5) = 1,\\
&&(f_2,  f_5) = (f_3,  f_5) = (f_1,  f_6) = (f_5,  f_6) = (f_5,  f_7) = (f_6,  f_7) = (f_2,  f_8) = 1,\\
&&(f_4,  f_8) = (f_5,  f_8) = (f_6,  f_8) = (f_7,  f_8) = 1,\quad  (f_1  f_7  f_2)^2 = (f_2  f_6  f_3)^2 = 1,\\
&&(f_1  f_7  f_3)^2 = (f_1  f_8  f_3)^2 = (f_2  f_6  f_4)^2 = (f_1  f_7  f_4)^2 = 1,\\
&& f_4  f_6  f_5  f_4  f_5  f_6 = f_4  f_7  f_5  f_4  f_5  f_7 =  f_3  f_8  f_6  f_3  f_6  f_8 = 1.\\
\end{eqnarray*}

$MT = \langle Mr_1, Mt_1, Ma_1\rangle$ is an elementary abelian
Sylow $3$-subgroup of $MH_3$ of order $27$, where $Ma_1 =
(Mt_1)^2Mh(Mt_1)^2Mh^2(Mt_1)^2Mh$.

The matrix $Ma_2 = Mf_2Mm_2Mm_3Mm_4Mm_6Mf_2Mm_1Mm_2Mm_6$ is an
involution commuting with the subgroup $\langle
Mr_1,Mr_2,Mt_1,Mt_2\rangle$ and satisfying $Ma_1^{Ma_2} = Ma_1^2$.

$MK_3 = \langle MT, Mr_2, Mt_2, Ma_1, Ma_2\rangle$ is a complement
of $MY_3$ in $MH_3$.

$MH_3$ is isomorphic to the finitely presented group
$$H_3 = \langle a_1,a_2,r_1,r_2,t_1,t_2,f_i \mid 1 \le i \le 8\rangle$$
with the set $\mathcal R(H_3)$ of defining relations consisting of
$\mathcal R(Y_3)$ and the following relations:
\begin{eqnarray*}
&&a_1^3 = a_2^2 = r_1^3 = r_2^2 = t_1^3 = t_2^2 = 1,\\
&&(a_i,r_j) = (a_i,t_k) = (r_j,t_k) = 1 \quad \mbox{for all} \quad 1 \le i,j,k \le 2,\\
&&a_1^{a_2} = a_1^2,\quad r_1^{r_2} = r_1^2,\quad t_1^{t_2} = t_1^2,\\
&&f_1^{r_1} = f_2f_3f_4,\quad f_1^{r_2} = f_2f_3f_4,\quad f_1^{t_1} = f_1f_4,\quad f_1^{t_2} = f_1f_4,\\
&&f_1^{a_1} = f_5f_6,\quad f_1^{a_2} = f_1f_5f_6,\quad f_2^{r_1} = f_4,\quad f_2^{r_2} = f_2,\\
&&f_2^{t_1} = f_1f_2f_3f_4,\quad f_2^{t_2} = f_2,\quad f_2^{a_1} = f_8,\quad f_2^{a_2} = f_2f_8,\\
&&f_3^{r_1} = f_1f_3f_4,\quad f_3^{r_2} = f_1f_4,\quad f_3^{t_1} = f_2f_4,\quad f_3^{t_2} = f_2f_3f_4,\\
&&f_3^{a_1} = f_4f_6f_4f_8,\quad f_3^{a_2} = f_6f_3f_8,\quad f_4^{r_1} = f_2f_4,\quad f_4^{r_2} = f_2f_4,\\
&&f_4^{t_1} = f_1,\quad f_4^{t_2} = f_4,\quad f_4^{a_1} = f_6f_7f_8,\quad f_4^{a_2} = f_4f_6f_7f_8,\\
&&f_5^{r_1} = f_4f_5f_4f_7,\quad f_5^{r_2} = f_5f_8,\quad f_5^{t_1} = f_7,\quad f_5^{t_2} = f_5f_8,\\
&&f_5^{a_1} = f_1f_2f_5f_8f_3f_8,\quad f_5^{a_2} = f_5,\quad f_6^{r_1} = f_5f_8,\\
&&f_6^{r_2} = f_4f_5f_4f_7,\quad f_6^{t_1} = f_5f_8,\quad f_6^{t_2} = f_7,\quad f_6^{a_1} = f_2f_6f_3,\\
&&f_6^{a_2} = f_6,\quad f_7^{r_1} = f_5,\quad f_7^{r_2} = f_4f_5f_4f_6f_8, \quad f_7^{t_1} = f_4f_5f_4f_7,\\
&&f_7^{t_2} = f_6,\quad f_7^{a_1} = f_3f_7f_4,\quad f_7^{a_2} = f_7,\quad f_8^{r_1} = f_6f_7f_8,\\
&&f_8^{r_2} = f_8,\quad f_8^{t_1} = f_4f_5f_4f_6f_7f_8,\quad f_8^{t_2} = f_8,\quad f_8^{a_1} = f_2f_8,\\
&&f_8^{a_2} = f_8.
\end{eqnarray*}

\item[\rm(k)] The subgroup
$$U = \langle a_2,r_1,r_2,t_1,t_2,f_j, q_i \mid 1 \le i, j \le 8 \rangle$$
of $D$ has center $Z(U) = \langle z = q_6^2 \rangle$ and satisfies
the set $\mathcal R(U)$ of defining relations consisting of
$\mathcal R(Q)$ and the following relations:
\begin{eqnarray*}
&&a_2^2 = r_1^3 = r_2^2 = t_1^3 = t_2^2 = 1,\\
&&f_j^2 = 1,\quad \mbox{for} \quad 1 \le j \le 8,\\
&&r_1^{r_2} = r_1^2,\quad t_1^{t_2} = t_1^2,\\
&&(a_2, r_1) = (a_2, r_2) = (a_2, t_1) = (a_2, t_2) = 1,\\
&&(r_1, t_1) = (r_1, t_2) = (r_2, t_1) = (r_2, t_2) = 1,\\
&&f_1^{r_1} ( f_2f_3f_4)^{-1}  = f_1^{r_2} ( f_2f_3f_4)^{-1}  = f_1^{t_1} ( f_1f_4)^{-1}  = f_1^{t_2} ( f_1f_4)^{-1}  = 1,\\
&&f_1^{a_2} ( f_1f_5f_6)^{-1}  = q_5,\\
&&f_2^{r_1} ( f_4)^{-1}  = f_2^{r_2} ( f_2)^{-1}  = f_2^{t_1} ( f_1f_2f_3f_4)^{-1}  = f_2^{t_2} ( f_2)^{-1}  = 1,\\
&&f_2^{a_2} ( f_2f_8)^{-1}  = f_3^{r_1} ( f_1f_3f_4)^{-1}  = f_3^{r_2} ( f_1f_4)^{-1}  = f_3^{t_1} ( f_2f_4)^{-1}  = 1,\\
&&f_3^{t_2} ( f_2f_3f_4)^{-1}  = 1,\quad f_3^{a_2} ( f_6f_3f_8)^{-1}  = q_4,\\
&&f_4^{r_1} ( f_2f_4)^{-1}  = f_4^{r_2} ( f_2f_4)^{-1}  = f_4^{t_1} ( f_1)^{-1}  = f_4^{t_2} ( f_4)^{-1}  = 1,\\
&&f_4^{a_2} ( f_4f_6f_7f_8)^{-1}  = q_1q_2q_4q_6q_7,\quad f_5^{r_1} ( f_4f_5f_4f_7)^{-1}  = q_2q_4q_7q_6,\\
&&f_5^{r_2} ( f_5f_8)^{-1}  = q_2q_4q_6q_7,\quad f_5^{t_1}(f_7)^{-1}  = q_8q_1q_8,\\
&&f_5^{t_2} ( f_5f_8)^{-1}  = q_2q_4q_5q_6q_7,\quad f_5^{a_2}  ( f_5)^{-1}  = 1,\\
\end{eqnarray*}
\begin{eqnarray*}
&&f_6^{r_1}  ( f_5f_8)^{-1}  = q_2q_5q_6q_7,\quad f_6^{r_2}  ( f_4f_5f_4f_7)^{-1}  = q_2q_4q_5q_7q_6,\\
&&f_6^{t_1}  ( f_5f_8)^{-1}  = q_2q_4q_6q_7,\quad f_6^{t_2}  ( f_7)^{-1}  = q_5q_8q_1q_8,\\
&&f_6^{a_2}  ( f_6)^{-1}  = 1,\quad f_7^{r_1}  ( f_5)^{-1}  = q_8q_1q_8,\\
&&f_7^{r_2}  ( f_4f_5f_4f_6f_8)^{-1}  = q_5,\quad f_7^{t_1}  ( f_4f_5f_4f_7)^{-1}  = q_1,\\
&&f_7^{t_2}  ( f_6)^{-1}  = q_1q_2q_4q_7q_6,\quad f_7^{a_2}  ( f_7)^{-1}  = 1,\\
&&f_8^{r_1}  ( f_6f_7f_8)^{-1}  = q_1q_5,\quad f_8^{r_2}  ( f_8)^{-1}  = 1,\\
&&f_8^{t_1}  ( f_4f_5f_4f_6f_7f_8)^{-1}  = q_1q_4,\quad f_8^{t_2}(f_8)^{-1}  = 1,\quad f_8^{a_2}(f_8)^{-1}  = 1,\\
&&(f_1,  f_2)  = (f_1,  f_3)  = (f_2,  f_3)  = (f_1,  f_4)  = (f_2,  f_4)  = (f_3,  f_4)  = 1,\\
&&(f_2,  f_5)  = (f_3,  f_5)  = (f_1,  f_6)  = (f_5,  f_6)  = (f_5,  f_7)  = (f_6,  f_7)  = 1,\\
&&(f_2,  f_8)  = 1,\quad (f_1 f_5)^2  = q_1q_4q_5,\quad (f_4f_8)^2  = q_2q_4q_5q_7q_6,\\
&&(f_5,  f_8)  = (f_6,  f_8)  = (f_7,  f_8)  = 1,\\
&&(f_1  f_7  f_2)^2  = q_1q_4q_5,\quad (f_2  f_6  f_3)^2  = q_2q_4q_5q_6q_7,\\
&&(f_1  f_7  f_3)^2  = q_1q_2q_6q_7,\quad (f_1  f_8  f_3)^2  = q_2q_4q_5q_6q_7,\\
&&(f_2  f_6  f_4)^2  = q_2q_4q_5q_6q_7,\quad (f_1  f_7  f_4)^2  = q_1q_2q_6q_7,\\
&&f_4  f_6  f_5  f_4  f_5  f_6  = q_2q_4q_5q_7q_6,\quad f_4  f_7  f_5  f_4  f_5  f_7  = q_2q_4q_5q_6q_7,\\
&&f_3  f_8  f_6  f_3  f_6  f_8 = 1,\\
&&q_1^{a_2} = q_1,\quad q_1^{r_1} = q_1q_5,\quad q_1^{r_2} = q_1q_5,\\
&&q_1^{t_1} = q_1q_2q_4q_5q_6q_7,\quad q_1^{t_2} = q_1q_2q_4q_5q_6q_7,\quad q_1^{f_1} = q_1,\\
&&q_1^{f_2} = q_1q_2q_4q_5q_6q_7,\quad q_1^{f_3} = q_2q_6q_7,\quad q_1^{f_4} = q_1,\\
&&q_1^{f_5} = q_1^{f_6} = q_1^{f_7} = q_1^{f_8} = q_1,\\
&&q_2^{a_2} = q_1q_6q_7,\quad q_2^{r_1} = q_2q_5,\quad q_2^{r_2} = q_2q_5,\\
&&q_2^{t_1} = q_2q_3q_4q_5q_6q_7,\quad q_2^{t_2} = q_3q_7q_6,\quad q_2^{f_1} = q_4q_5q_8q_6,\\
&&q_2^{f_2} = q_4q_6q_7q_8,\quad q_2^{f_3} = q_1q_7q_6,\quad q_2^{f_4} = q_4q_5q_8q_6,\\
&&q_2^{f_5} = q_1q_4q_6q_7,\quad q_2^{f_6} = q_1q_5q_6q_7,\quad q_2^{f_7} = q_1q_4q_5q_7q_6,\\
&&q_2^{f_8} = q_5q_8q_2q_8,\quad q_3^{a_2} = q_2q_3q_4q_5q_6q_7,\quad q_3^{r_1} = q_3q_5,\\
&&q_3^{r_2} = q_3q_5,\quad q_3^{t_1} = q_1q_3q_4q_5q_6q_7,\quad q_3^{t_2} = q_3,\\
&&q_3^{f_1} = q_2q_3q_4q_5q_6q_8,\quad q_3^{f_2} = q_2q_3q_4q_5q_7q_6,\\
&&q_3^{f_3} = q_1q_3q_4q_7,\quad q_3^{f_4} = q_3,\quad q_3^{f_5} = q_2q_3q_4q_6q_7,\\
&&q_3^{f_6} = q_2q_3q_5q_6q_7,\quad q_3^{f_7} = q_2q_3q_5q_7q_6,\quad q_3^{f_8} = q_3,\\
&&q_4^{a_2} = q_4,\quad q_4^{r_1} = q_5,\quad q_4^{r_2} = q_4q_5,\quad q_4^{t_1} = q_4,\\
&&q_4^{t_2} = q_4,\quad q_4^{f_1} = q_1q_5,\quad q_4^{f_2} = q_2q_5q_6q_7,\quad q_4^{f_3} = q_4,\\
&&q_4^{f_4} = q_2q_5q_6q_7,\quad q_4^{f_5} = q_4^{f_6} = q_4^{f_7} = q_4^{f_8} = q_4,\\
&&q_5^{a_2} = q_5,\quad q_5^{r_1} = q_4q_5,\quad q_5^{r_2} = q_5^{t_1} = q_5^{t_2} = q_5,\\
&&q_5^{f_1} = q_1q_4,\quad q_5^{f_2} = q_5,\quad q_5^{f_3} = q_1q_2q_5q_6q_7,\\
&&q_5^{f_4} = q_2q_4q_6q_7,\quad q_5^{f_5} = q_5^{f_6} = q_5^{f_7} = q_5^{f_8} = q_5,\\
&&q_6^{a_2} = q_1q_2q_4q_7,\quad q_6^{r_1} = q_2q_5q_8,\quad q_6^{r_2} = q_5q_6q_8,\\
&&q_6^{t_1} = q_3q_4q_5q_7,\quad q_6^{t_2} = q_2q_7q_3,\quad q_6^{f_1} = q_1q_2q_8,\\
&&q_6^{f_2} = q_5q_6q_8,\quad q_6^{f_3} = q_1q_2q_7,\quad q_6^{f_4} = q_6q_8q_7,\quad q_6^{f_5}=q_4q_6,\\
&&q_6^{f_6} = q_1q_4q_6,\quad q_6^{f_7} = q_1q_6^3,\quad q_6^{f_8} = q_2q_4q_7,\\
&&q_7^{a_2} = q_4q_7^3,\quad q_7^{r_1} = q_2q_5q_6q_8q_7,\quad q_7^{r_2} = q_5q_7q_8,\\
&&q_7^{t_1} = q_1q_6q_2,\quad q_7^{t_2} = q_7,\quad q_7^{f_1} = q_1q_4q_5q_7,\quad q_7^{f_2} = q_7,\\
&&q_7^{f_3} = q_1q_6q_2,\quad q_7^{f_4} = q_2q_4q_6q_5,\quad q_7^{f_5} = q_1q_2q_6,\\
\end{eqnarray*}
\begin{eqnarray*}
&&q_7^{f_6} = q_2q_4q_6q_5,\quad q_7^{f_7} = q_2q_4q_6q_5,\quad q_7^{f_8} = q_2q_4q_6q_5,\\
&&q_8^{a_2} = q_5q_7^2q_8,\quad q_8^{r_1} = q_2q_6q_8,\quad q_8^{r_2} = q_8^{t_1} = q_8^{t_2} = q_8,\\
&&q_8^{f_1} = q_1q_4q_5q_8,\quad q_8^{f_2} = q_8,\quad q_8^{f_3} = q_1q_2q_6q_7q_8,\\
&&q_8^{f_4} = q_2q_4q_5q_6q_7q_8,\quad q_8^{f_5} = q_2q_4q_5q_6q_8q_7,\\
&&q_8^{f_6} = q_1q_2q_6q_8q_7,\quad q_8^{f_7} = q_1q_4q_5q_8,\quad q_8^{f_8} = q_8.\\
\end{eqnarray*}
$U$ has center $Z(U) = \langle q_6^2 \rangle$ and a faithful
permutation representation $PU$ of degree $2048$ with stabilizer
$Y$ generated by $r_2r_1$, $(t_2q_6)^2$, $(q_6r_1f_4)^4$,
$(q_6f_8r_1)^4$, $(r_1^2q_6f_8)^3$, $(r_1t_1t_2f_4)^3$,
$(r_1t_1q_6f_8)^2$.

\item[\rm(l)] The finitely presented group $H_3$ has a non split
extension
$$H_2 = \langle a_1, a_2, r_1, r_2, t_1, t_2, f_j, v_i \mid 1 \le
i, j \le 8\rangle $$ by the elementary abelian group $V = \langle
v_i \mid 1 \le i \le 8\rangle$ with the following set $\mathcal
R(H_2)$ of defining relations:
\begin{eqnarray*}
&&a_1^3 = a_2^2 = r_1^3 = r_2^2 = t_1^3 = t_2^2 = 1,\\
&&a_1^{a_2} = a_1^2,\quad r_1^{r_2} = r_1^2,\quad t_1^{t_2} = t_1^2,\\
&&(a_k, r_1) = (a_k, r_2) = (a_k, t_1) = (a_k, t_2) = 1,\quad \mbox {for}\quad 1 \le k \le 2,\\
&&(r_1, t_1) = (r_1, t_2) = (r_2, t_1) = (r_2, t_2) = 1,\\
&&v_i^2 = 1 \quad \mbox{for} \quad 1 \le i \le 8,\\
&&f_j^2 = 1 \quad \mbox{for} \quad 1 \le j \le 8,\\
&&(v_i, v_j) = 1\quad \mbox{for} \quad 1 \le i, j \le 8,\\
&&a_1^{{-1}} v_1 a_1 v_1^{-1} v_2^{-1}v_4^{-1} v_5^{-1} v_6^{-1} v_8^{-1} = a_1^{-1} v_2 a_1 v_4^{-1} v_5^{-1}v_6^{-1} v_8^{-1} = 1,\\
&&a_1^{-1} v_3 a_1 v_2^{-1} v_3^{-1}v_4^{-1} v_5^{-1} v_6^{-1} v_8^{-1} = a_1^{-1} v_4 a_1 v_2^{-1} v_6^{-1} = 1,\\
&&a_1^{-1} v_5 a_1 v_8^{-1} = a_1^{-1} v_6 a_1 v_2^{-1} v_5^{-1}v_8^{-1} = a_1^{-1} v_8 a_1 v_5^{-1} v_8^{-1} = 1,\\
&&a_1^{-1} v_7 a_1 v_2^{-1} v_5^{-1}v_6^{-1} v_7^{-1} v_8^{-1} = 1,\\
&&(a_2, v_1^{-1}) = (a_2, v_2^{-1}) = (a_2, v_3^{-1}) = (a_2, v_4^{-1}) = (a_2, v_5^{-1}) = 1,\\
&&a_2^{-1} v_6 a_2 v_4^{-1} v_6^{-1} = a_2^{-1} v_7 a_2 v_4^{-1} v_7^{-1} = 1,\\
&&a_2^{-1} v_8 a_2 v_5^{-1} v_8^{-1} = r_1^{-1} v_1 r_1 v_1^{-1} v_5^{-1} = r_1^{-1} v_2 r_1 v_2^{-1} v_5^{-1} = 1,\\
&&r_1^{-1} v_3 r_1 v_3^{-1} v_5^{-1} = r_1^{-1} v_4 r_1 v_5^{-1} = r_1^{-1} v_5 r_1 v_4^{-1} v_5^{-1} = 1,\\
&&r_1^{-1} v_6 r_1 v_2^{-1} v_5^{-1}v_8^{-1} = r_1^{-1} v_7 r_1 v_2^{-1} v_5^{-1}v_6^{-1} v_7^{-1} v_8^{-1} = 1,\\
&&r_1^{-1} v_8 r_1 v_2^{-1} v_6^{-1}v_8^{-1} = r_2^{-1} v_1 r_2 v_1^{-1} v_5^{-1} = 1,\\
&&r_2^{-1} v_2 r_2 v_2^{-1} v_5^{-1} = r_2^{-1} v_3 r_2 v_3^{-1} v_5^{-1} = 1,\\
&&r_2^{-1} v_4 r_2 v_4^{-1} v_5^{-1} = 1,\quad (r_2, v_5^{-1}) = 1,\\
&&r_2^{-1} v_6 r_2 v_5^{-1} v_6^{-1}v_8^{-1} = 1,\\
&&r_2^{-1} v_7 r_2 v_5^{-1} v_7^{-1}v_8^{-1} = 1,\quad (r_2, v_8^{-1}) = 1,\\
&&t_1^{-1} v_1 t_1 v_1^{-1} v_2^{-1}v_4^{-1} v_5^{-1} v_6^{-1} v_7^{-1} = 1,\\
&&t_1^{-1} v_2 t_1 v_2^{-1} v_3^{-1}v_4^{-1} v_5^{-1} v_6^{-1} v_7^{-1} = 1,\\
&&t_1^{-1} v_3 t_1 v_1^{-1} v_3^{-1}v_4^{-1} v_5^{-1} v_6^{-1} v_7^{-1} = 1,\\
&&(t_1, v_4^{-1}) = (t_1, v_5^{-1}) = t_1^{-1} v_6 t_1 v_3^{-1} v_4^{-1}v_5^{-1} v_7^{-1} = 1,\\
&&t_1^{-1} v_7 t_1 v_1^{-1} v_2^{-1}v_6^{-1} = 1,\quad (t_1, v_8^{-1}) = 1,\\
\end{eqnarray*}
\begin{eqnarray*}
&&t_2^{-1} v_1 t_2 v_1^{-1} v_2^{-1}v_4^{-1} v_5^{-1} v_6^{-1} v_7^{-1} = 1,\\
&&t_2^{-1} v_2 t_2 v_3^{-1} v_6^{-1}v_7^{-1} = 1,\quad (t_2, v_3^{-1}) = (t_2, v_4^{-1}) = (t_2, v_5^{-1}) = 1,\\
&&t_2^{-1} v_6 t_2 v_2^{-1} v_3^{-1}v_7^{-1} = 1,\quad (t_2, v_7^{-1}) = (t_2, v_8^{-1}) = (f_1, v_1^{-1}) = 1,\\
&&f_1^{-1} v_2 f_1 v_4^{-1} v_5^{-1}v_6^{-1} v_8^{-1} = f_1^{-1} v_3 f_1 v_2^{-1} v_3^{-1}v_4^{-1} v_5^{-1} v_6^{-1} v_8^{-1} = 1,\\
&&f_1^{-1} v_4 f_1 v_1^{-1} v_5^{-1} = f_1^{-1} v_5 f_1 v_1^{-1} v_4^{-1} = f_1^{-1} v_6 f_1 v_1^{-1} v_2^{-1}v_8^{-1} = 1,\\
&&f_1^{-1} v_7 f_1 v_1^{-1} v_4^{-1}v_5^{-1} v_7^{-1} = f_1^{-1} v_8 f_1 v_1^{-1} v_4^{-1}v_5^{-1} v_8^{-1} = 1,\\
&&f_2^{-1} v_1 f_2 v_1^{-1} v_2^{-1}v_4^{-1} v_5^{-1} v_6^{-1} v_7^{-1} = 1,\\
&&f_2^{-1} v_2 f_2 v_4^{-1} v_6^{-1}v_7^{-1} v_8^{-1} = f_2^{-1} v_3 f_2 v_2^{-1} v_3^{-1}v_4^{-1} v_5^{-1} v_6^{-1} v_7^{-1} = 1,\\
&&f_2^{-1} v_4 f_2 v_2^{-1} v_5^{-1}v_6^{-1} v_7^{-1} = f_2^{-1} v_6 f_2 v_5^{-1} v_6^{-1}v_8^{-1} = 1,\\
&&(f_2, v_5) = (f_2, v_7) = (f_2, v_8) = (f_3, v_4) = (f_4, v_1) = (f_4, v_3) = 1,\\
&&(f_5, v_1) = (f_5, v_4) = (f_5, v_5) = (f_6, v_1) = (f_6, v_4) = (f_6, v_5) = 1,\\
&&(f_7, v_1) = (f_8, v_1) = (f_7, v_4) = (f_7, v_5) = (f_8, v_3) = (f_8, v_4) = 1,\\
&&(f_8, v_5) = (f_8, v_8) = 1,\quad f_3^{-1} v_1 f_3 v_2^{-1} v_6^{-1}v_7^{-1} = 1,\\
&&f_3^{-1} v_2 f_3 v_1^{-1} v_6^{-1}v_7^{-1} = f_3^{-1} v_3 f_3 v_1^{-1} v_3^{-1}v_4^{-1} v_7^{-1} = 1,\\
&&f_3^{-1} v_5 f_3 v_1^{-1} v_2^{-1}v_5^{-1} v_6^{-1} v_7^{-1} = f_3^{-1} v_6 f_3 v_1^{-1} v_2^{-1}v_7^{-1} = 1,\\
&&f_3^{-1} v_7 f_3 v_1^{-1} v_2^{-1}v_6^{-1} = f_3^{-1} v_8 f_3 v_1^{-1} v_2^{-1}v_6^{-1} v_7^{-1} v_8^{-1} = 1,\\
&&f_4^{-1} v_2 f_4 v_4^{-1} v_5^{-1}v_6^{-1} v_8^{-1} = f_4^{-1} v_4 f_4 v_2^{-1} v_5^{-1}v_6^{-1} v_7^{-1} = 1,\\
&&f_4^{-1} v_5 f_4 v_2^{-1} v_4^{-1}v_6^{-1} v_7^{-1} = f_4^{-1} v_6 f_4 v_6^{-1} v_7^{-1}v_8^{-1} = 1,\\
&&f_4^{-1} v_7 f_4 v_2^{-1} v_4^{-1}v_5^{-1} v_6^{-1} = f_4^{-1} v_8 f_4 v_2^{-1} v_4^{-1}v_5^{-1} v_6^{-1} v_7^{-1} v_8^{-1} = 1,\\
&&f_5^{-1} v_2 f_5 v_1^{-1} v_4^{-1}v_6^{-1} v_7^{-1} = f_5^{-1} v_3 f_5 v_2^{-1} v_3^{-1}v_4^{-1} v_6^{-1} v_7^{-1} = 1,\\
&&f_5^{-1} v_6 f_5 v_4^{-1} v_6^{-1} = f_5^{-1} v_7 f_5 v_1^{-1} v_2^{-1}v_6^{-1} = 1,\\
&&f_5^{-1} v_8 f_5 v_2^{-1} v_4^{-1}v_5^{-1} v_6^{-1} v_7^{-1} v_8^{-1} = 1,\\
&&f_6^{-1} v_2 f_6 v_1^{-1} v_5^{-1}v_6^{-1} v_7^{-1} = f_6^{-1} v_3 f_6 v_2^{-1} v_3^{-1}v_5^{-1} v_6^{-1} v_7^{-1} = 1,\\
&&f_6^{-1} v_6 f_6 v_1^{-1} v_4^{-1}v_6^{-1} = f_6^{-1} v_7 f_6 v_2^{-1} v_4^{-1}v_5^{-1} v_6^{-1} = 1,\\
&&f_6^{-1} v_8 f_6 v_1^{-1} v_2^{-1}v_6^{-1} v_7^{-1} v_8^{-1} = f_7^{-1} v_2 f_7 v_1^{-1} v_4^{-1}v_5^{-1} v_6^{-1} v_7^{-1} = 1,\\
&&f_7^{-1} v_3 f_7 v_2^{-1} v_3^{-1}v_5^{-1} v_6^{-1} v_7^{-1} = 1,\\
&&f_7^{-1} v_6 f_7 v_1^{-1} v_6^{-1} = f_7^{-1} v_7 f_7 v_2^{-1} v_4^{-1}v_5^{-1} v_6^{-1} = 1,\\
&&f_7^{-1} v_8 f_7 v_1^{-1} v_4^{-1}v_5^{-1} v_8^{-1} = f_8^{-1} v_2 f_8 v_2^{-1} v_5^{-1} = 1,\\
&&f_8^{-1} v_6 f_8 v_2^{-1} v_4^{-1}v_7^{-1} = f_8^{-1} v_7 f_8 v_2^{-1} v_4^{-1}v_5^{-1} v_6^{-1} = 1,\\
&&r_1^{-1}f_1r_1f_4^{-1}f_3^{-1}f_2^{-1}v_2^{-1}v_4^{-1}v_5^{-1}v_6^{-1}v_7^{-1} = 1,\\
&&r_2^{-1}f_1r_2f_4^{-1}f_3^{-1}f_2^{-1}v_2^{-1}v_4^{-1}v_5^{-1}v_6^{-1}v_7^{-1} = 1,\\
&&t_1^{-1} f_1 t_1 f_4^{-1} f_1^{-1} = t_2^{-1} f_1 t_2 f_4^{-1} f_1^{-1} = 1,\\
&&a_1^{-1} f_1 a_1 f_6^{-1} f_5^{-1}v_1^{-1} v_4^{-1} v_5^{-1} = 1,\\
&&a_2^{-1} f_1 a_2 f_6^{-1} f_5^{-1}f_1^{-1} v_1^{-1} v_4^{-1} v_5^{-1} = 1,\\
&&r_1^{-1} f_2 r_1 f_4^{-1} = 1,\quad (r_2, f_2^{-1}) = (t_2, f_2^{-1}) = 1,\\
&&t_1^{-1}f_2t_1f_4^{-1}f_3^{-1}f_2^{-1}f_1^{-1}v_1^{-1}v_2^{-1}v_6^{-1} v_7^{-1} = 1,\\
&&a_1^{-1} f_2 a_1 f_8^{-1} = a_2^{-1} f_2 a_2 f_8^{-1} f_2^{-1} = 1,\\
&&r_1^{-1} f_3 r_1 f_4^{-1} f_3^{-1}f_1^{-1} v_1^{-1} v_2^{-1} v_6^{-1}v_7^{-1} = 1,\\
&&r_2^{-1} f_3 r_2 f_4^{-1} f_1^{-1} = t_1^{-1} f_3 t_1 f_4^{-1} f_2^{-1}v_2^{-1} v_4^{-1} v_5^{-1} v_6^{-1}v_7^{-1} = 1,\\
&&t_2^{-1}f_3t_2f_4^{-1}f_3^{-1}f_2^{-1}v_2^{-1} v_4^{-1}v_5^{-1}v_6^{-1}v_7^{-1} = 1,\\
&&a_1^{-1}f_3a_1f_8^{-1}f_4^{-1}f_6^{-1}f_4^{-1}v_2^{-1}v_4^{-1}v_5^{-1}v_6^{-1}v_7^{-1} = 1,\\
&&a_2^{-1} f_3 a_2 f_8^{-1} f_3^{-1}f_6^{-1} = r_1^{-1} f_4 r_1 f_4^{-1} f_2^{-1}v_2^{-1}v_4^{-1}v_5^{-1} v_6^{-1}v_7^{-1} = 1,\\
\end{eqnarray*}
\begin{eqnarray*}
&&r_2^{-1} f_4 r_2 f_4^{-1} f_2^{-1}v_2^{-1}v_4^{-1}v_5^{-1}v_6^{-1}v_7^{-1} = 1,\\
&&t_1^{-1} f_4 t_1 f_1^{-1} = 1,\quad (t_2, f_4^{-1}) = 1,\\
&&a_1^{-1} f_4 a_1 f_8^{-1} f_7^{-1}f_6^{-1} = a_2^{-1} f_4 a_2 f_8^{-1} f_7^{-1}f_6^{-1} f_4^{-1} = 1,\\
&&r_1^{-1}f_5r_1f_7^{-1}f_4^{-1}f_5^{-1} f_4^{-1}v_1^{-1}v_2^{-1}v_6^{-1}v_7^{-1} = 1,\\
&&r_2^{-1} f_5 r_2 f_8^{-1} f_5^{-1}v_2^{-1} v_4^{-1} v_5^{-1} v_6^{-1}v_7^{-1} = t_1^{-1} f_5 t_1 f_7^{-1} = 1,\\
&&t_2^{-1} f_5 t_2 f_8^{-1} f_5^{-1}v_2^{-1} v_4^{-1} v_5^{-1} v_6^{-1}v_7^{-1} = 1,\\
&&a_1^{-1}f_5 a_1f_8^{-1}f_3^{-1}f_8^{-1}f_5^{-1}f_2^{-1}f_1^{-1}v_2^{-1}v_4^{-1}v_5^{-1}v_6^{-1}v_7^{-1} = 1,\\
&&(a_2, f_5^{-1}) = (a_2, f_6^{-1}) = (a_2, f_7^{-1}) = 1,\\
&&r_1^{-1} f_6 r_1 f_8^{-1} f_5^{-1}v_2^{-1} v_4^{-1} v_5^{-1} v_6^{-1}v_7^{-1} = 1,\\
&&r_2^{-1}f_6 r_2f_7^{-1}f_4^{-1}f_5^{-1}f_4^{-1}v_1^{-1}v_2^{-1}v_6^{-1}v_7^{-1} = 1,\\
&&t_1^{-1} f_6 t_1 f_8^{-1} f_5^{-1}v_2^{-1}v_4^{-1}v_5^{-1}v_6^{-1}v_7^{-1} = 1,\\
&&t_2^{-1} f_6 t_2 f_7^{-1} = a_1^{-1} f_6 a_1 f_3^{-1}f_6^{-1}f_2^{-1} = r_1^{-1} f_7 r_1 f_5^{-1} = 1,\\
&&r_2^{-1}f_7r_2f_8^{-1}f_6^{-1}f_4^{-1}f_5^{-1}f_4^{-1}v_1^{-1}v_4^{-1}v_5^{-1} = 1,\\
&&t_1^{-1} f_7 t_1 f_7^{-1} f_4^{-1}f_5^{-1}f_4^{-1}v_1^{-1}v_2^{-1}v_6^{-1}v_7^{-1} = t_2^{-1} f_7 t_2 f_6^{-1} = 1,\\
&&a_1^{-1} f_7 a_1 f_4^{-1} f_7^{-1}f_3^{-1} = r_1^{-1} f_8 r_1 f_8^{-1} f_7^{-1}f_6^{-1} = 1,\\
&&(r_2, f_8^{-1}) = (t_2, f_8^{-1}) = (a_2, f_8^{-1}) = 1,\\
&&t_1^{-1} f_8 t_1 f_8^{-1} f_7^{-1}f_6^{-1} f_4^{-1} f_5^{-1} f_4^{-1} = a_1^{-1} f_8 a_1 f_8^{-1} f_2^{-1} = 1,\\
&&(f_1 f_2)^2 = (f_1 f_3)^2 = (f_2 f_3)^2 = (f_1 f_4)^2 = (f_2 f_4)^2 = (f_3 f_4)^2 = 1,\\
&&f_1 f_5 f_1 f_5 v_1^{-1}v_4^{-1} v_5^{-1} = f_2 f_5 f_2 f_5 v_2^{-1}v_4^{-1} v_5^{-1} v_6^{-1} v_7^{-1} = 1,\\
&&f_3 f_5 f_3 f_5 v_1^{-1}v_2^{-1} v_6^{-1} v_7^{-1} = f_1 f_6 f_1 f_6 v_1^{-1}v_4^{-1} v_5^{-1} = 1,\\
&&(f_5 f_6)^2 = (f_5 f_7)^2 = (f_6 f_7)^2 = (f_2 f_8)^2 = 1,\\
&&f_4 f_8 f_4 f_8 v_2^{-1}v_4^{-1} v_5^{-1} v_6^{-1} v_7^{-1} = 1,\\
&&(f_5 f_8)^2 = (f_6 f_8)^2 = (f_7 f_8)^2 = 1,\\
&&f_1 f_7 f_2 f_1 f_7 f_2 v_1^{-1} v_4^{-1} v_5^{-1} = 1,\\
&&(f_2 f_6 f_3)^2 = f_1 f_7 f_3 f_1 f_7f_3v_1^{-1} v_2^{-1}v_6^{-1}v_7^{-1} = 1,\\
&&f_1 f_8 f_3 f_1 f_8 f_3 v_2^{-1}v_4^{-1}v_5^{-1}v_6^{-1} v_7^{-1} = 1,\\
&&f_2 f_6 f_4 f_2 f_6 f_4 v_2^{-1} v_4^{-1}v_5^{-1}v_6^{-1} v_7^{-1} = 1,\\
&&f_1 f_7 f_4 f_1 f_7f_4 v_1^{-1} v_2^{-1} v_6^{-1}v_7^{-1} = f_4 f_6 f_5 f_4 f_5f_6 = 1,\\
&&f_4 f_7 f_5 f_4 f_5f_7 v_2^{-1} v_4^{-1} v_5^{-1}v_6^{-1} v_7^{-1} = f_3 f_8 f_6 f_3 f_6f_8 = 1.\\
\end{eqnarray*}

\item[\rm(m)] $H_2$ has an extension $H_1$ by $Z(Q) = \langle z
\rangle$ which is isomorphic to the finitely presented group
$$ H_1 = \langle a_1,a_2,r_1,r_2,t_1,t_2,f_j,q_i \mid
1 \le i, j \le 8 \rangle$$ with center $Z(H_1) = \langle z = q_6^2
\rangle$ having the set $\mathcal R(H_1)$ of defining relations
consisting of $\mathcal R(U)$, and the following relations:
\begin{eqnarray*}
&&a_1^3 = 1,\quad a_1^{a_2} = a_1^2,\quad (a_1, r_1) = (a_1, r_2) = (a_1, t_1) = (a_1, t_2) = 1,\\
&&a_1^{-1}  q_1  a_1  (q_1q_2q_4q_5q_6q_8)^{-1} = a_1^{-1}  q_2  a_1  (q_4q_5q_6q_8)^{-1} = 1,\\
&&a_1^{-1}  q_3  a_1  (q_2q_3q_4q_5q_6q_8)^{-1}  q_6^2 = a_1^{-1}  q_4  a_1  (q_2q_6)^{-1} = 1,\\
\end{eqnarray*}
\begin{eqnarray*}
&&a_1^{-1}  q_5  a_1  (q_8)^{-1} = a_1^{-1}  q_6  a_1  (q_2q_5q_8)^{-1} = 1,\\
&&a_1^{-1}  q_7  a_1  (q_2q_5q_6q_7q_8)^{-1}  q_6^2 = a_1^{-1}  q_8  a_1  (q_5q_8)^{-1}  q_6^2 = 1,\\
&&a_1^{-1}  f_1  a_1  f_6^{-1}  f_5^{-1}  (q_5)^{-1}  q_6^2 = 1,\\
&&a_1^{-1}  f_2  a_1  f_8^{-1}  (q_2q_4q_5q_6q_7)^{-1}  q_6^2= 1,\\
&&a_1^{-1}  f_3  a_1  f_8^{-1}  f_4^{-1}  f_6^{-1}  f_4^{-1}  (q_1q_2q_4q_6q_7)^{-1}  q_6^2 = 1,\\
&&a_1^{-1}  f_4  a_1  f_8^{-1}  f_7^{-1}  f_6^{-1}  (q_1q_2q_4q_6q_7)^{-1}  q_6^2 = 1,\\
&&a_1^{-1}f_5a_1f_8^{-1}f_3^{-1}f_8^{-1}f_5^{-1}f_2^{-1}f_1^{-1}(q_2q_4q_5q_6q_8)^{-1}q_6^2 = 1,\\
&&a_1^{-1}  f_6  a_1  f_3^{-1}  f_6^{-1}  f_2^{-1}  (q_1q_4q_7)^{-1}  q_6^2 = 1,\\
&&a_1^{-1}  f_7  a_1  f_4^{-1}  f_7^{-1}  f_3^{-1}  (q_1q_4q_5)^{-1} = 1,\\
&&a_1^{-1}  f_8  a_1  f_8^{-1}  f_2^{-1}  (q_2q_4q_5q_6q_7)^{-1}  q_6^2 = 1.\\
\end{eqnarray*}
In particular, $H_1 \leftarrow U \rightarrow D$ is an amalgam of
$H_1$ and $D$ with common subgroup $U$.

\item[\rm(n)] $H_1$ and $U$ have faithful permutation
representations $PH_1$ and $PU$ of degrees $6144$ and $2048$,
respectively, with common stabilizer
$$Y = \langle r_2r_1,(t_2q_6)^2, (q_6r_1f_4)^4, (q_6f_8r_1)^4, (r_1^2q_6f_8)^3,
(r_1t_1t_2f_4)^3, (r_1t_1q_6f_8)^2\rangle.$$

\end{enumerate}
\end{proposition}

\begin{proof}
(a) By Lemma \ref{l. M24-extensions}(h) and (i) the split
extension $E = E_1$ of $\M_{24}$ has a faithful permutation
representation $PE$ of degree $2048$ and a unique conjugacy class
of involutions of highest defect. It is represented by $z = v_1$.
Its centralizer $D = C_{E}(z)$ has order $2^{21}\cdot 3^2\cdot
5\cdot 7$. Using the faithful permutation representation $PE$ and
MAGMA it has been checked that the $3$ generators $x$, $y$ and $e$
of $E$ given in the statement satisfy $D = C_G(z) = \langle x, y
\rangle$ and $E = \langle D, e \rangle$.

(b) Using Kratzer's Algorithm 5.3.18 of \cite{michler}, the
faithful permutation representation $PE$ and MAGMA it has been
checked that $E$ has $80$ conjugacy classes and that $z = a^2 =
(xy^3)^{14}$.

The character table of $E$ has been computed by means of MAGMA
using the faithful permutation representation $PE$ of $E$.

(c) The restriction of $PE$ to $D = \langle x, y\rangle$ is a
faithful permutation representation of $D$. Using it and the MAGMA
command $\verb"NormalSubgroups(PD)"$ it follows that $D$ has a
unique non abelian normal subgroup $Q$ of order $2^9$ with center
$Z(Q) = \langle z \rangle$ of order $2$. Furthermore, $Q$ is
extra-special. The short words in $x$ and $y$ for its $8$
generators $q_i$, $1 \le i \le 8$, have been found by means of the
first author's program $\verb"GetShortGens(PD, PQ)"$, which is introduced in \cite{kim1}. It turns
out that all $8$ generators $q_i$ satisfy the set $\mathcal R(Q)$
of defining relations given in (c).

(d) Since $Q$ is a characteristic subgroup of $D$ its center
$Z(Q)$ is normal in $D$. Let $\alpha: D \rightarrow D_1 = D/Z(Q)$
be the canonical epimorphism with kernel $ker(\alpha) = Z(Q)$. Let
$V = \alpha(Q)$ and let $v_i = \alpha(q_i) \in D_1$ for $i =
1,2,...,8$. In particular, $\mathcal B = \{v_i|1 \le i \le 8\}$ is
a basis of the vector space $V$ over $F = \GF(2)$. Calculating now
$\alpha(q_i^x) \in V$ and $\alpha(q_i^y) \in V$ as short words in
the elements of the basis $\mathcal B$ one gets the $2$ matrices
$Mx$ and $My$ of $x$ and $y$ in $\GL_8(2)$ stated in the
assertion.

Clearly, $D_1/V \cong D/Q$. Applying the MAGMA command
$\verb"HasComplement(D_1,V)"$ in the permutation representation of
$D_1$ it follows that $D_1$ is a non split extension of $D/Q$ by
$V$.

(e) Let $MD = \langle Mx, My \rangle$ be the subgroup of
$\GL_8(2)$ generated by these $2$ matrices. Then the map $\varphi: D
\rightarrow MD$ defined by $x \rightarrow Mx$ and $y \rightarrow
My$ is an epimorphism from $D$ onto $MD$ with kernel $Q$ by (c).
It is well known that $\GL_8(2)$ has a faithful permutation
representation $PL8$ of degree $255$ whose stabilizer is a
parabolic subgroup of order $2^7|\GL_7(2)|$. Using it and MAGMA
one checks that the Fitting subgroup $MY$ of $MD$ has order $2^6$,
and that $MC = \langle \varphi(x),\varphi(c)\rangle$ is a complement
of $MY$ in $MD$, where $c = (x^2yx^2yxyxy^2xy)^3$ has order $2$.
Again the $6$ generators $\varphi(m_j)$ of $MY$ stated in
assertion (e) have been found by means of the program
$\verb"GetShortGens"$.

An application of the isomorphism testing program of Cannon and Holt \cite{cannon} implemented in MAGMA yields that $MC \cong A_8$. Using the MAGMA command
$\verb"NormalSubgroups(PD)"$ again one sees that $D$ has a unique
normal subgroup $M$ of order $2^{15}$. It is generated by $Q$ and
the $6$ elements $m_j$ given in the statement. Furthermore, an
application of MAGMA command $\verb"HasComplement(PD,M)"$ yields
that $C = \langle x, c \rangle$ is a complement of $M$ in $D$.
Hence $MY = \varphi(M)$ and $MC = \varphi(C) \cong C \cong A_8$.
In particular, $D$ has a faithful permutation representation of
degree $2^{15}$ with stabilizer $C$.

(f) As $ker(\alpha) = Z(Q) \le Q = ker(\varphi)$ also $\alpha(C)$ is
a complement of $\alpha(M) = \langle V, \alpha(m_j) \mid 1 \le j
\le 6 \rangle$ in $D_1$. In particular, $D_1$ has a faithful
permutation representation of degree $2^{14}$. Using it one
verifies that $C_{D_1}(V) = V$.

Clearly, $D_1/V \cong D/Q$. Applying the MAGMA command
$\verb"HasComplement(D_1,V)"$ in the permutation representation
$PD_1$ it follows that $D_1$ is a non split extension of $D/Q$ by
$V$.

(g) Since $|MD| = 2^8\cdot 3^2\cdot 5\cdot 7$ is fairly small
MAGMA is able to provide the following set $\mathcal R(MD)$ of
defining relations of $MD = \langle Mx, My \rangle$ by means of
its command $\verb"FPMD := FPGroup(sub<GL(12,2)|Mx,My>)"$:
\begin{eqnarray*}
&&x_1^7 = y_1^4 = 1,\\
&&(y_1x_1^{-1})^7 = 1,\quad (x_1y_1x_1^{-1}y_1x_1)^4 = 1,\\
&&(x_1y_1x_1^{-1}y_1x_1^{-1}y_1)^4 = 1,\\
&&(x_1y_1x_1^{-1}y_1)^6 = 1,\quad (x_1y_1x_1^{-1}y_1x_1^{-1}y_1x_1)^4 = 1,\\
&&x_1^{-1}y_1x_1^{-2}y_1x_1^2y_1x_1y_1x_1^{-1}y_1x_1^3y_1x_1^{-1}y_1x_1y_1x_1^2y_1x_1^{-2}y_1x_1^{-1}y_1 = 1,\\
&&x_1^{-2}y_1x_1^{-1}y_1x_1y_1x_1^{-2}y_1x_1^{-1}y_1x_1^{-3}y_1x_1^2y_1x_1^2y_1x_1y_1x_1^{-1}y_1x_1^{-1}y_1x_1^{-1} = 1,\\
&&y_1x_1y_1x_1^{-1}y_1x_1^{-2}y_1x_1^2y_1x_1^{-1}y_1x_1^{-2}y_1x_1y_1x_1^{-2}y_1x_1^2y_1x_1^{-2}y_1x_1^2 = 1,\\
\end{eqnarray*}

where $x_1 = \varphi(x)$ and $y_1 = \varphi(x)$. As $ker(\varphi) = Q$
all equations of $\mathcal R(MD)$ can be lifted to $D$ by
replacing $x_1$ and $y_1$ with $x$ and $y$, respectively, and
calculating the right hand side of these equations in $x$ and $y$
as words in the generators of $q_i$ of $Q$. Thus we obtained the
printed $9$ relations of statement (g).

Since $Q$ is normal in $D$ all $q_i^x, q_i^y \in Q$. Hence an
application of the first author's program (introduced in \cite{kim1})
$$\verb"LookupWord(PQ, q_i^d)"$$
with $d \in \{x,y\}$ yields all the relations of (g) involving
$q_i^x$ and $q_i^y$ for $1 \le i \le 8$. Now the set $\mathcal
R(Q)$ stated in (c) completes the set $\mathcal R(D)$ of defining
relations of $D$.

(h) All statements of this assertion have been obtained
computationally using the faithful permutation representations
$PL8$ of $\GL_8(2)$ and $PE$, MAGMA and the program
$\verb"GetShortGens"$.

(i) Let $\psi$ be the group
isomorphism between $\GL_8(2)$ and its image $PL8$ in the
symmetric group of degree $255$. Then $\psi$ induces an embedding
of $MD$ into $PL8$ and $PX = C_{PL8}(\psi(Mt))$ is the image of
$MX = C_{\GL_8(2)}(Mt)$ in $PL8$. Using then the faithful
permutation representation $PL8$ and the MAGMA command
$$\verb"exists(h){h:h in PX|Order(sub<PL8|MD, h>) eq 174182400}"$$
we obtain the matrix $Mh \in \GL_8(2)$ of order $3$ given in the
statement.

It also has been checked that $MO = \langle MD, Mh\rangle$ has
order $2^{12}\cdot 3^5\cdot 5^2\cdot 7$ and that $C_{MO}(Mt) =
\langle C_{MD}(Mt), Mh \rangle$ has order $2^{12}\cdot 3^3$.

From (e) and (h) follows now immediately that $$MH_3 = C_{MO}(Mt)
= \langle Mm_j, Mf_i,Mr_1,Mr_2,Mt_1,Mt_2,Mh \mid 1 \le i \le 4, 1
\le j \le 6\rangle.$$

The isomorphism $MO \cong O_8^{+}(2)$ has been found by an
application of the MAGMA command $\verb"CompositionFactors(MO)"$.

(j) Another application of MAGMA yields that the Fitting subgroup
$MY_3$ of $MH_3$ is extra-special of order $2^9$, and that it contains $Mm_3$, $Mm_4$, $Mm_6$, $Mf_1$, $Mf_2$, $Mf_3$, $Mf_4$. Then an application of the first author's program $\verb"GetShortGens"$ yields that $Mf_8$ together with these seven elements generate $MY_3$. The set
$\mathcal R(Y_3)$ of defining relations has been calculated by
means of MAGMA. In particular $(Mf_4Mf_7)^2$ is the center of
$MY_3$.

Let $MR = \langle r_1, r_2, t_1, t_2 \rangle$ and $T_1 = \langle
r_1, t_1\rangle$. The third generator $Ma_1$ of the Sylow
$3$-subgroup $MT$ of $MH_3$ has been found by means of the command

$$\verb"exists(a){a:a in C_{MH_3}(MR)|Order(sub<MR | T_1, a >) eq 3^3}".$$

A similar existence command provided the matrix $Ma_2$ in
$C_{MD}(MR)$ satisfying $Ma_1^{Ma_2} = Ma_1^2$. The word of the
element $a_2 \in D$ satisfying $\varphi(a_2) = Ma_2$ has been
obtained by an application of $\verb"LookupWord(PD,a_2)"$ command.
Using the faithful permutation representation $PL8$ it has been
checked that $MK_3 = \langle MR, Ma_1, Ma_2 \rangle$ is a
complement of $MY_3$ in $MH_3$. Using the structure of the
subgroups $MK_3$ and its action on the normal subgroup $MY_3$ of
$MH_3$ it is easy to verify all the relations of the finitely
presented group $H_3$ in the matrix group $MH_3$. Another
application of MAGMA yields that the two groups $MH_3$ and $H_3$ are
isomorphic.

(k) By the previous statements we know that the generators of $$U
= \langle a_2,r_1,r_2,t_1,t_2,f_j,q_i \mid 1 \le i,j \le 8
\rangle$$ are words in the generators $x$ and $y$ of $D$.
Therefore $U$ is a subgroup of $D$ containing the extra-special
normal subgroup $Q$ of order $2^9$. Let $P = \langle Q, f_j \mid 1
\le j \le 8 \rangle$. Then by (j) $P/Q \cong MY_3$ is extra-special of
order $2^9$ as well. Using the restriction of the faithful permutation representation
$PE$ to $D$ it has been verified that $K = \langle a_2, r_1, r_2,
t_1, t_2 \rangle$ is a complement of $P$ in $U$. Furthermore, it
has been used together with MAGMA to determine the set $\mathcal
R(U)$ of defining relations of $U$ given in the statement.

The stabilizer $Y$ of the faithful permutation representation $PU$
of degree $2048$ of $U$ given in the statement has been obtained
by application of the MAGMA command $\verb"BasicStabilizer(pU,2)"$
to the restriction $pU$ of the faithful permutation representation
$PE$ of $E$ to $U$ (recall that $U$ is a subgroup of $D$, hence
also of $E$), and the first author's program $\verb"GetShortGens"$
using the $21$ generators given above.

(l) The subgroup $MH_3 \cong H_3$ of $\GL_8(2)$ has a faithful
permutation representation inside $PL8$. In view of its
presentation given in (j) it satisfies all hypotheses of Holt's
Algorithm 7.4.5 of \cite{michler} implemented in MAGMA. Its
application in MAGMA yields that the second cohomological
dimension $dim_F(H^2(H_3,V)) = 1$, where $V = Q/Z(Q) = \langle v_i
\mid 1 \le i \le 8\rangle$. Hence there is a unique non split
extension $H_2$ of $H_3$ by $V$. The presentation of the non split
extension $H_2$ has also been calculated by means of Holt's
Algorithm using the matrices of the generators of $H_3$ and the
presentation of $H_3$ given in (j).

Furthermore, $U/Z(Q)$ is isomorphic to the subgroup $U_2 = \langle
a_2,r_1,r_2,t_1,t_2,f_j, v_i \mid 1 \le i \le 8, 1 \le j \le 8
\rangle$ of $H_2$. In particular, $U_2$ has a presentation with
set $\mathcal R(U_2)$ consisting of all relations of $\mathcal
R(H_2)$ not involving $a_1$. All these relations can easily be
lifted to $U$ by replacing the generators $v_i$ by $q_i$ for $1
\le i \le 12$ and checking whether the lifted word of the left
hand side of a relation of $U_2$ equals either $1$ or $z$ in
$Z(Q)$. Hence $\mathcal R(U)$ consists of all these lifted
relations and $\mathcal R(Q)$.

(m) Since MAGMA was not able to calculate the second cohomology
group $H^2(H_2,1_{FH_2})$ with coefficients in the trivial
$FH_2$-module $1_{FH_2}$ we constructed the central extension
$H_1$ of $H_2$ by the center $Z(Q) = \langle z = q_6^2 \rangle$ as
follows.

Let $H_1$ be any central extension of the finitely presented group
$H_2$ by $Z(Q)$. Then the finitely presented group $U$ constructed
in (k) may be considered to be a subgroup of $H_1$ containing the
common normal subgroup $Q = \langle q_i \mid 1 \le i \le
8\rangle$. Hence $\mathcal R(U)$ is a subset of the set $\mathcal
R(H_1)$ of defining relations of $H_1$. By (j) $K_3 = \langle
r_1,r_2, t_1,t_2, a_1,a_2 \rangle$ is a complement in $H_3$ of its normal
subgroup $MY_3 \cong P/Q = \langle f_i \mid 1 \le i \le 8
\rangle$. Also $K = \langle a_2,r_1,r_2,t_1,t_2\rangle$ is a
complement of $P$ in $U$ by (k). Since $H_1/Q \cong H_2/V \cong
H_3 = MY_3:K_3$ it follows that $|H_1 : U| = |K_3 : K| = 3$. As
$Z(Q)$ is a common normal subgroup of order $2$ of $U$ and $H_1$
we can apply Theorem 1.4.15 of \cite{michler}. It asserts that $P$
has also a complement $K_1 \geq K$ in $H_1$. Thus $K_1 \cong H_1/P
\cong H_3/MY_3 \cong K_3$. Thus the set $\mathcal R(K_3)$ of
defining relations  of $K_3$ is part of $\mathcal R(H_1)$.
Therefore it remains to lift the $16$ relations of $H_2$ given in
(l) involving $a_1$ and being different from the set of defining
relations $\mathcal R(K_3)$.

For that purpose a faithful permutation representation $PU$ of $U$
of degree $512$ has been constructed by means of the faithful
permutation representation of $D$ described in (e). Its stabilizer
$Y$ is generated by the $7$ elements of $U$ stated in assertion
(n).

In particular, $Y$ is a subgroup of $H_1$. Therefore $Y$ is also
a stabilizer of a faithful permutation representation $PH_1$ of
degree $6144$ of $H_1$ because $|H_1 : U| = 3$. Any possible
central extension $H_1$ of $H_2$ by $Z(Q)$ has a set $\mathcal
R(H_1)$ of defining relations consisting of $\mathcal R(U)$,
$\mathcal R(K_3)$ and one of the following $2^{16}$ sets $\mathcal
R(i)$ of relations:
\begin{eqnarray*}
&&a_1^{-1} q_1 a_1 (q_1q_2q_4q_5q_6q_8)^{-1}[z[i][1]],\\
&&a_1^{-1} q_2 a_1 (q_4q_5q_6q_8)^{-1} [z[i][2]],\\
&&a_1^{-1} q_3 a_1 (q_2q_3q_4q_5q_6q_8)^{-1} [z[i][3]],\\
&&a_1^{-1} q_4 a_1 (q_2q_6)^{-1} [z[i][4]],\\
&&a_1^{-1} q_5 a_1 (q_8)^{-1} [z[i][5]],\\
&&a_1^{-1} q_6 a_1 (q_2q_5q_8)^{-1} [z[i][6]],\\
&&a_1^{-1} q_7 a_1 (q_2q_5q_6q_7q_8)^{-1} [z[i][7]],\\
&&a_1^{-1} q_8 a_1 (q_5q_8)^{-1} [z[i][8]],\\
\end{eqnarray*}
\begin{eqnarray*}
&&a_1^{-1} f_1 a_1 f_6^{-1} f_5^{-1} (q_5)^{-1} [z[i][9]],\\
&&a_1^{-1} f_2 a_1 f_8^{-1} (q_2q_4q_5q_6q_7)^{-1} [z[i][10]],\\
&&a_1^{-1} f_3 a_1 f_8^{-1} f_4^{-1} f_6^{-1} f_4^{-1}(q_1q_2q_4q_6q_7)^{-1} [z[i][11]],\\
&&a_1^{-1} f_4 a_1 f_8^{-1} f_7^{-1} f_6^{-1}(q_1q_2q_4q_6q_7)^{-1} [z[i][12]],\\
&&a_1^{-1} f_5 a_1 f_8^{-1} f_3^{-1} f_8^{-1} f_5^{-1} f_2^{-1}f_1^{-1} (q_2q_4q_5q_6q_8)^{-1} [z[i][13]],\\
&&a_1^{-1} f_6 a_1 f_3^{-1} f_6^{-1} f_2^{-1} (q_1q_4q_7)^{-1}[z[i][14]],\\
&&a_1^{-1} f_7 a_1 f_4^{-1} f_7^{-1} f_3^{-1} (q_1q_4q_5)^{-1}[z[i][15]],\\
&& a_1^{-1} f_8 a_1 f_8^{-1} f_2^{-1} (q_2q_4q_5q_6q_7)^{-1}[z[i][16]],\\
\end{eqnarray*}

where $z[i]$ runs through all possible sequences $z[i]$ whose
entries $[z[i][j]], 1 \le j \le 16,$ are either $1$ or $z$ in $Z(Q)
= \langle q_6^2\rangle$. If $\mathcal R(H_1(i)) = \mathcal
R(U)\cup R(K_3) \cup \mathcal R(i)$ describes a set of defining
relations of a central extension $H_1(i)$ of $H_2$ by $Z(Q)$, then
$Y$ is a stabilizer of a faithful permutation representation
$PH_1(i)$ of degree $6144$.

Using MAGMA we checked the degree of $PH_1(i)$ for each of the $
2^{16}$ presentations $H_1(i)$. It turned out that there is
exactly one index $i$ where the sequence $z[i]$ provides a central
extension $H_1(i)$ of $H_2$ by $Z(Q)$ having a faithful permutation representation of degree $6144$ with stabilizer $Y$.
Using MAGMA it has been checked that the Sylow $2$-subgroups of
this $H_1(i)$ and $D$ are isomorphic. Thus $H_1$ is uniquely
determined. Its presentation is given in the statement.

(n) The generators of the common stabilizer $Y$ in $H_1$ and $U$
have been constructed in the proof of (k). This completes the
proof.
\end{proof}

The set of all faithful characters of a finite group $U$ is
denoted by $f \mbox{char}_{\mathbb{C}} (U)$, and $mf \mbox{char}_
{\mathbb{C}}(U)$ denotes the set of all multiplicity-free faithful
characters of $U$.

\begin{proposition}\label{prop. H(Co_1)} Keep the notation of Lemma \ref{l.
M24-extensions} and Proposition \ref{prop. DCo_1}. Let $E =
\langle x, y, e\rangle$ be the split extension of $\M_{24}$ by its
simple module $V_1$ of dimension $11$ over $F = \GF(2)$ and let $D
= \langle x,y \rangle$. Let $U = \langle
a_2,r_1,r_2,t_1,t_2,f_j,q_i \mid 1 \le i,j \le 8 \rangle $ be the
subgroup of $D$ constructed in Proposition \ref{prop. DCo_1}(k).
Let $H_1 \leftarrow U \rightarrow D$ be the amalgam constructed in
Proposition \ref{prop. DCo_1}, where $H_1 = \langle U,
a_1\rangle$. Let $\sigma: U \rightarrow D$ denote the
corresponding monomorphism of $U$ into $D$. Then the following
statements hold:

\begin{enumerate}
\item[\rm(a)] $U = \langle y, j, k \rangle$ and $H_1 = \langle y,
j, k, h \rangle$, where $y = a_2q_7r_1t_1t_2q_7r_2q_7$, $j =
r_2f_4$, $k = t_2q_6$ and $h = a_1$.

\item[\rm(b)] The Goldschmidt index of the amalgam $H_1 \leftarrow
U \rightarrow D$ is $1$.

\item[\rm(c)] A system of representatives $r_i$ of the $376$
conjugacy classes of $H_1$ and the corresponding centralizers
orders $|C_{H_1}(r_i)|$ are given in Table \ref{Co_1 cc H1}.

\item[\rm(d)] A system of representatives $u_i$ of the $382$
conjugacy classes of $U$ and the corresponding centralizers orders
$|C_U(u_i)|$ are given in Table \ref{Co_1ct_U}.

\item[\rm(e)] A system of representatives $d_i$ of the $155$
conjugacy classes of $D$ and the corresponding centralizers orders
$|C_D(d_i)|$ are given in Table \ref{Co_1ct_D}.

\item[\rm(f)] The character tables of $H_1$, $U$ and $D$ are given
in Tables \ref{Co_1ct_H_1}, \ref{Co_1ct_U} and \ref{Co_1ct_D},
respectively.

\item[\rm(g)] There is exactly one compatible pair $(\chi, \tau)
\in mf \mbox{char}_{\mathbb{C}}(H_1) \times mf
\mbox{char}_{\mathbb{C}}(D)$ of degree $128$ of the groups
$H_1=\langle U, h \rangle$ and $D =\langle U, x \rangle$:
$$
 (\chi_{\bf 235} , \tau_{\bf 39})
$$
with common restriction
$$
\chi_{|U} = \tau_{|U} = \psi_{\bf 270},
$$
where irreducible characters with bold face indices denote
faithful irreducible characters.

\item[\rm(h)] Let $\mathfrak V$ and $\mathfrak W$ be the up to
isomorphism uniquely determined faithful irreducible
$128$-dimensional modules of $H_1$ and $D$ over $F = \GF(23)$
corresponding to the compatible pair characters $\chi$ and $\tau$, respectively.

Let $\kappa_\mathfrak V : H_1 \rightarrow \GL_{128}(23)$ and
$\kappa_\mathfrak W : D \rightarrow \GL_{128}(23)$ be the
representations of $H_1$ and $D$ afforded by the modules
$\mathfrak V$ and $\mathfrak W$, respectively.

Let $\mathfrak h = \kappa_\mathfrak V(h)$, $\mathfrak j =
\kappa_\mathfrak V(j)$, $\mathfrak k = \kappa_\mathfrak V(k)$, $\mathfrak y = \kappa_\mathfrak V(y)$ in $
\kappa_\mathfrak V(H_1) \le \GL_{128}(23)$. Then the following
assertions hold:

\begin{enumerate}
\item[\rm(1)] $\mathfrak V_{|U} \cong \mathfrak W_{|U}$, and there
is a transformation matrix $\mathcal T \in \GL_{128}(23)$ such
that
$$
\mathfrak j = \mathcal T^{-1} \kappa_\mathfrak W (\sigma(j))
\mathcal T,\quad
\mathfrak k = \mathcal T^{-1} \kappa_\mathfrak W (\sigma(k))
\mathcal T, \quad \mbox{and} \quad
$$
$$
\mathfrak y = \mathcal T^{-1} \kappa_\mathfrak W (\sigma(y))
\mathcal T.
$$

Let $\mathfrak x = \mathcal T^{-1} \kappa_\mathfrak W(x)
\mathcal T$ in $\GL_{128}(23)$.

\item[\rm(2)] Let $\mathfrak H = \langle \mathfrak h, \mathfrak y,
\mathfrak x \rangle$. Its generating matrices and the
transformation matrix $\mathcal T$ 
can be downloaded from the first author's website:\\
$\verb"http://www.math.yale.edu/~hk47/Co1/index.html"$.
\end{enumerate}

\item[\rm(i)] $\mathfrak H$ is isomorphic to the finitely
presented group $H = \langle h_i|1 \le i \le 8\rangle$ with the
following set ${\mathcal R}(H)$ of defining relations:
\begin{eqnarray*}
&&h_1^4 = h_2^4 = h_3^4 = h_4^3 = h_6^7 = h_7^6 = h_8^6 = 1,\\
    &&(h_7^{-1} h_6^{-1})^3 = 1,\\
    &&h_8^{-2} h_6^{-1} h_7^{-1} h_6 h_7^{-1} h_8^{-1} h_7 = 1,\\
    &&h_7 h_8 h_6^{-1} h_7^{-1} h_6 h_7^2 h_8^{-1} = 1,\\
    &&(h_6 h_7^{-1} h_8 h_7^{-1})^2 = 1,\quad (h_6 h_8^{-1})^4 = 1,\\
    &&(h_7^{-1} h_8^{-1} h_6 h_8^{-1})^2 = 1,\quad (h_7 h_8^{-1})^4 = 1,\\
    &&h_7 h_8 h_6^{-3} h_8^{-2} h_7^{-1} h_6 = h_6 h_8 h_7 h_6^3 h_7 h_8 h_7 h_8^{-1} = 1,\\
    &&h_1 h_6 h_7 h_6 h_7^4 h_6 h_7 h_8^3 h_1 = h_6^{-1} h_1^{-1} h_6^2 h_7^2 h_8 h_7 h_6^3 h_1 = 1,\\
    &&h_7^{-1} h_1^{-1} h_6^3 h_7 h_6 h_7^2 h_6 h_7h_8^2 h_6^2 h_1 = 1,\\
    &&h_8^{-1} h_1^{-1} h_6^4 h_7 h_6^5 h_7 h_6 h_8h_7^2 h_1 = 1,\\
    &&h_3 h_5 h_6 h_7 h_6 h_7^2 h_6 h_7^3 h_8h_1 h_5^{-1} h_3 = 1,\\
    &&h_5 h_2^{-2} h_6 h_7^2 h_6 h_7^2 h_6^2 h_8h_7^2 h_2 h_1 h_6^{-1} = 1,\\
    &&h_1 h_4^{-1} h_2^{-1} h_6 h_7 h_6^2 h_7 h_6h_7 h_6^2 h_8 h_1 h_2 h_1^{-1} h_4^{-1} = 1,\\
    &&h_2 h_5 h_2^{-1} h_6^3 h_7 h_6^2 h_7 h_6h_7^4 h_8 h_6 h_7 h_2^{-1} h_7^{-1} = 1,\\
    &&h_2^{-1} h_5^{-1} h_2 h_6 h_7 h_6 h_7^3 h_8^2h_6 h_8 h_2 h_1^{-1} h_6^{-1} = 1,\\
    &&h_3^{-1} h_8^{-1} h_2 h_7 h_6 h_7^2 h_6^2 h_7h_8 h_6 h_1 h_2 h_5^2 = 1,\\
    &&h_2^{-1} h_1 h_3^{-1} h_6^5 h_7^4 h_6 h_7 h_8h_6 h_7 h_1 h_2 h_7 h_6 = 1,\\
    &&h_5 h_2^{-1} h_3^{-1} h_7 h_6 h_7^3 h_6^2 h_7^2h_8 h_6^2 h_2 h_5 h_6^{-1} = 1,\\
    &&h_1 h_3^{-2} h_6^2 h_7^2 h_6^4 h_7 h_8 h_1h_3^2 h_1 = 1,\quad h_2^{-1} h_3^{-2} h_2 h_3^2 = 1,\\
     \end{eqnarray*}
    \begin{eqnarray*}
    &&h_4^{-1} h_3^{-2} h_6^4 h_7 h_6 h_7 h_6^2 h_7h_6^4 h_8 h_4^{-1} h_1 = 1,\\
    &&h_6^{-1} h_3^{-2} h_6^2 h_7^4 h_6 h_7 h_6 h_7h_3^2 = 1,\\
    &&h_8 h_3^{-2} h_6^4 h_8 h_6 h_8 h_3^2 h_1 = 1,\\
    &&h_1 h_4^{-1} h_3^{-1} h_6^5 h_7 h_6 h_7^3 h_6^2h_7^2 h_8 h_1 h_3 h_4 h_1 = 1,\\
    &&h_3^{-1} h_4 h_3^{-1} h_6^6 h_7^3 h_6 h_1 h_4h_1^{-1} = 1,\\
    &&h_3 h_5^{-1} h_3^{-1} h_6^4 h_7^3 h_6^3 h_3 h_5h_3 = 1,\\
    &&h_3 h_6^{-1} h_3^{-1} h_6^2 h_7^3 h_6^5 h_3 h_6h_3 = 1,\\
    &&h_1 h_6 h_3^{-1} h_6^2 h_7^2 h_6^2 h_8 h_6h_8 h_6 h_1 h_2 h_1^{-1} h_2 = 1,\\
    &&h_3 h_6 h_3^{-1} h_6^5 h_7 h_6^3 h_7^2 h_8h_1 h_3 h_6^{-1} h_3 = 1,\\
    &&h_7^{-2} h_3^{-1} h_7^3 h_3^{-1} h_7^{-1} = 1,\\
    &&h_4 h_1 h_3 h_6^5 h_7 h_6^4 h_7 h_6 h_7h_6 h_8 h_1 h_3 h_4 h_1 = 1,\\
    &&h_3^{-1} h_4^{-1} h_3^2 h_4 h_3^{-1} = 1,\\
    &&h_1 h_6^{-1} h_3 h_6^6 h_7^4 h_6 h_7^2 h_6h_7 h_8 h_6 h_7 h_2^{-1} h_8^{-1} h_5^{-1} = 1,\\
    &&h_2 h_7 h_3 h_7^2 h_6 h_7^2 h_6^4 h_7 h_8h_6 h_7 h_5^{-1} h_8 = 1,\\
    &&h_4^{-1} h_1^{-1} h_4^{-1} h_6^5 h_7 h_8 h_6 h_8h_3^2 h_7 = 1,\\
    &&h_1 h_3^{-1} h_4^{-1} h_6^5 h_7 h_6 h_7^3 h_6^2h_7^2 h_8 h_1 h_4 h_3 h_1 = 1,\\
    &&h_2^{-1} h_3^{-1} h_4^{-1} h_6 h_7 h_6^2 h_7 h_6h_7 h_6^2 h_8 h_1 h_4 h_2 h_3 = 1,\\
    &&h_4 h_3^{-1} h_4^{-1} h_6^5 h_7^3 h_6^2 h_2^{-1} h_3h_2 = 1,\\
    &&h_3 h_5^{-1} h_4^{-1} h_7^3 h_4 h_5 h_3 = 1,\\
    &&h_2 h_1 h_4 h_6 h_7 h_6^3 h_7^3 h_6 h_7h_8^2 h_3^{-1} h_8^{-1} h_4 = 1,\\
    &&h_6 h_1^{-1} h_5^{-1} h_6^4 h_7^3 h_6 h_7 h_6h_2^{-1} h_3 h_5 = 1,\\
    &&h_5^{-1} h_2^{-1} h_5^{-1} h_7 h_6 h_7^3 h_6^2 h_7h_6 h_1 h_3^{-1} h_5 h_7 = 1,\\
    &&h_6 h_2^{-1} h_5^{-1} h_6 h_7 h_6 h_7 h_6^2h_7 h_6^2 h_7^2 h_1 h_2 h_8^{-1} h_2 = 1,\\
    &&h_8^{-1} h_3 h_5^{-1} h_7^4 h_6^2 h_7^2 h_8^2 h_6h_3 h_7 h_5 = 1,\\
    &&h_1 h_4^{-1} h_5^{-1} h_7^3 h_5 h_1^{-1} h_4^{-1} = 1,\\
    &&h_1 h_5^{-2} h_6^4 h_7^3 h_6 h_8 h_6 h_8 h_1h_2^{-1} h_1 h_5^{-1} = 1,\\
    &&h_2^{-1} h_5^{-2} h_6^6 h_7 h_8 h_7^3 h_2 h_5^{-1}h_6 = 1,\\
    &&h_5 h_6^{-1} h_5^{-1} h_7 h_6 h_7^2 h_6 h_8^2h_6 h_3 h_8^{-1} h_1^{-1} = 1,\\
    &&h_7 h_6^{-1} h_5^{-1} h_7 h_6^2 h_7^2 h_6 h_8h_7 h_6 h_1 h_3 h_6 h_8 = 1,\\
    &&h_8 h_6^{-1} h_5^{-1} h_6^2 h_7^4 h_6 h_7^2 h_6h_7^2 h_8^2 h_3^{-1} h_1 h_7^{-1} = 1,\\
    &&h_8^{-1} h_7 h_5^{-1} h_6^2 h_7^3 h_6^4 h_7 h_6h_8 h_7 h_8 h_1 h_3 h_8^{-1} h_2^{-1} = 1,\\
    &&h_3^{-1} h_8 h_5^{-1} h_6 h_7 h_6 h_7^4 h_6^4h_7 h_8 h_2 h_5^{-1} h_3^{-1} = 1,\\
    &&h_1 h_2 h_5 h_6^2 h_7^4 h_6 h_7 h_6 h_7^2h_6^2 h_8 h_1 h_5^{-1} h_1 h_2 = 1,\\
    &&h_1^{-1} h_3^{-1} h_5 h_6^2 h_7^2 h_6^2 h_7^2 h_8h_7 h_3 h_5 h_8^{-1} = 1,\\
    &&h_2 h_4 h_5 h_6 h_7^2 h_6 h_7^2 h_6 h_7^2h_5^{-1} h_4^{-1} h_2 = 1,\\
    &&h_3 h_5 h_4^{-1} h_2^{-1} h_6^5 h_7 h_6^3 h_7^2h_8 h_1 h_2 h_4 h_5^{-1} h_3 = 1,\\
    &&h_7^{-2} h_4^{-1} h_2^{-1} h_6 h_7^4 h_6 h_7 h_6h_7 h_2 h_4 h_7^{-1} = 1,\\
    &&h_3 h_7 h_4^{-1} h_2^{-1} h_6 h_7 h_6 h_7^4h_6 h_7 h_2 h_4 h_7^{-1} h_3 = 1,\\
    &&h_2 h_6 h_4 h_2^{-1} h_6^3 h_7 h_6^5 h_7^2h_6 h_2 h_4^{-1} h_6^{-1} h_2 = 1,\\
    \end{eqnarray*}
    \begin{eqnarray*}
    &&h_7 h_4^{-1} h_7^{-1} h_2^{-1} h_7 h_6^3 h_7 h_6^3h_7 h_8 h_1 h_2 h_4^{-1} h_7^{-1} h_4 = 1,\\
    &&h_2^{-1} h_4 h_3^{-1} h_2 h_6 h_7 h_6^2 h_7h_6 h_7^3 h_8^2 h_1 h_2 h_5^{-1} h_1^{-1} h_4^{-1} = 1,\\
    &&h_3 h_8 h_4 h_2 h_6 h_7 h_6 h_7^4 h_6h_7 h_8^3 h_2^{-1} h_4^{-1} h_8^{-1} h_3 = 1,\\
    &&h_1 h_6^{-1} h_4^{-1} h_3^{-1} h_6^5 h_7 h_6 h_7^3h_6^2 h_7^2 h_8 h_1 h_3 h_4 h_6 h_1 = 1,\\
    &&h_1 h_6 h_4^{-1} h_3^{-1} h_6^5 h_7 h_6 h_7^3h_6^2 h_7^2 h_8 h_1 h_3 h_4 h_6^{-1} h_1 = 1,\\
    &&h_1 h_7^{-1} h_4^{-1} h_3^{-1} h_6^5 h_7 h_6 h_7^3h_6^2 h_7^2 h_8 h_1 h_3 h_4 h_7 h_1 = 1,\\
    &&h_1 h_7 h_4^{-1} h_3^{-1} h_6^5 h_7 h_6 h_7^3h_6^2 h_7^2 h_8 h_1 h_3 h_4 h_7^{-1} h_1 = 1,\\
    &&h_3^{-1} h_7 h_4^{-1} h_3^{-1} h_6^2 h_7^2 h_6^2 h_7h_6 h_8 h_1 h_2^2 h_4 h_7^{-1} = 1,\\
    &&h_1 h_8^{-1} h_4^{-1} h_3^{-1} h_6^5 h_7 h_6 h_7^3h_6^2 h_7^2 h_8 h_1 h_3 h_4 h_8 h_1 = 1,\\
    &&h_1 h_8 h_4^{-1} h_3^{-1} h_6^5 h_7 h_6 h_7^3h_6^2 h_7^2 h_8 h_1 h_3 h_4 h_8^{-1} h_1 = 1,\\
    &&h_6 h_5^{-1} h_1^{-1} h_4^{-1} h_6 h_7 h_6^5 h_7^2h_6 h_8^2 h_6 h_1 h_2 h_4^{-1} h_5^{-1} h_1^{-1} = 1,\\
    &&h_3^{-1} h_8^{-1} h_1 h_4^{-1} h_6^4 h_7 h_6 h_7h_6^2 h_7 h_6^4 h_8 h_4^{-1} h_8 h_3^{-1} = 1,\\
    &&h_3^{-1} h_7^{-1} h_3^{-1} h_4^{-1} h_6 h_7 h_6 h_7^4h_6 h_7 h_3^{-1} h_4 h_7 h_3^{-1} = 1,\\
    &&h_4^{-1} h_7^{-1} h_3 h_4^{-1} h_6^2 h_7^2 h_6^5 h_7h_1 h_2 h_6^{-1} h_4 h_7 = 1,\\
    &&h_1^{-2} h_5^{-1} h_4^{-1} h_6^5 h_7 h_6 h_7^3 h_6^2h_7^2 h_8 h_1 h_4 h_1^2 h_5 = 1,\\
    &&h_1 h_2^{-1} h_5 h_4^{-1} h_6^2 h_7 h_6^5 h_7h_6^4 h_7 h_4 h_5^{-1} h_2 h_1 = 1,\\
    &&h_1 h_4 h_5 h_4^{-1} h_6^4 h_7 h_6 h_7 h_6^2h_7 h_6^4 h_8 h_1 h_4 h_5^{-1} h_1^{-1} h_4 = 1,\\
    &&h_1 h_4^{-1} h_6^{-1} h_4^{-1} h_6^4 h_7 h_6 h_7h_6^2 h_7 h_6^4 h_8 h_1 h_4 h_6 h_4 h_1 = 1,\\
    &&h_1 h_4 h_6^{-1} h_4^{-1} h_6^5 h_7 h_6^4 h_7h_6 h_7 h_6 h_8 h_1 h_4 h_6 h_4^{-1} h_1 = 1,\\
    &&h_6 h_4 h_6^{-1} h_4^{-1} h_6 h_7^4 h_6 h_7h_8^2 h_6 h_2^{-1} h_4^{-1} h_2^{-1} h_4^{-1} = 1,\\
    &&h_6^{-1} h_2^{-1} h_6 h_4^{-1} h_7^4 h_6^2 h_7 h_6h_8^2 h_6 h_8 h_1 h_2 h_4^{-1} h_3 h_1^{-1} = 1,\\
    &&h_6^{-1} h_3^{-1} h_6 h_4^{-1} h_6^4 h_7^4 h_8 h_6h_7 h_3^{-1} h_6 h_4^{-1} h_1^{-1} = 1,\\
    &&h_3 h_5^{-1} h_6 h_4^{-1} h_7^3 h_4 h_6^{-1} h_5h_3 = 1,\\
    &&h_1 h_3^{-1} h_7^{-1} h_4^{-1} h_6^2 h_7 h_6^5 h_7h_6^4 h_7 h_4 h_7 h_3 h_1 = 1,\\
    &&h_6 h_1 h_7 h_4^{-1} h_6 h_7 h_6 h_7^3 h_6h_8 h_7^2 h_8 h_4 h_5 h_2 h_3 = 1,\\
    &&h_2^{-2} h_7 h_4^{-1} h_6^2 h_7 h_6 h_7 h_6^2h_7 h_6^2 h_4 h_7^{-1} h_3 = 1,\\
    &&h_3^{-1} h_2^{-1} h_7 h_4^{-1} h_6^3 h_7 h_6^3 h_7h_6^4 h_7 h_8 h_1 h_4 h_7^{-1} h_2^{-1} = 1,\\
    &&h_1 h_2^{-1} h_8 h_4^{-1} h_6^2 h_7 h_6^5 h_7h_6^4 h_7 h_4 h_8^{-1} h_2 h_1 = 1,\\
    &&h_7^{-1} h_2^{-1} h_1^{-1} h_4 h_6^4 h_8 h_7^3 h_2h_4 h_7 h_3^{-1} = 1,\\
    &&h_7^{-1} h_5^{-1} h_1^{-1} h_4 h_6^6 h_7 h_6 h_7^4h_6 h_7 h_6^2 h_1 h_4 h_5 h_7^{-2} = 1,\\
    &&h_4 h_6 h_1 h_4 h_6^6 h_7 h_6 h_7^3 h_6h_7 h_6 h_7 h_8 h_6 h_4 h_7^{-1} h_8^{-1} h_3 = 1,\\
    &&h_3^{-1} h_8 h_2^{-1} h_4 h_7 h_6^4 h_7 h_6^4h_8 h_6 h_2 h_4^{-1} h_3 h_8^{-1} = 1,\\
    &&h_6^{-1} h_5 h_3^{-1} h_4 h_7^2 h_6^5 h_7 h_8h_7 h_6 h_3^{-1} h_7^{-1} h_4^{-1} h_5 = 1,\\
    &&h_6 h_5 h_3 h_4 h_6 h_7 h_6^4 h_7^2 h_6h_7^2 h_8 h_6 h_2^{-1} h_4^{-1} h_5^{-1} h_8 = 1,\\
    &&h_3 h_1^{-1} h_5^{-1} h_4 h_6^5 h_7 h_6 h_7h_6^4 h_7 h_6^2 h_8 h_1 h_4^{-1} h_5 h_1 h_3 = 1,\\
    &&h_5 h_1 h_5^{-1} h_4 h_6 h_7 h_6^2 h_7 h_6^3h_7^2 h_8 h_7 h_8 h_1 h_3 h_7 h_3 h_4^{-1} = 1,\\
    &&h_1 h_2 h_5^{-1} h_4 h_6^5 h_7^3 h_6^2 h_4^{-1}h_5 h_2^{-1} h_1 = 1,\\
    &&h_2^{-1} h_3^{-1} h_5^{-1} h_4 h_6^3 h_7 h_6^3 h_7h_6^3 h_7 h_4^{-1} h_5 h_2 h_3^{-1} = 1,\\
    &&h_5^{-1} h_2 h_6^{-1} h_4 h_6^3 h_7^3 h_6 h_8h_7^3 h_4 h_2 h_3 h_8 = 1,\\
    &&h_2^{-1} h_3^{-1} h_6^{-1} h_4 h_6^3 h_7 h_6^3 h_7h_6^3 h_7 h_4^{-1} h_6 h_2 h_3^{-1} = 1,\\
    &&h_2 h_3^{-1} h_6^{-1} h_4 h_6 h_7^4 h_6 h_7h_6 h_7 h_4^{-1} h_6 h_2^{-1} h_3^{-1} = 1,\\
    &&h_5 h_3^{-1} h_6^{-1} h_4 h_7 h_6 h_7 h_6^4h_7 h_6^3 h_1 h_4^{-1} h_5 h_3^{-1} h_6^{-1} = 1,\\
    &&h_2 h_3^{-1} h_6 h_4 h_6^2 h_7^3 h_6^5 h_4^{-1}h_6^{-1} h_2^{-1} h_3 = 1,\\
    &&h_6^{-1} h_3^{-1} h_6 h_4 h_6^5 h_7 h_6 h_7^2h_6^2 h_7 h_6^2 h_8 h_7 h_3 h_2 h_4^{-1}h_2^{-1} = 1,\\
    &&h_6 h_1 h_8^{-1} h_4 h_6 h_7^2 h_6 h_7 h_6h_8 h_6^2 h_1 h_4^{-1} h_5^{-1} h_8 h_2^{-1} = 1,\\
    \end{eqnarray*}
    \begin{eqnarray*}
    &&h_1 h_5 h_8^{-1} h_4 h_6^5 h_7^3 h_6^2 h_4^{-1}h_8 h_5^{-1} h_1 = 1,\\
    &&h_1 h_5 h_4^{-1} h_5^{-1} h_6 h_7 h_6 h_7^4h_6 h_7 h_8^3 h_5 h_4 h_5^{-1} h_1 = 1,\\
    &&h_3 h_8 h_4^{-1} h_5^{-1} h_6^5 h_7 h_6 h_7^4h_6 h_7 h_6^3 h_5 h_4 h_8^{-1} h_3 = 1,\\
    &&h_2^{-1}h_5^{-1}h_4^{-1}h_1^{-1}h_2^{-1}h_6^4h_7^2h_6 h_7^2 h_6 h_7 h_8 h_6 h_4 h_6^{-1}h_3^{-1} h_8 = 1,\\
    &&h_3 h_6 h_4^{-1} h_1^{-1} h_2^{-1} h_6^5 h_7 h_6h_7 h_6^2 h_7 h_6^3 h_8 h_4^{-1} h_5 h_2 h_6 = 1.\\
\end{eqnarray*}
Moreover, $h_1,h_2,h_3,h_4,h_5$ are the isomorphic images of $\mathfrak y, \mathfrak j, \mathfrak k, \mathfrak h, \mathfrak x$, respectively.

\item[\rm(j)] $H$ has a faithful permutation representation of
degree $276480$ with stabilizer $\langle h_5,
(h_2^3h_5^2h_2h_5)^6, (h_2^2h_5^2h_2^2h_5h_2)^4,
 (h_2^2h_5h_2^3h_5^3)^6 \rangle$.

\item[\rm(k)] Each Sylow $2$-subgroup $S$ of $H$ has a unique
maximal elementary abelian normal subgroup $V$ of order $2^{11}$
and $N_H(V) \cong D = C_{E}(z)$.

\end{enumerate}
\end{proposition}

\begin{proof}
(a) The $2$ generators $j$ and $k$ of $U$ have been obtained by
means of the first author's program $\verb"GetShortGens"$ and the faithful permutation presentation $PU$ of
degree $2048$ constructed in Proposition \ref{prop. DCo_1}(k).

(b) The Goldschmidt index has been calculated by means of
Kratzer's Algorithm 7.1.10 of \cite{michler}, the faithful
permutation representations $PH_1$, $PU$, $PD$  and MAGMA.

(c), (d) and (e) The systems of representatives of the conjugacy
classes of $H_1$, $U$ and $D$ have been calculated by means of
$PH_1$, $PU$ and $PD$, MAGMA and Kratzer's Algorithm 5.3.18 of
\cite{michler}. Only for (c), the new set of generators $j_1 = yj,
k_1 = yk, h_1 = yh$ of $H_1$ are used in order to speed up the
application of Kratzer's Algorithm in our calculations.

(f) The three character tables have been obtained by means of
MAGMA using the faithful permutation representations of the $H_1$,
$U$ and $D$. The character tables of $H_1$, $D$ and $U$ are stated
in the appendix.

(g) Using MAGMA and the faithful permutation representations again
we have determined the fusion of the classes of $U$ in $H_1$ and
$\sigma(U)$ in $D$. A MAGMA calculation employing Kratzer's
Algorithm 7.3.10 of \cite{michler} shows then that the amalgam
$H_1 \leftarrow U \rightarrow D$ has a unique compatible pair
$(\chi, \tau)$ of degree $128$. It is stated in the assertion.

(h) In order to construct the faithful irreducible representation
$\mathfrak V$ corresponding to the character $\chi = \chi_{\bf
235}$ we determine the irreducible constituents of the faithful
permutation representation $(1_Y)^{H_1}$ of $H_1$ of degree $6144$
with stabilizer $Y$ given in Proposition  \ref{prop. DCo_1}(n).
Calculating inner products with the irreducible characters of
$H_1$ it follows $(1_Y)^{H_1}$ has $30$ irreducible constituents
including the character $\chi_{\bf 235}$ of degree $128$. Since
the permutation module $(1_Y)^{H_1}$ is only of degree $6144$ the
Meataxe algorithm can be applied to get the corresponding
irreducible module $\mathfrak V$, because $23$
does not divide the order of $H_1$ and the character $\chi_{\bf
235}$ has rational values.

In order to find the irreducible module $\mathfrak W$
corresponding to the irreducible character $\tau = \tau_{\bf 39}$ of
$D$ with degree $128$ we determined the $7$ irreducible
constituents of the permutation module $(1_T)^D$ of $D$ of degree
$2048$ with stabilizer $T = \verb"BasicStabilizer(PD,2)"$. It
follows that $\tau_{\bf 39}$ belongs to them. Applying the Meataxe
algorithm to the permutation module one obtains the irreducible
representation $\mathfrak W$ belonging to the constituent
$\tau_{\bf 39}$ of the compatible pair $(\chi, \tau)$.

By construction the $KU$-modules $\mathfrak V_{|U}$ and $\mathfrak
W_{|U}$ described by the $2$ triples $$(\kappa_ \mathfrak V(j),\kappa_\mathfrak
V(k),\kappa_\mathfrak V(y))\quad \mbox{and} \quad  (\kappa_\mathfrak
W(\sigma(j)),\kappa_\mathfrak W(\sigma(k)), \kappa_\mathfrak W(\sigma(y)))$$ of
matrices in $\GL_{128}(23)$ are isomorphic. Let $Y = \GL(128,23)$,
let $V$ and $W$ be the $KU$- and $K(\sigma(U))$-modules described by
the first and the second triple of matrices of $Y$. Applying
Parker's isomorphism test of Proposition 6.1.6 of \cite{michler}
by means of the MAGMA command
$$\verb"IsIsomorphic(GModule(sub<Y|W>),GModule(sub<Y|V>))",$$
one gets the transformation matrix $\mathcal T$ in $\GL_{128}(23)$
satisfying $\kappa_ \mathfrak V(y) = (\kappa_\mathfrak W(\sigma(y)))^{\mathcal
T}$, $\kappa_\mathfrak V(j) = (\kappa_\mathfrak W(\sigma(j)))^{\mathcal T}$ and
$\kappa_\mathfrak V(k) = (\kappa_\mathfrak W(\sigma(k)))^{\mathcal T}$.

Using the faithful permutation representation $PD$ of $D = \langle
x, y \rangle$ with stabilizer $T$ it has been checked that $x$
does not belong to $\sigma(U)$ and that $D = \langle \sigma(U), x
\rangle = \langle \sigma(y), x\rangle$. Let $\mathfrak x = (\kappa_ \mathfrak W(x))^{\mathcal T}$ and $\mathfrak y = \kappa_ \mathfrak V(y)$. Then
$\mathfrak D = \langle \mathfrak x, \mathfrak y \rangle \cong D$.
Let $\mathfrak H = \langle \mathfrak h,\mathfrak y, \mathfrak x
\rangle$.

(i) Denote the matrix group $\mathfrak H$ by $MH$. Subsequent
applications of the MAGMA commands $\verb"RandomSchreier(MH)"$ and
$\verb"BasicStabilizerChain(MH)"$ provided a base with $3$ base
points, and a chain of $3$ basic stabilizers $B(k)$ of orders
$$89181388800, 43008, 21504.$$ Using then the MAGMA command
$\verb"FPGroupStrong(MH)"$ it provided the presentation of a
finitely presented group $H \cong MH$ with the $8$ generators
$h_i$ and relations given in the statement. Moreover, this last MAGMA command preserves the original generators as first few generators of the resulting presentation, so we have that $h_1,h_2,h_3,h_4,h_5$ are the isomorphic images of $\mathfrak y, \mathfrak j, \mathfrak k, \mathfrak h, \mathfrak x$, respectively (in MAGMA, we actually put $\mathfrak H = \langle \mathfrak y, \mathfrak j, \mathfrak k, \mathfrak h, \mathfrak x\rangle$ and used it as $MH$).

(j) 
%
%
%
We took the third basic stabilizer $B(3)$ as obtained in (i), and checked by means of MAGMA that $B(3) = \langle h_6, h_7, h_8\rangle$ (we checked this in the matrix group $MH$). Then by applying P. Young's program
$$\verb"h_H_PH,PH:=MyCosetAction(H,sub<H|h_6,h_7,h_8>:maxsize:=50000000)"$$
we obtained the faithful permutation representation $P_1H$ of
degree $4147200$. In particular, we obtained $|H| = |P_1 H| = 2^{21} \cdot 3^5 \cdot 5^2 \cdot 7$, hence $|H : D| = 135$.

Let $pD$ be the restriction of the faithful permutation representation $PE$ of degree $2048$ of $E$ to $D$. By the MAGMA command $\verb"BasicStabilizer(pD,2)"$ and the first author's program $\verb"GetShortGens"$, we obtained the generators of stabilizer subgroup of $D$ giving $2048$ degree faithful permutation representation of $D$, in terms of $y,j,k,x$, i.e., in terms of $h_1,h_2,h_3,h_5$, as given in the statement. This subgroup works as a permutation stabilizer for $P_1H$, yielding the faithful permutation representation $PH$ of $P_1H$ (hence of $H$) of degree $2048 \cdot 135 = 276480$.

(k) By Proposition \ref{prop. DCo_1}(m) $H$ has a subgroup $U$ of
odd index which is isomorphic to the subgroup $\sigma(U)$ of $D$
with index $|D : \sigma(U)| = 35$. Hence it contains a Sylow $2$-subgroup $S_1$ of $H$ which is isomorphic to the ones of $D =
C_E(z)$. As $V_1$ is the unique maximal elementary abelian normal
subgroup of the Sylow $2$-subgroup $\sigma(S_1)$ in $D$ by Lemma
\ref{l. M24-extensions} also $S_1$ has a unique elementary abelian
normal subgroup $V$ of order $2^{11}$. Using the faithful
permutation representation $PH$ of $H$ constructed in (j) it has
been checked that $N_H(V) \cong D$. This completes the proof.
\end{proof}

\newpage

\section{Construction of Conway's simple group $\Co_1$}

By Proposition \ref{prop. H(Co_1)} the amalgam $H \leftarrow D
\rightarrow E$ constructed in section $3$ satisfies the main
condition of Algorithm 2.4 \cite{michler1}. Therefore we can apply
Algorithm 7.4.8 of \cite{michler} to give here a new existence
proof for Conway's sporadic group $\Co_1$, see \cite{conway}. It
is realized as a matrix group $\mathfrak G$ inside
$\GL_{276}(23)$. For that construction, we need a smaller faithful
permutation representation of the centralizer $H$ determined in
Proposition \ref{prop. H(Co_1)} than the one given there. This is
achieved by the following embedding of the centralizer of a
$2$-central involution of the second Conway group $\Co_2$
constructed in \cite{kim}. It is used also to embed $\Co_2$ into
$\Co_1$ which then allows us to construct a faithful permutation
representation of the matrix group $\mathfrak G$ in
$\GL_{276}(23)$ of degree $98280$.

\begin{proposition}\label{prop. H(Co_2 in Co_1)}
Keep the notation of Lemma \ref{l. M24-extensions} and
Propositions \ref{prop. H(Co_1)} and \ref{prop. DCo_1}. Let $H =
\langle h_i \mid 1 \le i \le 8 \rangle = \langle x, y, h \rangle $
be the finitely presented group of Proposition \ref{prop. H(Co_1)}
with center 
$z = (xy^3)^{14}$
, where $x = h_5$, $y = h_1$ and $h =
h_4$. Let
\begin{eqnarray*}
&&q_1 = y^2, \quad q_2 = (xy)^7,\quad q_3 = (yx)^7,\quad q_4 = (xy^2)^7,\quad q_5 = (yxy)^7, \\
&&q_6 = (x^5yx)^7, \quad q_7 = (x^4yx^2)^7, \quad q_8 = (x^4yxy)^6.\\
\end{eqnarray*}

Then the following statements hold:
\begin{enumerate}
\item[\rm(a)] $Q = \langle q_i \mid 1 \le i \le 8 \rangle$ is an
extra-special normal subgroup of $H$ with center $Z(Q) = \langle
q_6^2\rangle$.

\item[\rm(b)] The subgroup $K$ generated by
\begin{eqnarray*}
&&(xh^2x^2y)^6,\quad (h^2xyhx)^6,\quad(hxyhyxy)^6,\\
&&(xyhx^2yx^2)^4,\quad (xyhxh^2x^2)^6,\quad (xyh^2x^3h)^4\\
\end{eqnarray*}
is a simple group of order $2^9\cdot 3^4\cdot 5\cdot 7$. It is
isomorphic to $Sp_6(2)$ and $Q \cap K = 1$.

\item[\rm(c)] $H$ has a faithful permutation representation $PH$
of degree $61440$ with stabilizer $K$.

\item[\rm(d)] $C = QK$ has center $Z(C) = Z(Q)$, and $C$ is
isomorphic to the centralizer $H(\Co_2)$ of a $2$-central
involution $t$ of the simple group $\Co_2$ constructed in
\cite{kim}.

\item[\rm(e)] $H$ has $167$ conjugacy classes
${h_i}^H$ with representatives
$\mathfrak{h_i}$ and centralizer orders
$|C_H(h_i)|$ as given in Table \ref{Co_1 cc H}.

\item[\rm(f)] The character table of $H$ is stated in Table
\ref{Co_1ct H}.

\item[\rm(g)] The Goldschmidt index of the amalgam $H \leftarrow D
\rightarrow E$ is $1$.
\end{enumerate}
\end{proposition}

\begin{proof}
(a) This assertion is a restatement of Proposition \ref{prop.
DCo_1}(c).

(b) Using MAGMA and the faithful permutation representation of Proposition
\ref{prop. H(Co_1)}(j) one sees that $H$ has exactly
$3$ conjugacy classes of subgroups $R$ of order $2^{18}\cdot
3^4\cdot5\cdot7$ containing $Q$. Applying then the MAGMA command
$\verb"HasComplement(R,Q)"$ it follows that $Q$ has a complement
$K$ only in one class of subgroups $R$. The generators of $K$ have
been found by the first author's program
$\verb"GetShortGens"$.
Another application of MAGMA yields that $K$ is simple and
isomorphic to $Sp_6(2)$. In particular, $K \cap Q = 1$.

(c) Furthermore, $H$ has a faithful permutation representation
$PH$ of degree $|H:K| = 61440$ with stabilizer $K$.

(d) Clearly $C = QK$ is a subgroup of $H$ because $Q$ is normal in
$H$. Using $PH$ and MAGMA it is easy to see that $Z(C) = Z(Q)$.
The remaining statements follow from an isomorphism test between
$C = HK$ and $H(\Co_2)$ using MAGMA and the faithful permutation
representations $PH$ of $H$ and the one of $H(\Co_2)$ given in
Lemma 3.4 (a) of \cite{kim}.

(e) The system of representatives of the conjugacy classes of $H$
has been calculated by means of $PH$, MAGMA and Kratzer's
Algorithm 5.3.18 of \cite{michler}.

(f) The character table of $H$ has been computed by means of MAGMA
using the faithful permutation representation $PH$ of $H$.

(g) The Goldschmidt index of the amalgam $H \leftarrow D
\rightarrow E$ has been calculated by means of Kratzer's Algorithm
7.1.10 of \cite{michler}, the faithful permutation representations
$PH$, $PD$, $PE$  and MAGMA.
\end{proof}

\begin{theorem}\label{thm. existenceCo_1}
Keep the notation of Lemma \ref{l. M24-extensions} and
Propositions \ref{prop. H(Co_1)}, \ref{prop. DCo_1}, and \ref{prop. H(Co_2 in Co_1)}. Let $K =
\GF(23)$. Let $\sigma: D \rightarrow E$ denote the corresponding
monomorphism of $D$ into $E$. Using the notation of the $3$
character tables \ref{Co_1ct H}, \ref{Co_1ct_D} and \ref{Co_1ct_E}
of the groups $H$, $D$ and $E$, respectively, the following
statements hold:

\begin{enumerate}
\item[\rm(a)] There is exactly one compatible pair $(\chi, \tau)
\in mf \mbox{char}_{\mathbb{C}}(H) \times mf
\mbox{char}_{\mathbb{C}}(E)$ of degree $276$ of the groups
$H=\langle D, h \rangle$ and $E =\langle D, e \rangle$:
$$
\chi_{2}+ \chi_{10} + \chi_{\bf 11} = \tau_{\bf 10}
$$

with common restriction
$$
\tau_{|D} = \chi_{|D} = \psi_{11}+\psi_{37} +\psi_{\bf 39},
$$
where irreducible characters with bold face indices denote
faithful irreducible characters.

\item[\rm(b)] Let $\mathfrak V$ and $\mathfrak W$ be the up to
isomorphism uniquely determined faithful semi-simple
multiplicity-free $276$-dimensional modules of $H$ and $E$ over $F
= \GF(23)$ corresponding to the compatible pair characters $\chi$ and $\tau$,
respectively.

Let $\kappa_\mathfrak V : H \rightarrow \GL_{276}(23)$ and
$\kappa_\mathfrak W : E \rightarrow \GL_{276}(23)$ be the
representations of $H$ and $E$ afforded by the modules $\mathfrak
V$ and $\mathfrak W$, respectively.

Let $\mathfrak h = \kappa_\mathfrak V(h)$, $\mathfrak x =
\kappa_\mathfrak V(x)$, $\mathfrak y = \kappa_\mathfrak V(y)$ in $
\kappa_\mathfrak V(H) \le \GL_{276}(23)$. Then the following
assertions hold:

\begin{enumerate}
\item[\rm(1)] $\mathfrak V_{|D} \cong \mathfrak W_{|D}$, and there
is a transformation matrix $\mathcal T \in \GL_{276}(23)$ such
that
$$
\mathfrak x = \mathcal T^{-1} \kappa_\mathfrak W (\sigma(x))
\mathcal T, \mathfrak y = \mathcal T^{-1} \kappa_\mathfrak
W(\sigma(y)) \mathcal T.
$$
Let $\mathfrak e = \mathcal T^{-1} \kappa_\mathfrak W(e) \mathcal
T \in \GL_{276}(23)$.

Let $\mathfrak G = \langle \mathfrak h, \mathfrak x, \mathfrak y,
\mathfrak e \rangle$, and $\mathfrak H = \langle \mathfrak h, \mathfrak x, \mathfrak y\rangle = \kappa_ \mathfrak V (H)$.

\item[\rm(2)] The $4$ generating matrices of $\mathfrak G$ and the
transformation matrix $\mathcal T$ 
can be downloaded from the first author's website:\\
$\verb"http://www.math.yale.edu/~hk47/Co1/index.html"$.

\item[\rm(3)] In $\mathfrak G$ the subgroup $\mathfrak C = \langle
\mathfrak p, \mathfrak y, \mathfrak r, \mathfrak e \rangle$ is the
stabilizer of a $23$-dimensional $\mathfrak C$-invariant subspace
${\mathfrak A}$ of $\mathfrak U = K^{276}$ such that the
$\mathfrak G$-orbit ${\mathfrak A}^{{\mathfrak G}}$ has degree
$98280$, where $\mathfrak p = ((\mathfrak x\mathfrak y)^2\mathfrak
x)^3$ and $\mathfrak r = (\mathfrak x\mathfrak y\mathfrak h)^5$.

Furthermore, $\mathfrak C$ is isomorphic to the sporadic Conway
group $\Co_2$. The $2$-central involution centralizer of
$\mathfrak C$ is the subgroup $H(\mathfrak C) = \langle \mathfrak
p, \mathfrak y, \mathfrak r\rangle$ of $\mathfrak C$.


\item[\rm(4)] The character table of $\mathfrak G$ coincides with
that of $\Co_1$ in the Atlas \cite{atlas}, its p. 184 -185.

\end{enumerate}

\item[\rm(c)] $\mathfrak G$ is a finite simple group with
$2$-central involution $\mathfrak z = \kappa_\mathfrak V(z) =(\mathfrak x
\mathfrak y^3)^{14} $ such that
$$
C_{\mathfrak G} (\mathfrak z) = \mathfrak H = \kappa_\mathfrak V(H),~\text{and}~|\mathfrak G| = 2^{21}\cdot 3^9\cdot 5^4\cdot
7^2\cdot 11\cdot 13\cdot 23.
$$

\end{enumerate}
\end{theorem}

\begin{proof}
(a) Using MAGMA and the faithful permutation representations of
$PE$ and $PH$ determined in Lemma \ref{l. M24-extensions} and
Proposition \ref{prop. H(Co_2 in Co_1)}, respectively, we
determined the fusion of the classes of $D$ in $H$ and $\sigma(D)$ in $E$. A
MAGMA calculation employing Kratzer's Algorithm 7.3.10 of
\cite{michler} shows then that the amalgam $H \leftarrow D
\rightarrow E$ has a unique compatible pair $(\chi, \tau)$ of
degree $276$. It is stated in the assertion.


(b) The faithful irreducible constituent $\mathfrak V_1 =
\mathfrak V(\chi_{\bf 11})$ of dimension $128$ of the faithful
irreducible representation $\mathfrak V$ corresponding to the
semi-simple character $\chi_{\bf 11}$ of the compatible pair has already
been constructed in Proposition \ref{prop. H(Co_1)}(h). Its second
tensor power contains each of the irreducible representations
$\mathfrak V_2$ and $\mathfrak V_3$ corresponding to $\chi_2$ and
$\chi_{10}$, respectively, with multiplicity $1$. Since the second
tensor power of $\mathfrak V_1$ has degree $16384$ over $K
=\GF(23)$ the Meataxe algorithm can be applied to yield the $2$
irreducible representations $\mathfrak V_2$ and $\mathfrak V_3$.

In order to construct the faithful irreducible representation
$\mathfrak W$ corresponding to the character $\tau = \tau_{\bf
10}$ of degree $276$ we employed the MAGMA command
$\verb"LowIndexSubgroups(PE, 1000)"$ using the faithful
permutation presentation $PE$ of $E$ of degree $2048$. MAGMA found
a subgroup $R$ of index $552$ such that $\tau_{\bf 10}$ is a
constituent of the permutation character $(1_R)^E$. Using then a
stand alone program of the first author which is based on
Algorithm 5.7.1 of \cite{michler} we obtained the irreducible
representation $\mathfrak W$ over $K$ corresponding to $\tau_{\bf 10}$.

By construction the $KD$-modules $\mathfrak V_{|D}$ and $\mathfrak
W_{|D}$ described by the $2$ pairs $$(\kappa_\mathfrak V(x),\kappa_\mathfrak
V(y))\quad \mbox{and} \quad (\kappa_\mathfrak W(\sigma(x)),\kappa_\mathfrak
W(\sigma(y)))$$ of matrices in $\GL_{276}(23)$ are isomorphic
where $\sigma$ denotes the group isomorphism between the subgroup
$D = \langle x, y \rangle$ of $H$ and subgroup $D = \langle
\sigma(x),\sigma(y)\rangle $ of $E$. Let $Y = \GL_{128}(23)$, and let
$V$ and $W$ be the $KD$- and $K(\sigma(D))$-modules described by the
first and the second pair of matrices of $Y$. Applying Parker's
isomorphism test of Proposition 6.1.6 of \cite{michler} by means
of the MAGMA command
$$\verb"IsIsomorphic(GModule(sub<Y|W>),GModule(sub<Y|V>))",$$
one gets the transformation matrix $\mathcal T$ in $\GL_{276}(23)$
satisfying $\kappa_ \mathfrak V(y) = (\kappa_ \mathfrak W(\sigma(y)))^{\mathcal
T}$ and $\kappa_ \mathfrak V(x) = (\kappa_ \mathfrak W(\sigma(x)))^{\mathcal T}$.

Since $E = \langle \sigma(x), \sigma(y), e \rangle$ we let
$\mathfrak e = \mathcal T^{-1}\kappa_\mathfrak W(e)\mathcal T$ in
$\GL_{276}(23)$. Furthermore, $\mathfrak G = \langle \mathfrak x,
\mathfrak y, \mathfrak h, \mathfrak e \rangle$.

In order to construct a small faithful permutation representation
of $\mathfrak G$ we now embed the sporadic group $\Co_2$ into
$\mathfrak G$. By Proposition \ref{prop. H(Co_2 in Co_1)} $z =
(xy^3)^{14}$ generates the center $Z(H)$ of $H = \langle
x,y,h\rangle$ and there is a subgroup $C = QK$ of $H$ whose
generators are given in the statements (a) and (b) of that
proposition such that $C$ is isomorphic to the $2$-central
involution centralizer $H(\Co_2)$ constructed in \cite{kim}.
Let $D_C \le C$ be the image of $D(\Co_2) \le H(\Co_2)$ under this
isomorphism, where $D(\Co_2)$ is as constructed in \cite{kim}. By
the MAGMA command \\
$$\verb"exists(b){b:b in H|((D_C)^b meet D) eq (D_C)^b and y in (D_C)^b}"$$
we obtained an element $b\in H$ such that $(D_C)^b \le D = \langle
x , y \rangle$ and $y \in (D_C)^b$. Then, the first author's
program $\verb"GetShortGens"$ yielded the two elements $p =
((xy)^2 x)^3$ and $r = (xyh)^5$ such that $C^b = \langle
y,p,r\rangle$ and $(D_C)^b = \langle y,p\rangle$. Since $p$ is a
word in $x$ and $y$, we were able to construct a subgroup $D_E =
\langle y, p \rangle$ in $E = \langle x,y,e \rangle$. By the MAGMA
command\\
$$\verb"exists(b'){b':b' in E|Order(sub<E|D_E, b'>) eq 2^18*3^2*5*7*11}"$$
we obtained an element $b' \in E$ such that $|\langle D_E,
b'\rangle| = 2^{18} \cdot 3^2 \cdot 5 \cdot 7 \cdot 11 =
|E(\Co_2)|$, where $E(\Co_2)$ is as constructed in \cite{kim}. Let
$E_1 = \langle D_E, b' \rangle$. By means of MAGMA
command $\verb"IsIsomorphic"$ we checked that $E_1 \cong E(\Co_2)$ and
$E_1 = \langle y,p,e \rangle$.

By Theorem 4.2 of \cite{kim} the simple sporadic group $\Co_2$ is an
epimorphic image of the free product $C^b *_{(D_C)^b}E_1$ of $C^b$
and $E_1$ with amalgamated subgroup $(D_C)^b = \langle y, p\rangle$.
Therefore we let $\mathfrak C = \langle \mathfrak p, \mathfrak y,
\mathfrak r, \mathfrak e\rangle$. Then $\mathfrak C$ is a subgroup
of $\mathfrak G$, because its $4$ generators are words in the
generating matrices $\mathfrak x$, $\mathfrak y$, $\mathfrak h$
and $\mathfrak e$ of $\mathfrak G$.

Using the Meataxe algorithm it follows that $\mathfrak U|_{\mathfrak C}$ splits into $2$ irreducible $K\mathfrak C$-modules $\mathfrak A$ and $\mathfrak B$ of dimensions $23$ and
$253$, respectively. Applying then Weller's Algorithm described in
Theorem 6.2.1 of \cite{michler} and made compatible with MAGMA by
the first author we obtained a faithful permutation representation
$P\mathfrak G$ of $\mathfrak G$ of degree $98280$. Using then
Proposition 5.2.14 of \cite{michler} and MAGMA it follows that
$\mathfrak C$ is the stabilizer of $P\mathfrak G$. Furthermore,
$P\mathfrak G$ and MAGMA have been used to get the order
$|\mathfrak C| = 2^{18}\cdot 3^6\cdot 5^3 \cdot 7\cdot 11\cdot 23$. Hence
$|\mathfrak G| = |\mathfrak C|\cdot 98280$.

The system of representatives of the conjugacy classes of
$\mathfrak G$ has been calculated by means of $P\mathfrak G$,
MAGMA and Kratzer's Algorithm 5.3.18 of \cite{michler}.

The character table of $\mathfrak G$ has been computed by means of
MAGMA using the faithful permutation representation $P\mathfrak G$
of $\mathfrak G$.

(c) The matrix group $\mathfrak G$ is simple by its character
table. 
We checked by means of MAGMA and the faithful permutation representation $P\mathfrak G$ of $\mathfrak G$ that $C_\mathfrak G (\mathfrak z) = \mathfrak H$.
\end{proof}

In order to show that the simple  group $\mathfrak G$ constructed in
Theorem \ref{thm. existenceCo_1} is isomorphic to the original Conway group $\Co_1$
we use the following result of L. Soicher's article
\cite{soicher1}.

\begin{theorem}[L. Soicher]\label{thm. presentCo_1}
Conway's simple group $\Co_1$ is isomorphic to the finitely
presented group $G = \langle a,b,c,d,e,f,g,h,i\rangle$ with the
following set $\mathcal R(G)$ of defining relations:

\begin{eqnarray*}
&&a^2 = b^2 = c^2 = d^2 = e^2 = f^2 = g^2 = h^2 = i^2 = 1,\\
&&(ab)^3 = (bc)^3 = (cd)^8 = (de)^3 = (ef)^3 = (fg)^3 = (gh)^3 = (hi)^3 = 1,\\
&&(ac)^2 = (ad)^2 = (ae)^2 = (af)^2 = (ag)^2 = (ah)^2 = (ai)^2 = 1,\\
&&(bd)^2 = (be)^2 = (bf)^2 = (bg)^2 = (bh)^2 = (bi)^2 = 1,\\
&&(ce)^2 = (cf)^2 = (cg)^2 = (ch)^2 = (ci)^2 = (df)^2 = (dg)^2 = 1,\\
&&(dh)^2 = (di)^2 = (eg)^2 = (eh)^2 = (ei)^2 = (fh)^2 = (fi)^2 = (gi)^2 = 1,\\
&&a = (cd)^4,\quad (bcde)^8 = 1, \quad ((bcdcdefgh)^{13}i)^3 = 1.\\
\end{eqnarray*}

Furthermore, $G$ has a faithful permutation representation $PG$
of degree $1545600$ with stabilizer $Q = \langle a,b,c,d,e,f,g,h \rangle$.
\end{theorem}

\begin{proof}
Using the program $\verb"MyCosetAction(G,Q: maxsize:=10000000)"$
we obtained a faithful permutation
representation $PG$ of degree $1545600$ with stabilizer $Q$. Thus
we observed that $G$ has the same order as the simple group
$\mathfrak G$ constructed in Theorem \ref{thm. existenceCo_1}. The
isomorphism between the  original Conway group $\Co_1$ and $G$ is quoted from Soicher's article
\cite{soicher1}.
\end{proof}

\begin{corollary}\label{cor. identCo_1} 
Keep the notations for $\mathfrak x, \mathfrak y, \mathfrak h,
\mathfrak e, \mathfrak G, \mathfrak H, \mathfrak z$ from Theorem
\ref{thm. existenceCo_1} and for $G = \langle
a,b,c,d,e,f,g,h,i\rangle$ from Theorem \ref{thm. presentCo_1}. In
$G$ let $u_1=c$, $u_3=(cde)^3$, $u_5=(cdefg)^5$,
$u_7=(cdefghi)^7$, and $w_1=(u_3b)^4$, $w_2=(u_1u_7b)^4$,
$w_3=(u_1u_3u_5b)^4$, $w_4=(u_7bu_5b)^5$. Furthermore, let $h_0
=(w_4w_2w_1)^3$, $v=(w_3w_1w_2w_4)^5$, $w=(w_3w_1w_4w_2)^5$,
$x_0=h_0(vh_0^2)^2$, $y_0=h_0(h_0w^2v)^2$ , $q =
(bxbx^2bx^2b)^{30}$ and $e_0=(x^2yq)^5$. Then the following
assertions hold:

\begin{enumerate}
\item[\rm(a)] The generator $a$ is a $2$-central involution of the
finitely presented group $G$ such that
$$C_G(a) = \langle w_i\mid 1 \le i \le 4 \rangle \cong \mathfrak H = C_{\mathfrak G}(\mathfrak z),$$
where $\mathfrak z = (\mathfrak x \mathfrak y^3)^{14}$ is the $2$-central involution of the simple
group $\mathfrak G = \langle \mathfrak x, \mathfrak y, \mathfrak
h, \mathfrak e \rangle$ constructed in Theorem \ref{thm.
existenceCo_1}.

\item[\rm(b)] $\mathfrak G \cong G$.

\item[\rm(c)] $G$ has a faithful permutation representation of degree
$98280$ with stabilizer $C_0 = \langle h_0, x_0, y_0, e_0\rangle$.

\item[\rm(d)] $C_0$ is isomorphic to Conway's sporadic group $\Co_2$.
\end{enumerate}
\end{corollary}

\begin{proof}
(a) Let $PG$ be the faithful permutation representation of Soicher's
group $G$ defined in Theorem \ref{thm. presentCo_1}. Using it and
MAGMA we observed that $H_0 = C_G(a)$ is a centralizer of the
$2$-central involution $a$ of $G$ and $|H_0| = 2^{21}\cdot3^5\cdot5^2\cdot7$. Applying the program
$\verb"GetShortGens(G,H_0)"$ we found the given generators
of $H_0$.

Let $\mathfrak H = C_{\mathfrak G}(\mathfrak z)$ be the centralizer of the $2$-central involution $\mathfrak z$ of the simple
group $\mathfrak G$ constructed in Theorem \ref{thm. existenceCo_1}. By Proposition \ref{prop. H(Co_2
in Co_1)} it has a faithful permutation representation $PH = P\mathfrak H$ of degree $61440$.
Using the MAGMA command $\verb"DH_0:=DegreeReduction(H_0)"$ one can verify that $H_0$ has a faithful permutation representation $DH_0$ of degree $573440$. 
By means of the MAGMA command $\verb"IsIsomorphic(PH,DH_0)"$ we observed that $H_0 \cong \mathfrak H$. This last calculation took $508823$ seconds, which is $\approx5.9$ days.

(b) By Theorem \ref{thm. presentCo_1} the groups $\mathfrak G$ and
$G$ have the same character table stated in the Atlas
\cite{atlas}, its p. 184 -185. Therefore also $G$ has a uniquely
determined irreducible $276$-dimensional representation $\mathfrak
G_0$ over $\GF(23)$ by Brauer's Theorem 3.12.4 of \cite{michler}.
By (a) and the character table, both groups have isomorphic Sylow
$2$-subgroups and a unique conjugacy class of $2$-central
involutions. Let $S_0$ be a Sylow $2$-subgroup of $H_0 = C_G(a)$.
It has exactly one maximal elementary abelian normal subgroup $B$
by Proposition \ref{prop. H(Co_1)}(k). by means of the faithful
permutation representation $PG$ of Theorem \ref{thm. presentCo_1}
it has been checked computationally that $D_0 = N_{H_0}(B) =
C_{E_0}(a)$ where $E_0 = N_G(B)$. We also verified that $N_G(B)$
is isomorphic to the split extension $E$ constructed in Lemma
\ref{l. M24-extensions} and that $G = \langle H_0, N_G(B)\rangle$.
By Proposition \ref{prop. H(Co_2 in Co_1)} and (a) the amalgam
$H_0 \leftarrow D_0 \rightarrow E_0$ has Goldschmidt index $1$.
Therefore Theorem 7.5.1 of \cite{michler} asserts that $G \cong
\mathfrak G$.

(c) The existence of a faithful permuation representation $PG$ of $G$ of degree $98280$ is clear by Theorem \ref{thm. existenceCo_1} and (b). It remains to construct short generators of an appropriate stabilizer $C_0$ in $G$ in terms of the given generators of $G$. Using the faithful permutation representation $DH_0$ of degree $573440$ of $H_0 = \langle w_1,w_2,w_3,w_4\rangle$ and the MAGMA command $\verb"LowIndexSubgroups(DH_0, <120,120>)"$ we found $3$ conjugacy classes of such subgroups $R_i$ in $H_0$. Using the program $\verb"GetShortGens(H_0,R_i)"$ we got short generators for fixed representatives $R_i$ of each class. Thus we saw that only the subgroup $R_2 = \langle h_0, v, w\rangle$ is a split extension of $Sp_6(2)$ by an extra-special group $FR_2$ of order $2^9$,  where $FR_2$ is the Fitting Subgroup of $R_2$. The words for $h_0$, $v$ and $w$ are stated in the hypothesis. $R_2$ has a faithful permutation representation $pR_2$ of degree $14336$ obtained by the command $\verb"pR_2:=DegreeReduction(DR_2)"$, where $DR_2$ is the restriction of the faithful permutation representation $DH_0$ of $H_0$ to $R_2$. Using it and Corollary 4.3 of \cite{kim} we verified that $R_2 = \langle h_0,v,w\rangle$ is isomorphic to a centralizer $H(\Co_2)$ of Conway's second sporadic group $\mathfrak G_3$ constructed in Theorem 4.2 of \cite{kim}. In $R_2$ we chose a fixed Sylow $2$-subgroup $S_2$ and calculated a system of $5$ generators $s_1=(h_0wh_0)^6$, $s_2=(h_0vh_0w)^3$, $s_3=(h_0w^2v)^4$, $s_4=vwh_0v$ and
$s_5=(h_0v^2w^2)^3$, by means of the first author's program $\verb"GetShortGens(R_2,S_2)"$. Applying the MAGMA command
$$\verb"Subgroups(S_2: Al:=Normal, IsElementaryAbelian := true)"$$
we observed that $S_2$ has a unique maximal elementary abelian
normal subgroup $A_2$ of order $2^{10}$. The program
$\verb"GetShortGens(S_2,A_2)"$ showed that it is generated by
$a_1=(s_2)^2$, $a_2=(s_3)^2$, $a_3=(s_4)^4$, $a_4=(s_1s_3)^2$,
$a_5=(s_2s_4)^2$, $a_6=(s_4s_2)^2$, $a_7=(s_1s_2s_4)^4$,
$a_8=(s_1s_2s_3s_2)^2$, $a_9=(s_1s_2s_4s_3s_4)^2$ and
$a_{10}=s_1s_2s_4s_2s_3s_4$. Furthermore, $D_0 = N_{R_2}(A_2) =
\langle x_0, y_0\rangle$, where $x_0$ and $y_0$ of respective
orders $10$ and $12$ are defined in the hypothesis. They were
found by means of MAGMA and the program
$\verb"GetShortGens(R_2,D_0)"$. Using MAGMA and the large faithful
permutation representation $PG$ we determined $E_1 = N_G(A_2)$. It
has order $2^{19}\cdot3^2\cdot5\cdot7\cdot11$ and its Fitting
subgroup $F_2$ is elementary abelian of order $2^{11}$. By another
application of MAGMA and its command
$$\verb" exists(u){u: u in E_1|Order(sub<E_1|D_2,u>) eq 2^18*3^2*5*7*11}"$$
we obtained an element $u \in E_1$ of order $7$ such that $E_0 = \langle D_0, u\rangle$ has order $2^{18}\cdot3^2\cdot5\cdot7\cdot11$ and has $A_2$ as its Fitting subgroup. In particular, $E_0$ is a subgroup of $N = N_G(F_2)$. Using the MAGMA command $\verb"IsIsomorphic"$ it is checked that $N$ is isomorphic to the local subgroup $E = E(\Co_1)$ constructed in Lemma 2.1. In particular, $N \cong N_G(B)$. Similarly, we checked that $E_0$ is isomorphic to $E(\Co_2)$ constructed in Lemma 2.5(f) of \cite{kim}.

Using MAGMA and the faithful permutation representation $PG$ we ran the program $\verb"GetShortGens(G,N)"$. Thus we found the element $q \in N$ of order $2$ defined in the hypothesis, which satisfies $N = \langle x_0,y_0,q\rangle$. Another applications of $\verb"GetShortGens"$ yielded the elements $e_0=(x_0^2y_0q)^5$ and $e_1=(qx_0q)^5$ of orders $3$ and $2$ respectively, which satisfy $E_0 =\langle x_0,y_0,e_0\rangle$ and $E_1 = \langle x_0,y_0,e_0,e_1\rangle$. Furthermore, we checked by means of MAGMA that the subgroup $C_0 = \langle R_2, E_0\rangle = \langle h_0, x_0, y_0, e_0\rangle$ of $G$ is simple and has index $|G:C_0| = 98280$. Hence (c) holds.

(d) The subgroup $C_0 = \langle h_0,x_0,y_0,e_0\rangle$ of $G$ has a faithful permutation representation $DC_0$ of degree $2300$ which we obtained by means of the MAGMA command $\verb"DC_0:=DegreeReduction(PC_0)"$, where $PC_0$ is the restriction of $PG$ to $C_0$. Using it and the faithful permutation representation of $\Co_2$ stated in the proof of Corollary 4.3 of our article \cite{kim} an isomorphism test with MAGMA shows that $C_0 \cong \Co_2$.
\end{proof}

\newpage

\section{Representatives of conjugacy classes}

\setlength\textheight{48\baselineskip}
\addtolength\textheight{\topskip}
\begin{cclass}\label{Co_1cc E} Conjugacy classes of $E(\Co_1) = \langle x,y,e \rangle$
{ \setlength{\arraycolsep}{1mm}
\renewcommand{\baselinestretch}{0.5}
\scriptsize
 $$
\begin{array}{|c|r|r|c|c|c|c|c|c|
}\hline
 \mbox{Class} &
 \multicolumn{1}{c|}{\mbox{Representative}} &
 \multicolumn{1}{c|}{|\mbox{Centralizer}|}
 &
 2{\rm P}
 & 3{\rm P}
 & 5{\rm P}
 & 7{\rm P}
 & 11{\rm P}
 & 23{\rm P}
 \\\hline
 1 & 1 & 2^{21}\cdot3^{3}\cdot5\cdot7\cdot11\cdot23
 & 1 & 1 & 1 & 1 & 1 & 1 \\\hline

 2_{1} & (y)^{2} & 2^{21}\cdot3^{2}\cdot5\cdot7
 & 1 & 2_{1} & 2_{1} & 2_{1} & 2_{1} & 2_{1} \\\hline

 2_{2} & (xey^{2})^{11} & 2^{18}\cdot3^{3}\cdot5\cdot11
 & 1 & 2_{2} & 2_{2} & 2_{2} & 2_{2} & 2_{2} \\\hline

 2_{3} & e & 2^{17}\cdot3\cdot7
 & 1 & 2_{3} & 2_{3} & 2_{3} & 2_{3} & 2_{3} \\\hline

 2_{4} & (x^{2}ey)^{3} & 2^{15}\cdot3\cdot5
 & 1 & 2_{4} & 2_{4} & 2_{4} & 2_{4} & 2_{4} \\\hline

 2_{5} & (xexey)^{3} & 2^{15}\cdot3\cdot5
 & 1 & 2_{5} & 2_{5} & 2_{5} & 2_{5} & 2_{5} \\\hline

 2_{6} & (xyxeyxye)^{2} & 2^{17}\cdot3
 & 1 & 2_{6} & 2_{6} & 2_{6} & 2_{6} & 2_{6} \\
\hline

\end{array}
$$}
\normalsize

\newpage
Conjugacy classes of $E(\Co_1) = \langle x,y,e \rangle$(continued)
{ \setlength{\arraycolsep}{1mm}
\renewcommand{\baselinestretch}{0.5}
\scriptsize
 $$
 \begin{array}{|c|r|r|c|c|c|c|c|c|
}\hline
 \mbox{Class} &
 \multicolumn{1}{c|}{\mbox{Representative}} &
 \multicolumn{1}{c|}{|\mbox{Centralizer}|}
 &
 2{\rm P}
 & 3{\rm P}
 & 5{\rm P}
 & 7{\rm P}
 & 11{\rm P}
 & 23{\rm P}
 \\\hline
 14_{1} & xy^{2} & 2^{3}\cdot7
 & 7_{1} & 14_{2} & 14_{2} & 2_{1} & 14_{1} & 14_{1} \\\hline

 14_{2} & (xy^{2})^{3} & 2^{3}\cdot7
 & 7_{2} & 14_{1} & 14_{1} & 2_{1} & 14_{2} & 14_{2} \\\hline

 14_{3} & xy & 2^{2}\cdot7
 & 7_{1} & 14_{4} & 14_{4} & 2_{3} & 14_{3} & 14_{3} \\\hline

 14_{4} & (xy)^{3} & 2^{2}\cdot7
 & 7_{2} & 14_{3} & 14_{3} & 2_{3} & 14_{4} & 14_{4} \\\hline

 15_{1} & (x^{2}y)^{2} & 2\cdot3\cdot5
 & 15_{1} & 5 & 3_{1} & 15_{2} & 15_{2} & 15_{1} \\\hline

 15_{2} & (x^{2}y)^{14} & 2\cdot3\cdot5
 & 15_{2} & 5 & 3_{1} & 15_{1} & 15_{1} & 15_{2} \\\hline

 16_{1} & x^{3}yxe & 2^{6}
 & 8_{3} & 16_{1} & 16_{1} & 16_{1} & 16_{1} & 16_{1} \\
\hline

 \end{array}
 $$
}
\end{cclass}

\normalsize


\newpage

\begin{cclass}\label{Co_1 cc H} Conjugacy classes of $H(\Co_1) = \langle x,y,h \rangle$
{ \setlength{\arraycolsep}{1mm}
\renewcommand{\baselinestretch}{0.5}
\scriptsize
$$
 \begin{array}{|c|r|r|r|c|c|c|c|
}\hline
 \mbox{Class} &
 \multicolumn{1}{c|}{\mbox{Representative}} &
 \multicolumn{1}{c|}{|\mbox{Class}|} &
 \multicolumn{1}{c|}{|\mbox{Centralizer}|}
 &
 2{\rm P}
 & 3{\rm P}
 & 5{\rm P}
 & 7{\rm P}
 \\\hline
 1 & 1 & 1 & 2^{21}\cdot3^{5}\cdot5^{2}\cdot7
 & 1 & 1 & 1 & 1 \\\hline

 2_{1} & (yh)^{4} & 1 & 2^{21}\cdot3^{5}\cdot5^{2}\cdot7
 & 1 & 2_{1} & 2_{1} & 2_{1} \\\hline

 2_{2} & (y)^{2} & 270 & 2^{20}\cdot3^{2}\cdot5\cdot7
 & 1 & 2_{2} & 2_{2} & 2_{2} \\\hline

 2_{3} & (xyhxh^{2})^{6} & 12600 & 2^{18}\cdot3^{3}
 & 1 & 2_{3} & 2_{3} & 2_{3} \\\hline

 2_{4} & (xyh)^{5} & 60480 & 2^{15}\cdot3^{2}\cdot5
 & 1 & 2_{4} & 2_{4} & 2_{4} \\\hline

 2_{5} & (xh^{2}y)^{5} & 60480 & 2^{15}\cdot3^{2}\cdot5
 & 1 & 2_{5} & 2_{5} & 2_{5} \\\hline

 2_{6} & (x^{2}h^{2}xh^{2}xh)^{3} & 113400 & 2^{18}\cdot3
 & 1 & 2_{6} & 2_{6} & 2_{6} \\
\hline
 \end{array}
$$

 \newpage
 \normalsize
 Conjugacy classes of $H(\Co_1) = \langle x,y,h \rangle$
 (continued)
 \scriptsize
 $$
 \begin{array}{|c|r|r|r|c|c|c|c|
}\hline
 \mbox{Class} &
 \multicolumn{1}{c|}{\mbox{Representative}} &
 \multicolumn{1}{c|}{|\mbox{Class}|} &
 \multicolumn{1}{c|}{|\mbox{Centralizer}|}
 &
 2{\rm P}
 & 3{\rm P}
 & 5{\rm P}
 & 7{\rm P}
 \\\hline

 6_{22} & x^{2}h^{2}xh^{2}xh & 232243200 & 2^{7}\cdot3
 & 3_{4} & 2_{6} & 6_{22} & 6_{22} \\\hline

 6_{23} & xyhyxh & 928972800 & 2^{5}\cdot3
 & 3_{4} & 2_{8} & 6_{23} & 6_{23} \\\hline

 7 & x & 1592524800 & 2^{3}\cdot7
 & 7 & 7 & 7 & 1 \\\hline

 8_{1} & (x^{3}h)^{3} & 967680 & 2^{11}\cdot3^{2}\cdot5
 & 4_{1} & 8_{1} & 8_{1} & 8_{1} \\\hline

 8_{2} & yh & 14515200 & 2^{11}\cdot3
 & 4_{1} & 8_{2} & 8_{2} & 8_{2} \\\hline

 8_{3} & (x^{4}hxh)^{2} & 14515200 & 2^{11}\cdot3
 & 4_{3} & 8_{3} & 8_{3} & 8_{3} \\\hline

 8_{4} & (x^{2}yx^{2}hy)^{3} & 29030400 & 2^{10}\cdot3
 & 4_{6} & 8_{4} & 8_{4} & 8_{4} \\
\hline
 \end{array}
 $$

  \newpage
 \normalsize
 Conjugacy classes of $H(\Co_1) = \langle x,y,h \rangle$
 (continued)
 \scriptsize
 $$
 \begin{array}{|c|r|r|r|c|c|c|c|
}\hline
 \mbox{Class} &
 \multicolumn{1}{c|}{\mbox{Representative}} &
 \multicolumn{1}{c|}{|\mbox{Class}|} &
 \multicolumn{1}{c|}{|\mbox{Centralizer}|}
 &
 2{\rm P}
 & 3{\rm P}
 & 5{\rm P}
 & 7{\rm P}
 \\\hline

 12_{30} & x^{2}hxyxh & 1857945600 & 2^{4}\cdot3
 & 6_{12} & 4_{13} & 12_{30} & 12_{30} \\\hline

 12_{31} & x^{3}yhxh^{2} & 1857945600 & 2^{4}\cdot3
 & 6_{18} & 4_{19} & 12_{31} & 12_{31} \\\hline

 14_{1} & (xy^{3})^{2} & 1592524800 & 2^{3}\cdot7
 & 7 & 14_{1} & 14_{1} & 2_{1} \\\hline

 14_{2} & xy & 3185049600 & 2^{2}\cdot7
 & 7 & 14_{3} & 14_{3} & 2_{2} \\\hline

 14_{3} & xy^{2} & 3185049600 & 2^{2}\cdot7
 & 7 & 14_{2} & 14_{2} & 2_{2} \\\hline

 15_{1} & (x^{2}hy)^{2} & 743178240 & 2^{3}\cdot3\cdot5
 & 15_{1} & 5_{1} & 3_{1} & 15_{1} \\\hline

 15_{2} & (x^{2}y)^{2} & 2972712960 & 2\cdot3\cdot5
 & 15_{2} & 5_{2} & 3_{2} & 15_{2} \\
\hline

 \end{array}
 $$

}
\end{cclass}


\newpage
\normalsize

\begin{cclass}\label{Co_1 cc D} Conjugacy classes of $D(\Co_1) = \langle x,y \rangle$
{ \setlength{\arraycolsep}{1mm}
\renewcommand{\baselinestretch}{0.5}
\scriptsize
$$
 \begin{array}{|c|r|r|r|c|c|c|c|
}\hline
 \mbox{Class} &
 \multicolumn{1}{c|}{\mbox{Representative}} &
 \multicolumn{1}{c|}{|\mbox{Class}|} &
 \multicolumn{1}{c|}{|\mbox{Centralizer}|}
 &
 2{\rm P}
 & 3{\rm P}
 & 5{\rm P}
 & 7{\rm P}
 \\\hline
 1 & 1 & 1 & 2^{21}\cdot3^{2}\cdot5\cdot7
 & 1 & 1 & 1 & 1 \\\hline

 2_{1} & (x^{2}y)^{15} & 1 & 2^{21}\cdot3^{2}\cdot5\cdot7
 & 1 & 2_{1} & 2_{1} & 2_{1} \\\hline

 2_{2} & (y)^{2} & 30 & 2^{20}\cdot3\cdot7
 & 1 & 2_{2} & 2_{2} & 2_{2} \\\hline

 2_{3} & (xy)^{7} & 240 & 2^{17}\cdot3\cdot7
 & 1 & 2_{3} & 2_{3} & 2_{3} \\\hline

 2_{4} & (x^{5}yxyx^{2}y)^{4} & 280 & 2^{18}\cdot3^{2}
 & 1 & 2_{4} & 2_{4} & 2_{4} \\\hline

 2_{5} & (x^{4}yxyxy)^{5} & 448 & 2^{15}\cdot3^{2}\cdot5
 & 1 & 2_{5} & 2_{5} & 2_{5} \\\hline

 2_{6} & (x^{4}yxyxy^{3})^{5} & 448 & 2^{15}\cdot3^{2}\cdot5
 & 1 & 2_{6} & 2_{6} & 2_{6} \\
\hline
\end{array}
$$

 \newpage
 \normalsize
 Conjugacy classes of $D(\Co_1) = \langle x,y \rangle$
 (continued)
 \scriptsize

 $$
 \begin{array}{|c|r|r|r|c|c|c|c|
}\hline
 \mbox{Class} &
 \multicolumn{1}{c|}{\mbox{Representative}} &
 \multicolumn{1}{c|}{|\mbox{Class}|} &
 \multicolumn{1}{c|}{|\mbox{Centralizer}|}
 &
 2{\rm P}
 & 3{\rm P}
 & 5{\rm P}
 & 7{\rm P}
 \\\hline
 4_{42} & x^{4}yx^{4}yxyxy^{2}xy & 645120 & 2^{10}
 & 2_{9} & 4_{42} & 4_{42} & 4_{42} \\\hline

 4_{43} & x^{4}yx^{3}yxyxyx^{4}yxyxy & 645120 & 2^{10}
 & 2_{8} & 4_{43} & 4_{43} & 4_{43} \\\hline

 4_{44} & x^{3}yxyx^{2}y^{3}xyxy & 1290240 & 2^{9}
 & 2_{10} & 4_{44} & 4_{44} & 4_{44} \\\hline

 4_{45} & x^{4}yx^{4}yxyx^{2}y^{3} & 1290240 & 2^{9}
 & 2_{10} & 4_{45} & 4_{45} & 4_{45} \\\hline

 4_{46} & x^{4}yxyx^{3}y^{3} & 2580480 & 2^{8}
 & 2_{14} & 4_{46} & 4_{46} & 4_{46} \\\hline

 4_{47} & x^{4}yx^{2}yxyx^{3}y & 2580480 & 2^{8}
 & 2_{13} & 4_{47} & 4_{47} & 4_{47} \\\hline

 4_{48} & x^{4}yxyx^{3}y^{3}xy & 2580480 & 2^{8}
 & 2_{15} & 4_{48} & 4_{48} & 4_{48} \\
\hline
\end{array}
 $$

  \newpage
 \normalsize
 Conjugacy classes of $D(\Co_1) = \langle x,y \rangle$
 (continued)
 \scriptsize

 $$
 \begin{array}{|c|r|r|r|c|c|c|c|
}\hline
 \mbox{Class} &
 \multicolumn{1}{c|}{\mbox{Representative}} &
 \multicolumn{1}{c|}{|\mbox{Class}|} &
 \multicolumn{1}{c|}{|\mbox{Centralizer}|}
 &
 2{\rm P}
 & 3{\rm P}
 & 5{\rm P}
 & 7{\rm P}
 \\\hline
 10_{1} & (x^{2}y)^{3} & 5505024 & 2^{3}\cdot3\cdot5
 & 5 & 10_{1} & 2_{1} & 10_{1} \\\hline

 10_{2} & x^{4}yxyxy & 16515072 & 2^{3}\cdot5
 & 5 & 10_{2} & 2_{5} & 10_{2} \\\hline

 10_{3} & x^{4}yxyxy^{3} & 16515072 & 2^{3}\cdot5
 & 5 & 10_{3} & 2_{6} & 10_{3} \\\hline

 12_{1} & (x^{3}yx^{2}yx^{2}y)^{2} & 3440640 & 2^{6}\cdot3
 & 6_{2} & 4_{1} & 12_{1} & 12_{1} \\\hline

 12_{2} & x^{2}yxyx^{2}yxyxy^{3}xy & 3440640 & 2^{6}\cdot3
 & 6_{2} & 4_{2} & 12_{2} & 12_{2} \\\hline

 12_{3} & x^{3}yx^{2}y & 6881280 & 2^{5}\cdot3
 & 6_{2} & 4_{5} & 12_{3} & 12_{3} \\\hline

 12_{4} & x^{4}yx^{2}y & 6881280 & 2^{5}\cdot3
 & 6_{10} & 4_{6} & 12_{4} & 12_{4} \\
\hline

 \end{array}
 $$

}
\end{cclass}

\newpage
\normalsize

\begin{cclass}\label{Co_1 cc H1} Conjugacy classes of $H_1(\Co_1) = \langle j_1,k_1,h_1\rangle$\\
{ \setlength{\arraycolsep}{1mm}
\renewcommand{\baselinestretch}{0.5}
\scriptsize
$$
 \begin{array}{|c|r|r|c|c|
}\hline
 \mbox{Class} &
 \multicolumn{1}{c|}{\mbox{Representative}} &
 \multicolumn{1}{c|}{|\mbox{Centralizer}|}
 &
 2{\rm P}
 & 3{\rm P}
 \\\hline
 1 & 1 & 2^{21}\cdot3^{3}
 & 1 & 1 \\\hline

 2_{1} & (j_1)^{6} & 2^{21}\cdot3^{3}
 & 1 & 2_{1} \\\hline

 2_{2} & (j_1^{2}h_1)^{4} & 2^{20}\cdot3^{2}
 & 1 & 2_{2} \\\hline

 2_{3} & (j_1^{2}k_1^{2})^{6} & 2^{18}\cdot3^{3}
 & 1 & 2_{3} \\\hline

 2_{4} & (j_1k_1)^{4} & 2^{18}\cdot3
 & 1 & 2_{4} \\\hline

 2_{5} & (j_1^{2}k_1^{2}h_1j_1^{2}h_1)^{3} & 2^{18}\cdot3
 & 1 & 2_{5} \\\hline

 2_{6} & (k_1h_1)^{6} & 2^{15}\cdot3^{2}
 & 1 & 2_{6} \\
\hline
 \end{array}
$$

 \newpage
 \normalsize
 Conjugacy classes of $H_1(\Co_1) = \langle
j_1,k_1,h_1\rangle$(continued)
 \scriptsize
$$
 \begin{array}{|c|r|r|c|c|
}\hline
 \mbox{Class} &
 \multicolumn{1}{c|}{\mbox{Representative}} &
 \multicolumn{1}{c|}{|\mbox{Centralizer}|}
 &
 2{\rm P}
 & 3{\rm P}
 \\\hline
 4_{12} & (j_1^{3}k_1j_1h_1k_1^{2}h_1)^{2} & 2^{16}
 & 2_{1} & 4_{12} \\\hline

 4_{13} & (j_1^{2}k_1h_1)^{6} & 2^{14}\cdot3
 & 2_{1} & 4_{13} \\\hline

 4_{14} & (j_1k_1^{2}h_1)^{6} & 2^{14}\cdot3
 & 2_{1} & 4_{14} \\\hline

 4_{15} & (j_1^{2}k_1j_1h_1j_1h_1)^{3} & 2^{14}\cdot3
 & 2_{4} & 4_{15} \\\hline

 4_{16} & (j_1^{2}k_1j_1h_1k_1^{3})^{3} & 2^{12}\cdot3^{2}
 & 2_{1} & 4_{16} \\\hline

 4_{17} & (j_1^{4}h_1)^{2} & 2^{15}
 & 2_{2} & 4_{17} \\\hline

 4_{18} & (j_1^{3}k_1^{2}j_1h_1)^{2} & 2^{15}
 & 2_{1} & 4_{18} \\
\hline

\end{array}
$$

\newpage
\normalsize
Conjugacy classes of $H_1(\Co_1) = \langle
j_1,k_1,h_1\rangle$(continued)
\scriptsize
 $$
 \begin{array}{|c|r|r|c|c|
}\hline
 \mbox{Class} &
 \multicolumn{1}{c|}{\mbox{Representative}} &
 \multicolumn{1}{c|}{|\mbox{Centralizer}|}
 &
 2{\rm P}
 & 3{\rm P}
 \\\hline
 4_{65} & j_1^{2}k_1^{4}j_1k_1^{2}j_1h_1 & 2^{12}
 & 2_{10} & 4_{65} \\\hline

 4_{66} & j_1^{2}k_1h_1j_1k_1j_1k_1j_1k_1^{2} & 2^{12}
 & 2_{4} & 4_{66} \\\hline

 4_{67} & j_1^{2}k_1h_1k_1^{4}j_1^{2}h_1 & 2^{12}
 & 2_{4} & 4_{67} \\\hline

 4_{68} & j_1^{4}k_1j_1^{2}h_1j_1h_1k_1h_1 & 2^{12}
 & 2_{4} & 4_{68} \\\hline

 4_{69} & j_1^{4}k_1j_1k_1j_1k_1j_1^{2}k_1h_1 & 2^{12}
 & 2_{10} & 4_{69} \\\hline

 4_{70} & j_1^{4}h_1j_1k_1^{3}h_1k_1h_1k_1 & 2^{12}
 & 2_{4} & 4_{70} \\\hline

 4_{71} & j_1^{3}h_1^{2} & 2^{11}
 & 2_{4} & 4_{71} \\
\hline

\end{array}
$$

\newpage
\normalsize
Conjugacy classes of $H_1(\Co_1) = \langle
j_1,k_1,h_1\rangle$(continued)
\scriptsize
 $$
 \begin{array}{|c|r|r|c|c|
}\hline
 \mbox{Class} &
 \multicolumn{1}{c|}{\mbox{Representative}} &
 \multicolumn{1}{c|}{|\mbox{Centralizer}|}
 &
 2{\rm P}
 & 3{\rm P}
 \\\hline
 4_{118} & j_1^{4}h_1j_1h_1 & 2^{9}
 & 2_{17} & 4_{118} \\\hline

 4_{119} & j_1^{2}k_1j_1h_1j_1k_1 & 2^{9}
 & 2_{22} & 4_{119} \\\hline

 4_{120} & j_1^{2}k_1^{2}j_1k_1^{2} & 2^{9}
 & 2_{22} & 4_{120} \\\hline

 4_{121} & j_1k_1j_1k_1^{2}h_1k_1 & 2^{9}
 & 2_{12} & 4_{121} \\\hline

 4_{122} & j_1k_1j_1k_1h_1k_1^{2} & 2^{9}
 & 2_{12} & 4_{122} \\\hline

 4_{123} & j_1^{2}h_1k_1h_1j_1h_1k_1 & 2^{9}
 & 2_{17} & 4_{123} \\\hline

 4_{124} & j_1k_1h_1k_1^{2}h_1^{3} & 2^{9}
 & 2_{18} & 4_{124} \\
\hline

\end{array}
$$

\newpage
\normalsize
Conjugacy classes of $H_1(\Co_1) = \langle
j_1,k_1,h_1\rangle$(continued)
\scriptsize
 $$
 \begin{array}{|c|r|r|c|c|
}\hline
 \mbox{Class} &
 \multicolumn{1}{c|}{\mbox{Representative}} &
 \multicolumn{1}{c|}{|\mbox{Centralizer}|}
 &
 2{\rm P}
 & 3{\rm P}
 \\\hline
 6_{37} & j_1h_1k_1 & 2^{5}\cdot3^{2}
 & 3_{5} & 2_{7} \\\hline

 6_{38} & k_1^{3}h_1 & 2^{5}\cdot3^{2}
 & 3_{4} & 2_{9} \\\hline

 6_{39} & j_1^{3}h_1k_1 & 2^{5}\cdot3^{2}
 & 3_{3} & 2_{7} \\\hline

 6_{40} & j_1h_1k_1h_1^{2} & 2^{5}\cdot3^{2}
 & 3_{5} & 2_{7} \\\hline

 6_{41} & j_1k_1j_1h_1j_1h_1k_1h_1 & 2^{5}\cdot3^{2}
 & 3_{3} & 2_{8} \\\hline

 6_{42} & j_1^{3}k_1h_1k_1j_1^{2}k_1 & 2^{5}\cdot3^{2}
 & 3_{5} & 2_{16} \\\hline

 6_{43} & j_1^{2}h_1k_1h_1j_1h_1k_1^{2} & 2^{5}\cdot3^{2}
 & 3_{5} & 2_{15} \\
\hline

\end{array}
$$

\newpage
\normalsize
Conjugacy classes of $H_1(\Co_1) = \langle
j_1,k_1,h_1\rangle$(continued)
\scriptsize
 $$
 \begin{array}{|c|r|r|c|c|
}\hline
 \mbox{Class} &
 \multicolumn{1}{c|}{\mbox{Representative}} &
 \multicolumn{1}{c|}{|\mbox{Centralizer}|}
 &
 2{\rm P}
 & 3{\rm P}
 \\\hline
 8_{36} & j_1^{5}k_1h_1^{2} & 2^{9}
 & 4_{18} & 8_{36} \\\hline

 8_{37} & j_1^{4}h_1j_1h_1k_1^{2} & 2^{9}
 & 4_{18} & 8_{37} \\\hline

 8_{38} & j_1^{3}k_1h_1k_1^{4} & 2^{9}
 & 4_{12} & 8_{38} \\\hline

 8_{39} & j_1k_1h_1j_1h_1k_1^{4} & 2^{9}
 & 4_{17} & 8_{39} \\\hline

 8_{40} & j_1^{4}k_1h_1j_1h_1k_1h_1 & 2^{9}
 & 4_{12} & 8_{40} \\\hline

 8_{41} & j_1^{2}k_1j_1k_1h_1j_1h_1k_1h_1 & 2^{9}
 & 4_{18} & 8_{41} \\\hline

 8_{42} & j_1^{2}k_1^{3} & 2^{8}
 & 4_{27} & 8_{42} \\
\hline

\end{array}
$$

\newpage
\normalsize
Conjugacy classes of $H_1(\Co_1) = \langle
j_1,k_1,h_1\rangle$(continued)
\scriptsize
 $$
 \begin{array}{|c|r|r|c|c|
}\hline
 \mbox{Class} &
 \multicolumn{1}{c|}{\mbox{Representative}} &
 \multicolumn{1}{c|}{|\mbox{Centralizer}|}
 &
 2{\rm P}
 & 3{\rm P}
 \\\hline
 12_{23} & j_1^{3}k_1j_1h_1 & 2^{5}\cdot3^{2}
 & 6_{7} & 4_{55} \\\hline

 12_{24} & j_1^{2}k_1^{2}h_1^{2} & 2^{5}\cdot3^{2}
 & 6_{8} & 4_{8} \\\hline

 12_{25} & j_1k_1j_1k_1^{2}h_1 & 2^{5}\cdot3^{2}
 & 6_{21} & 4_{9} \\\hline

 12_{26} & j_1k_1j_1h_1k_1^{2} & 2^{5}\cdot3^{2}
 & 6_{6} & 4_{10} \\\hline

 12_{27} & j_1k_1j_1h_1k_1h_1 & 2^{5}\cdot3^{2}
 & 6_{14} & 4_{6} \\\hline

 12_{28} & j_1k_1^{2}h_1j_1h_1 & 2^{5}\cdot3^{2}
 & 6_{14} & 4_{8} \\\hline

 12_{29} & k_1^{2}h_1k_1h_1^{2} & 2^{5}\cdot3^{2}
 & 6_{21} & 4_{11} \\
\hline

\end{array}
$$

\newpage
\normalsize
Conjugacy classes of $H_1(\Co_1) = \langle
j_1,k_1,h_1\rangle$(continued)
\scriptsize
 $$
 \begin{array}{|c|r|r|c|c|
}\hline
 \mbox{Class} &
 \multicolumn{1}{c|}{\mbox{Representative}} &
 \multicolumn{1}{c|}{|\mbox{Centralizer}|}
 &
 2{\rm P}
 & 3{\rm P}
 \\\hline
 24_{9} & k_1 & 2^{4}\cdot3
 & 12_{17} & 8_{1} \\\hline

 24_{10} & j_1k_1^{2} & 2^{4}\cdot3
 & 12_{18} & 8_{12} \\\hline

 24_{11} & j_1^{3}h_1k_1h_1 & 2^{4}\cdot3
 & 12_{19} & 8_{13} \\\hline

 24_{12} & j_1^{4}h_1k_1h_1 & 2^{4}\cdot3
 & 12_{17} & 8_{4} \\\hline

 24_{13} & j_1k_1^{2}h_1 & 2^{3}\cdot3
 & 12_{46} & 8_{29} \\\hline
 \end{array}
 $$
}
\end{cclass}

\normalsize


\newpage
\begin{cclass}\label{Co_1 cc U}Conjugacy classes of $U(\Co_1) = \langle y,j,k\rangle$
{ \setlength{\arraycolsep}{1mm}
\renewcommand{\baselinestretch}{0.5}
\scriptsize
$$
\begin{array}{|c|r|r|r|c|c|
}\hline
 \mbox{Class} &
 \multicolumn{1}{c|}{\mbox{Representative}} &
 \multicolumn{1}{c|}{|\mbox{Centralizer}|}
 &
 2{\rm P}
 & 3{\rm P}
 \\\hline
 1 & 1 & 2^{21}\cdot3^{2}
 & 1 & 1 \\\hline

 2_{1} & (yj)^{6} & 2^{21}\cdot3^{2}
 & 1 & 2_{1} \\\hline

 2_{2} & (y^{3}j)^{6} & 2^{20}\cdot3
 & 1 & 2_{2} \\\hline

 2_{3} & (yj^{2})^{4} & 2^{18}\cdot3^{2}
 & 1 & 2_{3} \\\hline

 2_{4} & (y)^{2} & 2^{18}\cdot3
 & 1 & 2_{4} \\\hline

 2_{5} & (yjyj^{2}ykjyk)^{3} & 2^{18}\cdot3
 & 1 & 2_{5} \\\hline

 2_{6} & (yjyk)^{4} & 2^{17}\cdot3
 & 1 & 2_{6} \\
\hline
 \end{array}
 $$
 \normalsize

  \newpage
Conjugacy classes of $U(\Co_1) = \langle y,j,k\rangle$(continued)
\scriptsize
 $$
 \begin{array}{|c|r|r|r|c|c|
}\hline
 \mbox{Class} &
 \multicolumn{1}{c|}{\mbox{Representative}} &
 \multicolumn{1}{c|}{|\mbox{Centralizer}|}
 &
 2{\rm P}
 & 3{\rm P}
 \\\hline
 4_{5} & (yj^{3}kykykj)^{2} & 2^{16}
 & 2_{1} & 4_{5} \\\hline

 4_{6} & (yjyj^{2}kykjyk)^{2} & 2^{16}
 & 2_{2} & 4_{6} \\\hline

 4_{7} & (yj^{2}k)^{3} & 2^{14}\cdot3
 & 2_{4} & 4_{7} \\\hline

 4_{8} & (yjkyjk^{2})^{3} & 2^{14}\cdot3
 & 2_{2} & 4_{8} \\\hline

 4_{9} & (yjyjykyk)^{3} & 2^{14}\cdot3
 & 2_{3} & 4_{9} \\\hline

 4_{10} & (yjyj^{3}k^{2})^{3} & 2^{14}\cdot3
 & 2_{2} & 4_{10} \\\hline

 4_{11} & (yj^{3}kyjk)^{3} & 2^{14}\cdot3
 & 2_{3} & 4_{11} \\
\hline
 \end{array}
 $$
 \normalsize

  \newpage
Conjugacy classes of $U(\Co_1) = \langle y,j,k\rangle$(continued)
\scriptsize
 $$
 \begin{array}{|c|r|r|r|c|c|
}\hline
 \mbox{Class} &
 \multicolumn{1}{c|}{\mbox{Representative}} &
 \multicolumn{1}{c|}{|\mbox{Centralizer}|}
 &
 2{\rm P}
 & 3{\rm P}
 \\\hline
 4_{58} & yjyjyj^{3}ykj^{2}k & 2^{13}
 & 2_{1} & 4_{58} \\\hline

 4_{59} & yjyj^{2}yj^{3}k^{2}yk & 2^{13}
 & 2_{2} & 4_{59} \\\hline

 4_{60} & yjyj^{2}yjkjkjyk & 2^{13}
 & 2_{1} & 4_{60} \\\hline

 4_{61} & yjyj^{2}kjyjkjk^{2} & 2^{13}
 & 2_{9} & 4_{61} \\\hline

 4_{62} & yjyjkj^{2}k^{2}ykj^{2} & 2^{13}
 & 2_{9} & 4_{62} \\\hline

 4_{63} & yj^{2}yj^{2}yjkjkyk & 2^{13}
 & 2_{2} & 4_{63} \\\hline

 4_{64} & yj^{2}kyj^{2}kyjkjk^{2} & 2^{13}
 & 2_{1} & 4_{64} \\
\hline
 \end{array}
 $$

 \normalsize

  \newpage
Conjugacy classes of $U(\Co_1) = \langle y,j,k\rangle$(continued)
\scriptsize
 $$
 \begin{array}{|c|r|r|r|c|c|
}\hline
 \mbox{Class} &
 \multicolumn{1}{c|}{\mbox{Representative}} &
 \multicolumn{1}{c|}{|\mbox{Centralizer}|}
 &
 2{\rm P}
 & 3{\rm P}
 \\\hline
 4_{111} & yjyjkyj^{3}k^{3} & 2^{11}
 & 2_{9} & 4_{111} \\\hline

 4_{112} & yjyjkyjkj^{2}k^{2} & 2^{11}
 & 2_{9} & 4_{112} \\\hline

 4_{113} & yjyjk^{3}yj^{3}k & 2^{11}
 & 2_{9} & 4_{113} \\\hline

 4_{114} & yjyjk^{3}yjkj^{2} & 2^{11}
 & 2_{9} & 4_{114} \\\hline

 4_{115} & yj^{3}kyk^{3}yjk & 2^{11}
 & 2_{9} & 4_{115} \\\hline

 4_{116} & y^{3}jykyjkyj^{2}k & 2^{11}
 & 2_{9} & 4_{116} \\\hline

 4_{117} & y^{3}jykyjkykj^{2} & 2^{11}
 & 2_{9} & 4_{117} \\
\hline
 \end{array}
 $$

 \normalsize

  \newpage
Conjugacy classes of $U(\Co_1) = \langle y,j,k\rangle$(continued)
\scriptsize
 $$
 \begin{array}{|c|r|r|r|c|c|
}\hline
 \mbox{Class} &
 \multicolumn{1}{c|}{\mbox{Representative}} &
 \multicolumn{1}{c|}{|\mbox{Centralizer}|}
 &
 2{\rm P}
 & 3{\rm P}
 \\\hline
 4_{164} & yjyjkykjyk & 2^{9}
 & 2_{12} & 4_{164} \\\hline

 4_{165} & yj^{2}ykykykj & 2^{9}
 & 2_{12} & 4_{165} \\\hline

 4_{166} & y^{3}jykykjk^{2} & 2^{9}
 & 2_{12} & 4_{166} \\\hline

 4_{167} & yjykyjkyjkj & 2^{9}
 & 2_{10} & 4_{167} \\\hline

 4_{168} & yj^{2}yj^{3}kyjk & 2^{9}
 & 2_{27} & 4_{168} \\\hline

 4_{169} & yjyjyj^{2}kjkyk & 2^{9}
 & 2_{19} & 4_{169} \\\hline

 4_{170} & yjyj^{3}yj^{2}yjk & 2^{9}
 & 2_{11} & 4_{170} \\
\hline
  \end{array}
 $$
 \normalsize

  \newpage
Conjugacy classes of $U(\Co_1) = \langle y,j,k\rangle$(continued)
\scriptsize
 $$
 \begin{array}{|c|r|r|r|c|c|
}\hline
 \mbox{Class} &
 \multicolumn{1}{c|}{\mbox{Representative}} &
 \multicolumn{1}{c|}{|\mbox{Centralizer}|}
 &
 2{\rm P}
 & 3{\rm P}
 \\\hline
 6_{29} & yj^{3}kj & 2^{5}\cdot3
 & 3_{2} & 2_{35} \\\hline

 6_{30} & yjkykjk & 2^{5}\cdot3
 & 3_{1} & 2_{22} \\\hline

 6_{31} & yjykyk^{3} & 2^{5}\cdot3
 & 3_{1} & 2_{36} \\\hline

 6_{32} & yjyjyj^{2}kyk & 2^{5}\cdot3
 & 3_{3} & 2_{24} \\\hline

 6_{33} & yjyjyjkjyk & 2^{5}\cdot3
 & 3_{3} & 2_{21} \\\hline

 8_{1} & (yk)^{3} & 2^{11}\cdot3
 & 4_{1} & 8_{1} \\\hline

 8_{2} & (yjkj)^{3} & 2^{11}\cdot3
 & 4_{1} & 8_{2} \\
\hline
 \end{array}
 $$
 \normalsize

  \newpage
Conjugacy classes of $U(\Co_1) = \langle y,j,k\rangle$(continued)
\scriptsize
 $$
 \begin{array}{|c|r|r|r|c|c|
}\hline
 \mbox{Class} &
 \multicolumn{1}{c|}{\mbox{Representative}} &
 \multicolumn{1}{c|}{|\mbox{Centralizer}|}
 &
 2{\rm P}
 & 3{\rm P}
 \\\hline
 8_{49} & yj^{2}ykyk & 2^{8}
 & 4_{22} & 8_{49} \\\hline

 8_{50} & yjyjyj^{3} & 2^{8}
 & 4_{21} & 8_{50} \\\hline

 8_{51} & yjkyjkyk & 2^{8}
 & 4_{25} & 8_{51} \\\hline

 8_{52} & y^{3}jyjyj^{2} & 2^{8}
 & 4_{46} & 8_{52} \\\hline

 8_{53} & yjyj^{2}yj^{3} & 2^{8}
 & 4_{21} & 8_{53} \\\hline

 8_{54} & yjykyjk^{3} & 2^{8}
 & 4_{26} & 8_{54} \\\hline

 8_{55} & yj^{2}k^{2}yjk^{2} & 2^{8}
 & 4_{27} & 8_{55} \\
\hline
\end{array}
 $$
 \normalsize

\newpage
Conjugacy classes of $U(\Co_1) = \langle y,j,k\rangle$(continued)
\scriptsize
 $$
 \begin{array}{|c|r|r|r|c|c|
}\hline
 \mbox{Class} &
 \multicolumn{1}{c|}{\mbox{Representative}} &
 \multicolumn{1}{c|}{|\mbox{Centralizer}|}
 &
 2{\rm P}
 & 3{\rm P}
 \\\hline

 12_{24} & yjk^{2} & 2^{4}\cdot3
 & 6_{23} & 4_{131} \\\hline

 12_{25} & y^{3}jk & 2^{4}\cdot3
 & 6_{22} & 4_{132} \\\hline

 12_{26} & yjyjk & 2^{4}\cdot3
 & 6_{23} & 4_{129} \\\hline

 12_{27} & ykykj & 2^{4}\cdot3
 & 6_{24} & 4_{65} \\\hline

 12_{28} & yj^{2}kyjk & 2^{4}\cdot3
 & 6_{24} & 4_{66} \\\hline

 16_{1} & yj^{2}k^{2} & 2^{6}
 & 8_{3} & 16_{1} \\\hline

 16_{2} & yjyjkyk & 2^{6}
 & 8_{4} & 16_{2} \\\hline

 24_{1} & yk & 2^{4}\cdot3
 & 12_{1} & 8_{1} \\\hline

 24_{2} & yjkj & 2^{4}\cdot3
 & 12_{1} & 8_{2} \\\hline

 24_{3} & yjykyk & 2^{4}\cdot3
 & 12_{2} & 8_{9} \\\hline

 24_{4} & yjyj^{3}k & 2^{4}\cdot3
 & 12_{2} & 8_{10} \\\hline
 \end{array}
 $$
}
\end{cclass}

\setlength\textheight{50\baselineskip}
\addtolength\textheight{\topskip}

\newpage

\section{Character tables}

\begin{ctab}\label{Co_1ct_E} Character table of $E(\Co_1) = \langle x,y,e \rangle$
{ \setlength{\arraycolsep}{1mm}
\renewcommand{\baselinestretch}{0.5}

{\tiny
\setlength{\arraycolsep}{0,3mm}
$$
\begin{array}{r|rrrrrrrrrrrrrrrrrrr}
2&21&21&18&17&15&15&17&8&6&17&17&14&14&15&15&10&11&11\\
3&3&2&3&1&1&1&1&3&2&1&1&1&1&.&.&1&1&1\\
5&1&1&1&.&1&1&.&1&.&.&.&.&.&.&.&1&.&.\\
7&1&1&.&1&.&.&.&.&1&1&.&.&.&.&.&.&.&.\\
11&1&.&1&.&.&.&.&.&.&.&.&.&.&.&.&.&.&.\\
23&1&.&.&.&.&.&.&.&.&.&.&.&.&.&.&.&.&.\\
\hline
&1a&2a&2b&2c&2d&2e&2f&3a&3b&4a&4b&4c&4d&4e&4f&4g&4h&4i\\
\hline
2P&1a&1a&1a&1a&1a&1a&1a&3a&3b&2a&2a&2a&2a&2a&2a&2b&2c&2c\\
\end{array}
$$
} }

\newpage
Character table of $E(\Co_1) = \langle x,y,e\rangle$ (continued){
\renewcommand{\baselinestretch}{0.5}

{\tiny
\setlength{\arraycolsep}{0,3mm}
$$
\begin{array}{r|rrrrrrrrrrrrrrrrrrrrrrrrrrr}
2&11&11&10&10&8&8&5&7&7&8&6&6&6&5&5&3&3&11&11&11&11&10&10&8&8&8\\
3&.&.&.&.&1&1&1&2&3&2&2&1&1&1&1&1&1&1&1&.&.&.&.&.&.&.\\
5&.&.&.&.&.&.&1&1&.&.&.&.&.&.&.&.&.&.&.&.&.&.&.&.&.&.\\
7&.&.&.&.&.&.&.&.&.&.&.&.&.&.&.&1&1&.&.&.&.&.&.&.&.&.\\
11&.&.&.&.&.&.&.&.&.&.&.&.&.&.&.&.&.&.&.&.&.&.&.&.&.&.\\
23&.&.&.&.&.&.&.&.&.&.&.&.&.&.&.&.&.&.&.&.&.&.&.&.&.&.\\
\hline
&4j&4k&4l&4m&4n&4o&5a&6a&6b&6c&6d&6e&6f&6g&6h&7a&7b&8a&8b&8c&8d
&8e&8f&8g&8h&8i\\
\hline
2P&2c&2c&2f&2f&2d&2e&5a&3a&3a&3a&3b&3a&3a&3b&3b&7a&7b&4a&4a&4b&4b
&4a&4b&4e&4c&4f\\
\end{array}
$$
} }

\newpage
Character table of $E(\Co_1 = \langle x,y,e\rangle)$ (continued){
\renewcommand{\baselinestretch}{0.5}

{\tiny
\setlength{\arraycolsep}{0,3mm}
$$
\begin{array}{r|rrrrrrrrrrrrrrrrrrrrrrrrrr}
2&6&6&5&3&4&4&1&6&6&5&5&4&4&4&3&3&3&3&2&2&1&1&6&6&3\\
3&.&.&.&1&.&.&.&1&1&1&1&1&1&1&1&1&.&.&.&.&1&1&.&.&.\\
5&.&.&1&1&1&1&.&.&.&.&.&.&.&.&.&.&.&.&.&.&1&1&.&.&1\\
7&.&.&.&.&.&.&.&.&.&.&.&.&.&.&.&.&1&1&1&1&.&.&.&.&.\\
11&.&.&.&.&.&.&1&.&.&.&.&.&.&.&.&.&.&.&.&.&.&.&.&.&.\\
23&.&.&.&.&.&.&.&.&.&.&.&.&.&.&.&.&.&.&.&.&.&.&.&.&.\\
\hline
&8j&8k&10a&10b&10c&10d&11a&12a&12b&12c&12d&12e&12f&12g&12h&12i
&14a&14b&14c&14d&15a&15b&16a&16b&20a\\
\hline
2P&4j&4k&5a&5a&5a&5a&11a&6c&6c&6c&6c&6d&6e&6e&6g&6h&7a&7b&7a&7b
&15a&15b&8c&8d&10a\\
\end{array}
$$
} }

\newpage
Character table of $E(\Co_1)=\langle x,y,e\rangle$ (continued){
\renewcommand{\baselinestretch}{0.5}

{\tiny
\setlength{\arraycolsep}{0,3mm}
$$
\begin{array}{r|rrrrrrrrrrrr}
2&.&.&1&.&.&4&4&2&2&1&1\\
3&1&1&.&.&.&1&1&.&.&1&1\\
5&.&.&.&.&.&.&.&.&.&1&1\\
7&1&1&.&.&.&.&.&1&1&.&.\\
11&.&.&1&.&.&.&.&.&.&.&.\\
23&.&.&.&1&1&.&.&.&.&.&.\\
\hline
&21a&21b&22a&23a&23b&24a&24b&28a&28b&30a&30b\\
\hline
2P&21a&21b&11a&23a&23b&12a&12a&14a&14b&15a&15b\\
\end{array}
$$
} }
$
A = \zeta(7)^4 + \zeta(7)^2 + \zeta(7),
B = 2*\zeta(15)_3*\zeta(15)_5^3 + 2*\zeta(15)_3*\zeta(15)_5^2 + \zeta(15)_3 +
\zeta(15)_5^3 + \zeta(15)_5^2,
C = \zeta(23)^{18} + \zeta(23)^{16} + \zeta(23)^{13} + \zeta(23)^{12} +
\zeta(23)^9 + \zeta(23)^8 + \zeta(23)^6 + \zeta(23)^4 + \zeta(23)^3 +
\zeta(23)^2 + \zeta(23),
D = -2*\zeta(7)^4 - 2*\zeta(7)^2 - 2*\zeta(7) - 2,
E = -3*\zeta(7)^4 - 3*\zeta(7)^2 - 3*\zeta(7) - 3
$.
\end{ctab}

\newpage

\begin{ctab}\label{Co_1ct_H_1}Character table of $H_1(\Co_1) = \langle
j_1,k_1,h_1\rangle$
{ \setlength{\arraycolsep}{1mm}
\renewcommand{\baselinestretch}{0.5}

{\tiny \setlength{\arraycolsep}{0,3mm}
$$
\begin{array}{r|rrrrrrrrrrrrrrrrrrrrrrrrrrr}
2&21&21&20&18&18&18&15&14&14&14&17&17&15&15&15&13&13&16&16&14&14
&14&15&13&13&13\\
3&3&3&2&3&1&1&2&2&2&2&.&.&1&1&1&2&2&.&.&1&1&1&.&1&1&1\\
\hline
&1a&2_{1}&2_{2}&2_{3}&2_{4}&2_{5}&2_{6}&2_{7}&2_{8}&2_{9}&2_{10}
&2_{11}&2_{12}&2_{13}&2_{14}&2_{15}&2_{16}&2_{17}&2_{18}&2_{19}&2_{20}
&2_{21}&2_{22}&2_{23}&2_{24}&2_{25}\\
\hline
2P&1a&1a&1a&1a&1a&1a&1a&1a&1a&1a&1a&1a&1a&1a&1a&1a&1a&1a&1a&1a&1a
&1a&1a&1a&1a&1a\\
3P&1a&2_{1}&2_{2}&2_{3}&2_{4}&2_{5}&2_{6}&2_{7}&2_{8}&2_{9}&2_{10}
&2_{11}&2_{12}&2_{13}&2_{14}&2_{15}&2_{16}&2_{17}&2_{18}&2_{19}&2_{20}
&2_{21}&2_{22}&2_{23}&2_{24}&2_{25}\\
\hline
X.1&1&1&1&1&1&1&1&1&1&1&1&1&1&1&1&1&1&1&1&1&1&1&1&1&1&1\\
X.2&1&1&1&1&1&1&1&-1&1&1&1&1&1&1&1&-1&-1&1&1&1&1&-1&1&1&1&1\\
X.3&1&1&1&1&1&1&1&-1&1&-1&1&1&1&1&1&-1&-1&1&1&1&1&-1&1&1&-1&-1\\
\end{array}
$$
} }

\newpage
Character table of $H_1(\Co_1) = \langle
j_1,k_1,h_1\rangle$(continued){
\renewcommand{\baselinestretch}{0.5}

{\tiny \setlength{\arraycolsep}{0,3mm}
$$
\begin{array}{r|rrrrrrrrrrrrrrrrrrrrrrrrrrr}
2&13&13&14&14&12&12&12&13&11&13&8&8&8&7&7&6&17&17&15&16&16&14&14
&14&14&14\\
3&1&1&.&.&1&1&1&.&.&3&3&3&3&3&3&3&2&2&2&1&1&2&2&2&2&2\\
\hline
&2_{26}&2_{27}&2_{28}&2_{29}&2_{30}&2_{31}&2_{32}&2_{33}&2_{34}
&3a&3b&3c&3d&3e&3f&3g&4_{1}&4_{2}&4_{3}&4_{4}&4_{5}&4_{6}&4_{7}
&4_{8}&4_{9}&4_{10}\\
\hline
2P&1a&1a&1a&1a&1a&1a&1a&1a&1a&3a&3b&3c&3d&3e&3f&3g&2_{1}&2_{1}
&2_{1}&2_{2}&2_{2}&2_{3}&2_{2}&2_{3}&2_{2}&2_{1}\\
3P&2_{26}&2_{27}&2_{28}&2_{29}&2_{30}&2_{31}&2_{32}&2_{33}&2_{34}
&1a&1a&1a&1a&1a&1a&1a&4_{1}&4_{2}&4_{3}&4_{4}&4_{5}&4_{6}&4_{7}
&4_{8}&4_{9}&4_{10}\\
\hline
X.1&1&1&1&1&1&1&1&1&1&1&1&1&1&1&1&1&1&1&1&1&1&1&1&1&1&1\\
X.2&1&-1&1&1&-1&1&-1&-1&-1&1&1&1&1&1&1&1&1&1&1&1&1&1&1&1&1&1\\
X.3&-1&-1&1&-1&-1&-1&1&-1&1&1&1&1&1&1&1&1&1&1&1&1&1&1&1&1&-1&-1\\
\end{array}
$$
} }

\newpage
Character table of $H_1(\Co_1) = \langle
j_1,k_1,h_1\rangle$(continued){
\renewcommand{\baselinestretch}{0.5}

{\tiny \setlength{\arraycolsep}{0,3mm}
$$
\begin{array}{r|rrrrrrrrrrrrrrrrrrrrrrrrrr}
2&14&16&14&14&14&12&15&15&15&13&13&13&13&13&13&13&14&14&14&14&14
&14&14&14&14\\
3&2&.&1&1&1&2&.&.&.&1&1&1&1&1&1&1&.&.&.&.&.&.&.&.&.\\
\hline
&4_{11}&4_{12}&4_{13}&4_{14}&4_{15}&4_{16}&4_{17}&4_{18}&4_{19}
&4_{20}&4_{21}&4_{22}&4_{23}&4_{24}&4_{25}&4_{26}&4_{27}&4_{28}&4_{29}
&4_{30}&4_{31}&4_{32}&4_{33}&4_{34}&4_{35}\\
\hline
2P&2_{2}&2_{1}&2_{1}&2_{1}&2_{4}&2_{1}&2_{2}&2_{1}&2_{1}&2_{1}
&2_{2}&2_{1}&2_{3}&2_{2}&2_{1}&2_{2}&2_{2}&2_{2}&2_{2}&2_{1}&2_{2}
&2_{4}&2_{2}&2_{2}&2_{2}\\
3P&4_{11}&4_{12}&4_{13}&4_{14}&4_{15}&4_{16}&4_{17}&4_{18}&4_{19}
&4_{20}&4_{21}&4_{22}&4_{23}&4_{24}&4_{25}&4_{26}&4_{27}&4_{28}&4_{29}
&4_{30}&4_{31}&4_{32}&4_{33}&4_{34}&4_{35}\\
\hline
X.1&1&1&1&1&1&1&1&1&1&1&1&1&1&1&1&1&1&1&1&1&1&1&1&1&1\\
X.2&1&1&1&1&-1&1&1&1&1&1&1&1&1&1&-1&1&1&1&1&1&1&-1&1&1&1\\
X.3&-1&1&1&1&-1&1&1&1&1&1&-1&-1&1&1&-1&-1&1&1&1&-1&1&-1&1&-1&-1\\
\end{array}
$$
} }

\newpage
Character table of $H_1(\Co_1) = \langle
j_1,k_1,h_1\rangle$(continued){
\renewcommand{\baselinestretch}{0.5}

{\tiny \setlength{\arraycolsep}{0,3mm}
$$
\begin{array}{r|rrrrrrrrrrrrrrrrrrrrrrrrrr}
2&12&12&12&12&12&12&12&12&12&12&12&12&12&13&13&13&13&13&13&9&12
&12&12&12&12\\
3&1&1&1&1&1&1&1&1&1&1&1&1&1&.&.&.&.&.&.&2&.&.&.&.&.\\
\hline
&4_{36}&4_{37}&4_{38}&4_{39}&4_{40}&4_{41}&4_{42}&4_{43}&4_{44}
&4_{45}&4_{46}&4_{47}&4_{48}&4_{49}&4_{50}&4_{51}&4_{52}&4_{53}&4_{54}
&4_{55}&4_{56}&4_{57}&4_{58}&4_{59}&4_{60}\\
\hline
2P&2_{1}&2_{1}&2_{4}&2_{4}&2_{2}&2_{4}&2_{2}&2_{2}&2_{2}&2_{2}
&2_{2}&2_{4}&2_{4}&2_{2}&2_{3}&2_{3}&2_{10}&2_{1}&2_{10}&2_{6}&2_{4}
&2_{4}&2_{4}&2_{2}&2_{10}\\
3P&4_{36}&4_{37}&4_{38}&4_{39}&4_{40}&4_{41}&4_{42}&4_{43}&4_{44}
&4_{45}&4_{46}&4_{47}&4_{48}&4_{49}&4_{50}&4_{51}&4_{52}&4_{53}&4_{54}
&4_{55}&4_{56}&4_{57}&4_{58}&4_{59}&4_{60}\\
\hline
X.1&1&1&1&1&1&1&1&1&1&1&1&1&1&1&1&1&1&1&1&1&1&1&1&1&1\\
X.2&-1&-1&-1&-1&1&-1&-1&-1&-1&1&1&1&1&1&1&1&-1&-1&-1&1&1&1&1&-1&-1\\
X.3&1&-1&-1&-1&-1&-1&1&-1&1&1&-1&1&1&1&1&1&1&-1&1&-1&1&1&1&-1&1\\
\end{array}
$$
} }

\newpage
Character table of $H_1(\Co_1) = \langle
j_1,k_1,h_1\rangle$(continued){
\renewcommand{\baselinestretch}{0.5}

{\tiny \setlength{\arraycolsep}{0,3mm}
$$
\begin{array}{r|rrrrrrrrrrrrrrrrrrrrrrrrrr}
2&12&12&12&12&12&12&12&12&12&12&11&11&11&11&11&11&11&11&11&11&11
&11&11&11&11\\
3&.&.&.&.&.&.&.&.&.&.&.&.&.&.&.&.&.&.&.&.&.&.&.&.&.\\
\hline
&4_{61}&4_{62}&4_{63}&4_{64}&4_{65}&4_{66}&4_{67}&4_{68}&4_{69}
&4_{70}&4_{71}&4_{72}&4_{73}&4_{74}&4_{75}&4_{76}&4_{77}&4_{78}&4_{79}
&4_{80}&4_{81}&4_{82}&4_{83}&4_{84}&4_{85}\\
\hline
2P&2_{4}&2_{3}&2_{4}&2_{4}&2_{10}&2_{4}&2_{4}&2_{4}&2_{10}&2_{4}
&2_{4}&2_{10}&2_{10}&2_{10}&2_{10}&2_{10}&2_{2}&2_{10}&2_{10}&2_{10}
&2_{10}&2_{10}&2_{10}&2_{4}&2_{10}\\
3P&4_{61}&4_{62}&4_{63}&4_{64}&4_{65}&4_{66}&4_{67}&4_{68}&4_{69}
&4_{70}&4_{71}&4_{72}&4_{73}&4_{74}&4_{75}&4_{76}&4_{77}&4_{78}&4_{79}
&4_{80}&4_{81}&4_{82}&4_{83}&4_{84}&4_{85}\\
\hline
X.1&1&1&1&1&1&1&1&1&1&1&1&1&1&1&1&1&1&1&1&1&1&1&1&1&1\\
X.2&1&1&1&-1&-1&-1&-1&1&-1&-1&-1&1&1&1&-1&-1&-1&1&1&1&1&1&-1&-1&-1\\
X.3&1&1&1&-1&1&-1&-1&1&1&-1&1&-1&-1&1&-1&-1&1&1&-1&-1&-1&1&1&1&-1\\
\end{array}
$$
} }

\newpage
Character table of $H_1(\Co_1) = \langle
j_1,k_1,h_1\rangle$(continued){
\renewcommand{\baselinestretch}{0.5}

{\tiny \setlength{\arraycolsep}{0,3mm}
$$
\begin{array}{r|rrrrrrrrrrrrrrrrrrrrrrrr}
2&11&11&11&11&11&11&11&11&11&11&11&11&11&11&11&11&11&11&9&9&9&9
&10\\
3&.&.&.&.&.&.&.&.&.&.&.&.&.&.&.&.&.&.&1&1&1&1&.\\
\hline
&4_{86}&4_{87}&4_{88}&4_{89}&4_{90}&4_{91}&4_{92}&4_{93}&4_{94}
&4_{95}&4_{96}&4_{97}&4_{98}&4_{99}&4_{100}&4_{101}&4_{102}&4_{103}
&4_{104}&4_{105}&4_{106}&4_{107}&4_{108}\\
\hline
2P&2_{10}&2_{10}&2_{10}&2_{10}&2_{10}&2_{10}&2_{4}&2_{4}&2_{4}
&2_{10}&2_{10}&2_{10}&2_{10}&2_{10}&2_{10}&2_{10}&2_{4}&2_{10}&2_{12}
&2_{6}&2_{12}&2_{6}&2_{17}\\
3P&4_{86}&4_{87}&4_{88}&4_{89}&4_{90}&4_{91}&4_{92}&4_{93}&4_{94}
&4_{95}&4_{96}&4_{97}&4_{98}&4_{99}&4_{100}&4_{101}&4_{102}&4_{103}
&4_{104}&4_{105}&4_{106}&4_{107}&4_{108}\\
\hline
X.1&1&1&1&1&1&1&1&1&1&1&1&1&1&1&1&1&1&1&1&1&1&1&1\\
X.2&-1&-1&-1&-1&-1&-1&-1&1&1&1&1&1&1&-1&-1&-1&-1&-1&-1&1&-1&-1&-1\\
X.3&1&1&-1&-1&1&-1&1&-1&-1&-1&-1&-1&1&1&-1&-1&1&1&-1&-1&-1&1&-1\\
\end{array}
$$
} }

\newpage
Character table of $H_1(\Co_1) = \langle
j_1,k_1,h_1\rangle$(continued){
\renewcommand{\baselinestretch}{0.5}

{\tiny \setlength{\arraycolsep}{0,3mm}
$$
\begin{array}{r|rrrrrrrrrrrrrrrrrrrrr}
2&10&10&10&10&8&9&9&9&9&9&9&9&9&9&9&9&9&9&9&9\\
3&.&.&.&.&1&.&.&.&.&.&.&.&.&.&.&.&.&.&.&.\\
\hline
&4_{109}&4_{110}&4_{111}&4_{112}&4_{113}&4_{114}&4_{115}&4_{116}
&4_{117}&4_{118}&4_{119}&4_{120}&4_{121}&4_{122}&4_{123}&4_{124}&4_{125}
&4_{126}&4_{127}&4_{128}\\
\hline
2P&2_{18}&2_{5}&2_{18}&2_{17}&2_{20}&2_{12}&2_{6}&2_{22}&2_{12}
&2_{17}&2_{22}&2_{22}&2_{12}&2_{12}&2_{17}&2_{18}&2_{17}&2_{18}&2_{22}
&2_{11}\\
3P&4_{109}&4_{110}&4_{111}&4_{112}&4_{113}&4_{114}&4_{115}&4_{116}
&4_{117}&4_{118}&4_{119}&4_{120}&4_{121}&4_{122}&4_{123}&4_{124}&4_{125}
&4_{126}&4_{127}&4_{128}\\
\hline
X.1&1&1&1&1&1&1&1&1&1&1&1&1&1&1&1&1&1&1&1&1\\
X.2&1&1&1&-1&1&1&-1&1&1&-1&-1&-1&-1&-1&1&1&-1&-1&1&1\\
X.3&1&1&1&-1&-1&-1&1&-1&-1&1&1&1&1&1&1&-1&-1&-1&-1&-1\\
\end{array}
$$
} }

\newpage
Character table of $H_1(\Co_1) = \langle
j_1,k_1,h_1\rangle$(continued){
\renewcommand{\baselinestretch}{0.5}

{\tiny \setlength{\arraycolsep}{0,3mm}
$$
\begin{array}{r|rrrrrrrrrrrrrrrrrrrrrrrrrr}
2&9&8&8&8&8&8&13&12&10&8&8&8&9&7&7&7&7&7&8&6&6&6&6&9&7\\
3&.&.&.&.&.&.&3&2&3&3&3&3&2&3&3&3&3&3&2&3&3&3&3&1&2\\
\hline
&4_{129}&4_{130}&4_{131}&4_{132}&4_{133}&4_{134}&6_{1}&6_{2}
&6_{3}&6_{4}&6_{5}&6_{6}&6_{7}&6_{8}&6_{9}&6_{10}&6_{11}&6_{12}&6_{13}
&6_{14}&6_{15}&6_{16}&6_{17}&6_{18}&6_{19}\\
\hline
2P&2_{18}&2_{28}&2_{20}&2_{14}&2_{28}&2_{28}&3a&3a&3a&3b&3c&3d&3a
&3e&3f&3e&3e&3c&3a&3f&3g&3g&3g&3a&3d\\
3P&4_{129}&4_{130}&4_{131}&4_{132}&4_{133}&4_{134}&2_{1}&2_{2}
&2_{3}&2_{1}&2_{1}&2_{1}&2_{6}&2_{3}&2_{1}&2_{1}&2_{3}&2_{3}&2_{9}
&2_{3}&2_{3}&2_{3}&2_{1}&2_{20}&2_{8}\\
\hline
X.1&1&1&1&1&1&1&1&1&1&1&1&1&1&1&1&1&1&1&1&1&1&1&1&1&1\\
X.2&-1&1&-1&-1&-1&-1&1&1&1&1&1&1&1&1&1&1&1&1&1&1&1&1&1&1&1\\
X.3&-1&-1&1&-1&1&-1&1&1&1&1&1&1&1&1&1&1&1&1&-1&1&1&1&1&1&1\\
\end{array}
$$
} }

\newpage
Character table of $H_1(\Co_1) = \langle
j_1,k_1,h_1\rangle$(continued){
\renewcommand{\baselinestretch}{0.5}

{\tiny \setlength{\arraycolsep}{0,3mm}
$$
\begin{array}{r|rrrrrrrrrrrrrrrrrrrrrrrrrrr}
2&7&7&7&5&5&8&8&6&6&6&6&6&6&7&7&7&7&5&5&5&5&5&5&5&5&6\\
3&2&2&2&3&3&1&1&2&2&2&2&2&2&1&1&1&1&2&2&2&2&2&2&2&2&1\\
\hline
&6_{20}&6_{21}&6_{22}&6_{23}&6_{24}&6_{25}&6_{26}&6_{27}&6_{28}
&6_{29}&6_{30}&6_{31}&6_{32}&6_{33}&6_{34}&6_{35}&6_{36}&6_{37}&6_{38}
&6_{39}&6_{40}&6_{41}&6_{42}&6_{43}&6_{44}&6_{45}\\
\hline
2P&3b&3d&3b&3d&3b&3a&3a&3d&3b&3b&3b&3c&3c&3d&3c&3b&3c&3e&3d&3c&3e
&3c&3e&3e&3b&3e\\
3P&2_{2}&2_{2}&2_{7}&2_{3}&2_{3}&2_{14}&2_{24}&2_{6}&2_{6}&2_{15}
&2_{16}&2_{16}&2_{15}&2_{19}&2_{4}&2_{21}&2_{5}&2_{7}&2_{9}&2_{7}&2_{7}
&2_{8}&2_{16}&2_{15}&2_{9}&2_{12}\\
\hline
X.1&1&1&1&1&1&1&1&1&1&1&1&1&1&1&1&1&1&1&1&1&1&1&1&1&1&1\\
X.2&1&1&-1&1&1&1&1&1&1&-1&-1&-1&-1&1&1&-1&1&-1&1&-1&-1&1&-1&-1&1&1\\
X.3&1&1&-1&1&1&1&-1&1&1&-1&-1&-1&-1&1&1&-1&1&-1&-1&-1&-1&1&-1&-1&-1
&1\\
\end{array}
$$
} }

\newpage
Character table of $H_1(\Co_1) = \langle
j_1,k_1,h_1\rangle$(continued){
\renewcommand{\baselinestretch}{0.5}

{\tiny \setlength{\arraycolsep}{0,3mm}
$$
\begin{array}{r|rrrrrrrrrrrrrrrrrrrrrrrrrrrr}
2&6&6&6&6&6&4&5&5&5&11&11&11&11&9&10&10&11&11&11&11&9&9&10&10&10
&10&10\\
3&1&1&1&1&1&2&1&1&1&1&1&1&1&2&1&1&.&.&.&.&1&1&.&.&.&.&.\\
\hline
&6_{46}&6_{47}&6_{48}&6_{49}&6_{50}&6_{51}&6_{52}&6_{53}&6_{54}
&8_{1}&8_{2}&8_{3}&8_{4}&8_{5}&8_{6}&8_{7}&8_{8}&8_{9}&8_{10}&8_{11}
&8_{12}&8_{13}&8_{14}&8_{15}&8_{16}&8_{17}&8_{18}\\
\hline
2P&3e&3d&3d&3d&3b&3f&3c&3d&3b&4_{1}&4_{2}&4_{1}&4_{1}&4_{3}&4_{5}
&4_{4}&4_{2}&4_{1}&4_{1}&4_{1}&4_{3}&4_{3}&4_{4}&4_{17}&4_{17}&4_{4}
&4_{17}\\
3P&2_{13}&2_{23}&2_{25}&2_{26}&2_{27}&2_{8}&2_{30}&2_{31}&2_{32}
&8_{1}&8_{2}&8_{3}&8_{4}&8_{5}&8_{6}&8_{7}&8_{8}&8_{9}&8_{10}&8_{11}
&8_{12}&8_{13}&8_{14}&8_{15}&8_{16}&8_{17}&8_{18}\\
\hline
X.1&1&1&1&1&1&1&1&1&1&1&1&1&1&1&1&1&1&1&1&1&1&1&1&1&1&1&1\\
X.2&1&1&1&1&-1&1&-1&1&-1&-1&1&1&-1&1&1&1&1&1&-1&-1&-1&1&-1&-1&-1&-1
&-1\\
X.3&1&1&-1&-1&-1&1&-1&-1&1&-1&1&1&-1&-1&-1&-1&1&1&-1&-1&1&-1&1&1&-1
&1&-1\\
\end{array}
$$
} }

\newpage
Character table of $H_1(\Co_1) = \langle
j_1,k_1,h_1\rangle$(continued){
\renewcommand{\baselinestretch}{0.5}

{\tiny \setlength{\arraycolsep}{0,3mm}
$$
\begin{array}{r|rrrrrrrrrrrrrrrrrrrrrrrrrr}
2&10&10&10&10&10&10&10&10&10&8&8&9&9&9&9&9&9&9&9&9&9&9&9&8&8\\
3&.&.&.&.&.&.&.&.&.&1&1&.&.&.&.&.&.&.&.&.&.&.&.&.&.\\
\hline
&8_{19}&8_{20}&8_{21}&8_{22}&8_{23}&8_{24}&8_{25}&8_{26}&8_{27}
&8_{28}&8_{29}&8_{30}&8_{31}&8_{32}&8_{33}&8_{34}&8_{35}&8_{36}&8_{37}
&8_{38}&8_{39}&8_{40}&8_{41}&8_{42}&8_{43}\\
\hline
2P&4_{4}&4_{12}&4_{4}&4_{12}&4_{17}&4_{5}&4_{12}&4_{12}&4_{4}
&4_{13}&4_{14}&4_{1}&4_{4}&4_{3}&4_{17}&4_{18}&4_{4}&4_{18}&4_{18}
&4_{12}&4_{17}&4_{12}&4_{18}&4_{27}&4_{28}\\
3P&8_{19}&8_{20}&8_{21}&8_{22}&8_{23}&8_{24}&8_{25}&8_{26}&8_{27}
&8_{28}&8_{29}&8_{30}&8_{31}&8_{32}&8_{33}&8_{34}&8_{35}&8_{36}&8_{37}
&8_{38}&8_{39}&8_{40}&8_{41}&8_{42}&8_{43}\\
\hline
X.1&1&1&1&1&1&1&1&1&1&1&1&1&1&1&1&1&1&1&1&1&1&1&1&1&1\\
X.2&1&-1&-1&-1&-1&1&-1&-1&-1&1&1&-1&-1&-1&-1&-1&-1&1&-1&-1&-1&1&1
&-1&1\\
X.3&-1&-1&-1&-1&1&-1&1&1&-1&-1&1&1&1&1&1&1&-1&-1&1&-1&-1&1&-1&-1&-1\\
\end{array}
$$
} }

\newpage
Character table of $H_1(\Co_1) = \langle
j_1,k_1,h_1\rangle$(continued){
\renewcommand{\baselinestretch}{0.5}

{\tiny \setlength{\arraycolsep}{0,3mm}
$$
\begin{array}{r|rrrrrrrrrrrrrrrrrrrrrrrrrr}
2&8&8&8&8&8&8&8&8&8&8&8&8&8&7&7&7&7&7&7&7&7&7&7&10&10\\
3&.&.&.&.&.&.&.&.&.&.&.&.&.&.&.&.&.&.&.&.&.&.&.&2&2\\
\hline
&8_{44}&8_{45}&8_{46}&8_{47}&8_{48}&8_{49}&8_{50}&8_{51}&8_{52}
&8_{53}&8_{54}&8_{55}&8_{56}&8_{57}&8_{58}&8_{59}&8_{60}&8_{61}&8_{62}
&8_{63}&8_{64}&8_{65}&8_{66}&12_{1}&12_{2}\\
\hline
2P&4_{13}&4_{27}&4_{29}&4_{50}&4_{27}&4_{28}&4_{19}&4_{14}&4_{29}
&4_{29}&4_{28}&4_{50}&4_{14}&4_{56}&4_{57}&4_{58}&4_{58}&4_{56}&4_{57}
&4_{49}&4_{51}&4_{6}&4_{23}&6_{1}&6_{1}\\
3P&8_{44}&8_{45}&8_{46}&8_{47}&8_{48}&8_{49}&8_{50}&8_{51}&8_{52}
&8_{53}&8_{54}&8_{55}&8_{56}&8_{57}&8_{58}&8_{59}&8_{60}&8_{61}&8_{62}
&8_{63}&8_{64}&8_{65}&8_{66}&4_{2}&4_{1}\\
\hline
X.1&1&1&1&1&1&1&1&1&1&1&1&1&1&1&1&1&1&1&1&1&1&1&1&1&1\\
X.2&-1&1&1&-1&-1&-1&-1&-1&-1&-1&-1&-1&-1&1&1&-1&-1&-1&-1&-1&-1&-1
&-1&1&1\\
X.3&1&-1&1&1&1&1&-1&1&1&1&-1&1&1&-1&-1&1&1&1&1&-1&1&1&1&1&1\\
\end{array}
$$
} }

\newpage
Character table of $H_1(\Co_1) = \langle
j_1,k_1,h_1\rangle$(continued){
\renewcommand{\baselinestretch}{0.5}

{\tiny \setlength{\arraycolsep}{0,3mm}
$$
\begin{array}{r|rrrrrrrrrrrrrrrrrrrrrr}
2&9&8&8&8&8&9&9&9&7&7&8&8&8&8&6&6&6&6&5&5&5\\
3&2&2&2&2&2&1&1&1&2&2&1&1&1&1&2&2&2&2&2&2&2\\
\hline &12_{3}&12_{4}&12_{5}&12_{6}&12_{7}&12_{8}&12_{9}&12_{10}
&12_{11}&12_{12}&12_{13}&12_{14}&12_{15}&12_{16}&12_{17}&12_{18}&12_{19}
&12_{20}&12_{21}&12_{22}&12_{23}\\
\hline
2P&6_{1}&6_{2}&6_{1}&6_{2}&6_{2}&6_{1}&6_{2}&6_{2}&6_{3}&6_{3}
&6_{2}&6_{1}&6_{2}&6_{2}&6_{5}&6_{4}&6_{6}&6_{5}&6_{8}&6_{20}&6_{7}\\
3P&4_{3}&4_{9}&4_{10}&4_{11}&4_{7}&4_{13}&4_{5}&4_{4}&4_{6}&4_{8}
&4_{21}&4_{22}&4_{24}&4_{26}&4_{1}&4_{3}&4_{3}&4_{2}&4_{6}&4_{7}&4_{55}\\
\hline
X.1&1&1&1&1&1&1&1&1&1&1&1&1&1&1&1&1&1&1&1&1&1\\
X.2&1&1&1&1&1&1&1&1&1&1&1&1&1&1&1&1&1&1&1&1&1\\
X.3&1&-1&-1&-1&1&1&1&1&1&1&-1&-1&1&-1&1&1&1&1&1&1&-1\\
\end{array}
$$
} }

\newpage
Character table of $H_1(\Co_1) = \langle
j_1,k_1,h_1\rangle$(continued){
\renewcommand{\baselinestretch}{0.5}

{\tiny \setlength{\arraycolsep}{0,3mm}
$$
\begin{array}{r|rrrrrrrrrrrrrrrrrrrrr}
2&5&5&5&5&5&5&5&5&5&5&5&5&5&5&5&5&6&6&4&4\\
3&2&2&2&2&2&2&2&2&2&2&2&2&2&2&2&2&1&1&2&2\\
\hline
&12_{24}&12_{25}&12_{26}&12_{27}&12_{28}&12_{29}&12_{30}&12_{31}
&12_{32}&12_{33}&12_{34}&12_{35}&12_{36}&12_{37}&12_{38}&12_{39}&12_{40}
&12_{41}&12_{42}&12_{43}\\
\hline
2P&6_{8}&6_{21}&6_{6}&6_{14}&6_{14}&6_{21}&6_{15}&6_{17}&6_{17}
&6_{21}&6_{6}&6_{5}&6_{4}&6_{20}&6_{15}&6_{20}&6_{6}&6_{4}&6_{27}
&6_{28}\\
3P&4_{8}&4_{9}&4_{10}&4_{6}&4_{8}&4_{11}&4_{8}&4_{1}&4_{2}&4_{7}
&4_{16}&4_{16}&4_{10}&4_{11}&4_{6}&4_{9}&4_{20}&4_{25}&4_{55}&4_{55}\\
\hline
X.1&1&1&1&1&1&1&1&1&1&1&1&1&1&1&1&1&1&1&1&1\\
X.2&1&1&1&1&1&1&1&1&1&1&1&1&1&1&1&1&1&-1&1&1\\
X.3&1&-1&-1&1&1&-1&1&1&1&1&1&1&-1&-1&1&-1&1&-1&-1&-1\\
\end{array}
$$
} }

\newpage
Character table of $H_1(\Co_1) = \langle
j_1,k_1,h_1\rangle$(continued){
\renewcommand{\baselinestretch}{0.5}

{\tiny \setlength{\arraycolsep}{0,3mm}
$$
\begin{array}{r|rrrrrrrrrrrrrrrrrrrrr}
2&4&5&5&5&5&5&5&5&5&5&5&5&5&5&5&5&5&4&4&4\\
3&2&1&1&1&1&1&1&1&1&1&1&1&1&1&1&1&1&1&1&1\\
\hline
&12_{44}&12_{45}&12_{46}&12_{47}&12_{48}&12_{49}&12_{50}&12_{51}
&12_{52}&12_{53}&12_{54}&12_{55}&12_{56}&12_{57}&12_{58}&12_{59}&12_{60}
&12_{61}&12_{62}&12_{63}\\
\hline
2P&6_{9}&6_{4}&6_{9}&6_{5}&6_{34}&6_{34}&6_{18}&6_{21}&6_{34}
&6_{34}&6_{20}&6_{20}&6_{20}&6_{21}&6_{21}&6_{34}&6_{34}&6_{45}&6_{27}
&6_{45}\\
3P&4_{16}&4_{36}&4_{14}&4_{37}&4_{38}&4_{39}&4_{113}&4_{40}&4_{41}
&4_{15}&4_{42}&4_{43}&4_{44}&4_{45}&4_{46}&4_{47}&4_{48}&4_{104}
&4_{105}&4_{106}\\
\hline
X.1&1&1&1&1&1&1&1&1&1&1&1&1&1&1&1&1&1&1&1&1\\
X.2&1&-1&1&-1&-1&-1&1&1&-1&-1&-1&-1&-1&1&1&1&1&-1&1&-1\\
X.3&1&1&1&-1&-1&-1&-1&-1&-1&-1&1&-1&1&1&-1&1&1&-1&-1&-1\\
\end{array}
$$
} }

\newpage
Character table of $H_1(\Co_1) = \langle
j_1,k_1,h_1\rangle$(continued){
\renewcommand{\baselinestretch}{0.5}

{\tiny \setlength{\arraycolsep}{0,3mm}
$$
\begin{array}{r|rrrrrrrrrrrrrrrrrr}
2&4&4&6&6&5&6&6&4&4&5&5&5&4&4&4&4&3\\
3&1&1&.&.&2&1&1&2&2&1&1&1&1&1&1&1&1\\
\hline
&12_{64}&12_{65}&16a&16b&24a&24b&24c&24d&24e&24f&24g&24h&24i&24j
&24k&24l&24m\\
\hline
2P&6_{28}&6_{16}&8_{2}&8_{8}&12_{3}&12_{1}&12_{2}&12_{19}&12_{18}
&12_{8}&12_{9}&12_{10}&12_{17}&12_{18}&12_{19}&12_{17}&12_{46}\\
3P&4_{107}&4_{23}&16a&16b&8_{5}&8_{2}&8_{3}&8_{5}&8_{5}&8_{28}
&8_{6}&8_{7}&8_{1}&8_{12}&8_{13}&8_{4}&8_{29}\\
\hline
X.1&1&1&1&1&1&1&1&1&1&1&1&1&1&1&1&1&1\\
X.2&-1&1&-1&-1&1&1&1&1&1&1&1&1&-1&-1&1&-1&1\\
X.3&1&1&1&1&-1&1&1&-1&-1&-1&-1&-1&-1&1&-1&-1&1\\
\end{array}
$$
} }

\newpage
Character table of $H_1(\Co_1) = \langle
j_1,k_1,h_1\rangle$(continued){
\renewcommand{\baselinestretch}{0.5}

{\tiny \setlength{\arraycolsep}{0,3mm}
$$
\begin{array}{r|rrrrrrrrrrrrrrrrrrrrrrr}
2&21&21&20&18&18&18&15&14&14&14&17&17&15&15&15&13&13&16&16&14&14
&14\\
3&3&3&2&3&1&1&2&2&2&2&.&.&1&1&1&2&2&.&.&1&1&1\\
\hline
&1a&2_{1}&2_{2}&2_{3}&2_{4}&2_{5}&2_{6}&2_{7}&2_{8}&2_{9}&2_{10}
&2_{11}&2_{12}&2_{13}&2_{14}&2_{15}&2_{16}&2_{17}&2_{18}&2_{19}&2_{20}
&2_{21}\\
\hline
2P&1a&1a&1a&1a&1a&1a&1a&1a&1a&1a&1a&1a&1a&1a&1a&1a&1a&1a&1a&1a&1a
&1a\\
3P&1a&2_{1}&2_{2}&2_{3}&2_{4}&2_{5}&2_{6}&2_{7}&2_{8}&2_{9}&2_{10}
&2_{11}&2_{12}&2_{13}&2_{14}&2_{15}&2_{16}&2_{17}&2_{18}&2_{19}&2_{20}
&2_{21}\\
\hline
X.95&27&27&27&27&27&27&27&9&-9&-9&11&11&3&3&3&9&9&11&11&-9&3&9\\
X.96&32&32&32&-32&32&-32&32&.&8&-8&.&.&.&.&.&.&.&.&.&8&.&.\\
X.97&32&32&32&-32&32&-32&32&8&.&8&.&.&.&.&.&-8&-8&.&.&.&.&8\\
\end{array}
$$
} }

\newpage
Character table of $H_1(\Co_1) = \langle
j_1,k_1,h_1\rangle$(continued){
\renewcommand{\baselinestretch}{0.5}

{\tiny \setlength{\arraycolsep}{0,3mm}
$$
\begin{array}{r|rrrrrrrrrrrrrrrrrrrrrrrr}
2&15&13&13&13&13&13&14&14&12&12&12&13&11&13&8&8&8&7&7&6&17&17&15\\
3&.&1&1&1&1&1&.&.&1&1&1&.&.&3&3&3&3&3&3&3&2&2&2\\
\hline
&2_{22}&2_{23}&2_{24}&2_{25}&2_{26}&2_{27}&2_{28}&2_{29}&2_{30}
&2_{31}&2_{32}&2_{33}&2_{34}&3a&3b&3c&3d&3e&3f&3g&4_{1}&4_{2}&4_{3}\\
\hline
2P&1a&1a&1a&1a&1a&1a&1a&1a&1a&1a&1a&1a&1a&3a&3b&3c&3d&3e&3f&3g
&2_{1}&2_{1}&2_{1}\\
3P&2_{22}&2_{23}&2_{24}&2_{25}&2_{26}&2_{27}&2_{28}&2_{29}&2_{30}
&2_{31}&2_{32}&2_{33}&2_{34}&1a&1a&1a&1a&1a&1a&1a&4_{1}&4_{2}&4_{3}\\
\hline
X.95&11&-9&-9&3&3&9&-5&-9&-3&3&-3&-3&1&.&.&.&.&.&.&.&27&27&27\\
X.96&.&8&8&-8&-8&.&.&-8&.&8&.&.&.&-4&2&2&-1&8&-4&2&32&-32&32\\
X.97&.&.&-8&.&.&8&.&8&.&.&-8&.&.&-4&-1&2&2&-4&8&2&32&-32&32\\
\end{array}
$$
} }

\newpage
Character table of $H_1(\Co_1) = \langle
j_1,k_1,h_1\rangle$(continued){
\renewcommand{\baselinestretch}{0.5}

{\tiny \setlength{\arraycolsep}{0,3mm}
$$
\begin{array}{r|rrrrrrrrrrrrrrrrrrrrrr}
2&16&16&14&14&14&14&14&14&16&14&14&14&12&15&15&15&13&13&13&13&13\\
3&1&1&2&2&2&2&2&2&.&1&1&1&2&.&.&.&1&1&1&1&1\\
\hline
&4_{4}&4_{5}&4_{6}&4_{7}&4_{8}&4_{9}&4_{10}&4_{11}&4_{12}&4_{13}
&4_{14}&4_{15}&4_{16}&4_{17}&4_{18}&4_{19}&4_{20}&4_{21}&4_{22}&4_{23}
&4_{24}\\
\hline
2P&2_{2}&2_{2}&2_{3}&2_{2}&2_{3}&2_{2}&2_{1}&2_{2}&2_{1}&2_{1}
&2_{1}&2_{4}&2_{1}&2_{2}&2_{1}&2_{1}&2_{1}&2_{2}&2_{1}&2_{3}&2_{2}\\
3P&4_{4}&4_{5}&4_{6}&4_{7}&4_{8}&4_{9}&4_{10}&4_{11}&4_{12}&4_{13}
&4_{14}&4_{15}&4_{16}&4_{17}&4_{18}&4_{19}&4_{20}&4_{21}&4_{22}&4_{23}
&4_{24}\\
\hline
X.95&3&3&-9&27&-9&-9&-9&-9&11&3&3&-3&-9&3&11&3&-9&-9&-9&-9&3\\
X.96&.&.&.&-32&.&8&8&-8&.&.&.&.&-8&.&.&.&8&8&-8&.&.\\
X.97&.&.&.&-32&.&-8&-8&8&.&.&.&.&.&.&.&.&.&-8&8&.&.\\
\end{array}
$$
} }

\newpage
Character table of $H_1(\Co_1) = \langle
j_1,k_1,h_1\rangle$(continued){
\renewcommand{\baselinestretch}{0.5}

{\tiny \setlength{\arraycolsep}{0,3mm}
$$
\begin{array}{r|rrrrrrrrrrrrrrrrrrrrrrr}
2&13&13&14&14&14&14&14&14&14&14&14&12&12&12&12&12&12&12&12&12&12
&12\\
3&1&1&.&.&.&.&.&.&.&.&.&1&1&1&1&1&1&1&1&1&1&1\\
\hline
&4_{25}&4_{26}&4_{27}&4_{28}&4_{29}&4_{30}&4_{31}&4_{32}&4_{33}
&4_{34}&4_{35}&4_{36}&4_{37}&4_{38}&4_{39}&4_{40}&4_{41}&4_{42}&4_{43}
&4_{44}&4_{45}&4_{46}\\
\hline
2P&2_{1}&2_{2}&2_{2}&2_{2}&2_{2}&2_{1}&2_{2}&2_{4}&2_{2}&2_{2}
&2_{2}&2_{1}&2_{1}&2_{4}&2_{4}&2_{2}&2_{4}&2_{2}&2_{2}&2_{2}&2_{2}
&2_{2}\\
3P&4_{25}&4_{26}&4_{27}&4_{28}&4_{29}&4_{30}&4_{31}&4_{32}&4_{33}
&4_{34}&4_{35}&4_{36}&4_{37}&4_{38}&4_{39}&4_{40}&4_{41}&4_{42}&4_{43}
&4_{44}&4_{45}&4_{46}\\
\hline
X.95&9&-9&-5&3&3&-9&11&-3&11&-9&-9&-3&-3&9&-3&3&9&-3&9&-3&-9&3\\
X.96&.&-8&.&.&.&8&.&.&.&8&-8&.&.&.&.&8&.&.&.&.&-8&-8\\
X.97&8&8&.&.&.&-8&.&.&.&-8&8&8&.&8&.&.&-8&8&-8&-8&.&.\\
\end{array}
$$
} }

\newpage
Character table of $H_1(\Co_1) = \langle
j_1,k_1,h_1\rangle$(continued){
\renewcommand{\baselinestretch}{0.5}

{\tiny \setlength{\arraycolsep}{0,3mm}
$$
\begin{array}{r|rrrrrrrrrrrrrrrrrrrrrrrr}
2&12&12&13&13&13&13&13&13&9&12&12&12&12&12&12&12&12&12&12&12&12
&12&12\\
3&1&1&.&.&.&.&.&.&2&.&.&.&.&.&.&.&.&.&.&.&.&.&.\\
\hline
&4_{47}&4_{48}&4_{49}&4_{50}&4_{51}&4_{52}&4_{53}&4_{54}&4_{55}
&4_{56}&4_{57}&4_{58}&4_{59}&4_{60}&4_{61}&4_{62}&4_{63}&4_{64}&4_{65}
&4_{66}&4_{67}&4_{68}&4_{69}\\
\hline
2P&2_{4}&2_{4}&2_{2}&2_{3}&2_{3}&2_{10}&2_{1}&2_{10}&2_{6}&2_{4}
&2_{4}&2_{4}&2_{2}&2_{10}&2_{4}&2_{3}&2_{4}&2_{4}&2_{10}&2_{4}&2_{4}
&2_{4}&2_{10}\\
3P&4_{47}&4_{48}&4_{49}&4_{50}&4_{51}&4_{52}&4_{53}&4_{54}&4_{55}
&4_{56}&4_{57}&4_{58}&4_{59}&4_{60}&4_{61}&4_{62}&4_{63}&4_{64}&4_{65}
&4_{66}&4_{67}&4_{68}&4_{69}\\
\hline
X.95&-9&-9&-5&-1&-1&1&-3&1&-9&3&-5&3&-3&1&11&-1&-9&-3&1&9&-3&-9&1\\
X.96&8&-8&.&.&.&.&.&.&-8&.&.&.&.&.&.&.&8&.&.&.&.&-8&.\\
X.97&.&.&.&.&.&.&.&.&8&.&.&.&.&.&.&.&.&.&.&8&.&.&.\\
\end{array}
$$
} }

\newpage
Character table of $H_1(\Co_1) = \langle
j_1,k_1,h_1\rangle$(continued){
\renewcommand{\baselinestretch}{0.5}

{\tiny \setlength{\arraycolsep}{0,3mm}
$$
\begin{array}{r|rrrrrrrrrrrrrrrrrrrrrrrrrr}
2&12&11&11&11&11&11&11&11&11&11&11&11&11&11&11&11&11&11&11&11&11
&11&11&11&11\\
3&.&.&.&.&.&.&.&.&.&.&.&.&.&.&.&.&.&.&.&.&.&.&.&.&.\\
\hline
&4_{70}&4_{71}&4_{72}&4_{73}&4_{74}&4_{75}&4_{76}&4_{77}&4_{78}
&4_{79}&4_{80}&4_{81}&4_{82}&4_{83}&4_{84}&4_{85}&4_{86}&4_{87}&4_{88}
&4_{89}&4_{90}&4_{91}&4_{92}&4_{93}&4_{94}\\
\hline
2P&2_{4}&2_{4}&2_{10}&2_{10}&2_{10}&2_{10}&2_{10}&2_{2}&2_{10}
&2_{10}&2_{10}&2_{10}&2_{10}&2_{10}&2_{4}&2_{10}&2_{10}&2_{10}&2_{10}
&2_{10}&2_{10}&2_{10}&2_{4}&2_{4}&2_{4}\\
3P&4_{70}&4_{71}&4_{72}&4_{73}&4_{74}&4_{75}&4_{76}&4_{77}&4_{78}
&4_{79}&4_{80}&4_{81}&4_{82}&4_{83}&4_{84}&4_{85}&4_{86}&4_{87}&4_{88}
&4_{89}&4_{90}&4_{91}&4_{92}&4_{93}&4_{94}\\
\hline
X.95&9&-3&3&3&-1&1&-3&1&-1&-1&3&3&-1&-3&1&1&-3&-3&-3&1&1&-3&1&3&3\\
X.96&.&.&.&.&.&.&.&.&.&.&.&.&.&.&.&.&.&.&.&.&.&.&.&8&-8\\
X.97&-8&8&.&.&.&.&.&.&.&.&.&.&.&.&.&.&.&.&.&.&.&.&.&.&.\\
\end{array}
$$
} }

\newpage
Character table of $H_1(\Co_1) = \langle
j_1,k_1,h_1\rangle$(continued){
\renewcommand{\baselinestretch}{0.5}

{\tiny \setlength{\arraycolsep}{0,3mm}
$$
\begin{array}{r|rrrrrrrrrrrrrrrrrrrrrr}
2&11&11&11&11&11&11&11&11&11&9&9&9&9&10&10&10&10&10&8&9&9\\
3&.&.&.&.&.&.&.&.&.&1&1&1&1&.&.&.&.&.&1&.&.\\
\hline &4_{95}&4_{96}&4_{97}&4_{98}&4_{99}&4_{100}&4_{101}&4_{102}
&4_{103}&4_{104}&4_{105}&4_{106}&4_{107}&4_{108}&4_{109}&4_{110}&4_{111}
&4_{112}&4_{113}&4_{114}&4_{115}\\
\hline
2P&2_{10}&2_{10}&2_{10}&2_{10}&2_{10}&2_{10}&2_{10}&2_{4}&2_{10}
&2_{12}&2_{6}&2_{12}&2_{6}&2_{17}&2_{18}&2_{5}&2_{18}&2_{17}&2_{20}
&2_{12}&2_{6}\\
3P&4_{95}&4_{96}&4_{97}&4_{98}&4_{99}&4_{100}&4_{101}&4_{102}
&4_{103}&4_{104}&4_{105}&4_{106}&4_{107}&4_{108}&4_{109}&4_{110}&4_{111}
&4_{112}&4_{113}&4_{114}&4_{115}\\
\hline
X.95&-1&-1&-1&-1&-3&-3&1&-3&1&-3&3&-3&-3&1&-1&-1&-1&1&3&-1&1\\
X.96&.&.&.&.&.&.&.&.&.&.&-8&.&.&.&.&.&.&.&.&.&.\\
X.97&.&.&.&.&.&.&.&-8&.&.&.&.&8&.&.&.&.&.&.&.&.\\
\end{array}
$$
} }

\newpage
Character table of $H_1(\Co_1) = \langle
j_1,k_1,h_1\rangle$(continued){
\renewcommand{\baselinestretch}{0.5}

{\tiny \setlength{\arraycolsep}{0,3mm}
$$
\begin{array}{r|rrrrrrrrrrrrrrrrrrrrr}
2&9&9&9&9&9&9&9&9&9&9&9&9&9&9&8&8&8&8&8&13\\
3&.&.&.&.&.&.&.&.&.&.&.&.&.&.&.&.&.&.&.&3\\
\hline
&4_{116}&4_{117}&4_{118}&4_{119}&4_{120}&4_{121}&4_{122}&4_{123}
&4_{124}&4_{125}&4_{126}&4_{127}&4_{128}&4_{129}&4_{130}&4_{131}&4_{132}
&4_{133}&4_{134}&6_{1}\\
\hline
2P&2_{22}&2_{12}&2_{17}&2_{22}&2_{22}&2_{12}&2_{12}&2_{17}&2_{18}
&2_{17}&2_{18}&2_{22}&2_{11}&2_{18}&2_{28}&2_{20}&2_{14}&2_{28}&2_{28}
&3a\\
3P&4_{116}&4_{117}&4_{118}&4_{119}&4_{120}&4_{121}&4_{122}&4_{123}
&4_{124}&4_{125}&4_{126}&4_{127}&4_{128}&4_{129}&4_{130}&4_{131}&4_{132}
&4_{133}&4_{134}&2_{1}\\
\hline
X.95&3&-1&-3&1&-3&1&1&-1&3&-3&-3&-1&-1&1&-1&1&1&1&1&.\\
X.96&.&.&.&.&.&.&.&.&.&.&.&.&.&.&.&.&.&.&.&-4\\
X.97&.&.&.&.&.&.&.&.&.&.&.&.&.&.&.&.&.&.&.&-4\\
\end{array}
$$
} }

\newpage
Character table of $H_1(\Co_1) = \langle
j_1,k_1,h_1\rangle$(continued){
\renewcommand{\baselinestretch}{0.5}

{\tiny \setlength{\arraycolsep}{0,3mm}
$$
\begin{array}{r|rrrrrrrrrrrrrrrrrrrrrrrrrr}
2&12&10&8&8&8&9&7&7&7&7&7&8&6&6&6&6&9&7&7&7&7&5&5&8&8\\
3&2&3&3&3&3&2&3&3&3&3&3&2&3&3&3&3&1&2&2&2&2&3&3&1&1\\
\hline
&6_{2}&6_{3}&6_{4}&6_{5}&6_{6}&6_{7}&6_{8}&6_{9}&6_{10}&6_{11}
&6_{12}&6_{13}&6_{14}&6_{15}&6_{16}&6_{17}&6_{18}&6_{19}&6_{20}&6_{21}
&6_{22}&6_{23}&6_{24}&6_{25}&6_{26}\\
\hline
2P&3a&3a&3b&3c&3d&3a&3e&3f&3e&3e&3c&3a&3f&3g&3g&3g&3a&3d&3b&3d&3b
&3d&3b&3a&3a\\
3P&2_{2}&2_{3}&2_{1}&2_{1}&2_{1}&2_{6}&2_{3}&2_{1}&2_{1}&2_{3}
&2_{3}&2_{9}&2_{3}&2_{3}&2_{3}&2_{1}&2_{20}&2_{8}&2_{2}&2_{2}&2_{7}
&2_{3}&2_{3}&2_{14}&2_{24}\\
\hline
X.95&.&.&.&.&.&.&.&.&.&.&.&.&.&.&.&.&.&.&.&.&.&.&.&.&.\\
X.96&-4&4&2&2&-1&-4&-8&-4&8&-8&-2&-2&4&-2&-2&2&.&-1&2&-1&.&1&-2&.&2\\
X.97&-4&4&-1&2&2&-4&4&8&-4&4&-2&2&-8&-2&-2&2&.&.&-1&2&-1&-2&1&.&-2\\
\end{array}
$$
} }

\newpage
Character table of $H_1(\Co_1) = \langle
j_1,k_1,h_1\rangle$(continued){
\renewcommand{\baselinestretch}{0.5}

{\tiny \setlength{\arraycolsep}{0,3mm}
$$
\begin{array}{r|rrrrrrrrrrrrrrrrrrrrrrrrrr}
2&6&6&6&6&6&6&7&7&7&7&5&5&5&5&5&5&5&5&6&6&6&6&6&6&4\\
3&2&2&2&2&2&2&1&1&1&1&2&2&2&2&2&2&2&2&1&1&1&1&1&1&2\\
\hline
&6_{27}&6_{28}&6_{29}&6_{30}&6_{31}&6_{32}&6_{33}&6_{34}&6_{35}
&6_{36}&6_{37}&6_{38}&6_{39}&6_{40}&6_{41}&6_{42}&6_{43}&6_{44}&6_{45}
&6_{46}&6_{47}&6_{48}&6_{49}&6_{50}&6_{51}\\
\hline
2P&3d&3b&3b&3b&3c&3c&3d&3c&3b&3c&3e&3d&3c&3e&3c&3e&3e&3b&3e&3e&3d
&3d&3d&3b&3f\\
3P&2_{6}&2_{6}&2_{15}&2_{16}&2_{16}&2_{15}&2_{19}&2_{4}&2_{21}
&2_{5}&2_{7}&2_{9}&2_{7}&2_{7}&2_{8}&2_{16}&2_{15}&2_{9}&2_{12}&2_{13}
&2_{23}&2_{25}&2_{26}&2_{27}&2_{8}\\
\hline
X.95&.&.&.&.&.&.&.&.&.&.&.&.&.&.&.&.&.&.&.&.&.&.&.&.&.\\
X.96&-1&2&.&.&.&.&-1&2&.&-2&.&1&.&.&2&.&.&-2&.&.&-1&1&1&.&2\\
X.97&2&-1&1&1&-2&-2&.&2&-1&-2&2&2&2&2&.&-2&-2&-1&.&.&.&.&.&-1&.\\
\end{array}
$$
} }

\newpage
Character table of $H_1(\Co_1) = \langle
j_1,k_1,h_1\rangle$(continued){
\renewcommand{\baselinestretch}{0.5}

{\tiny \setlength{\arraycolsep}{0,3mm}
$$
\begin{array}{r|rrrrrrrrrrrrrrrrrrrrrrrrr}
2&5&5&5&11&11&11&11&9&10&10&11&11&11&11&9&9&10&10&10&10&10&10&10
&10\\
3&1&1&1&1&1&1&1&2&1&1&.&.&.&.&1&1&.&.&.&.&.&.&.&.\\
\hline
&6_{52}&6_{53}&6_{54}&8_{1}&8_{2}&8_{3}&8_{4}&8_{5}&8_{6}&8_{7}
&8_{8}&8_{9}&8_{10}&8_{11}&8_{12}&8_{13}&8_{14}&8_{15}&8_{16}&8_{17}
&8_{18}&8_{19}&8_{20}&8_{21}\\
\hline
2P&3c&3d&3b&4_{1}&4_{2}&4_{1}&4_{1}&4_{3}&4_{5}&4_{4}&4_{2}&4_{1}
&4_{1}&4_{1}&4_{3}&4_{3}&4_{4}&4_{17}&4_{17}&4_{4}&4_{17}&4_{4}&4_{12}
&4_{4}\\
3P&2_{30}&2_{31}&2_{32}&8_{1}&8_{2}&8_{3}&8_{4}&8_{5}&8_{6}&8_{7}
&8_{8}&8_{9}&8_{10}&8_{11}&8_{12}&8_{13}&8_{14}&8_{15}&8_{16}&8_{17}
&8_{18}&8_{19}&8_{20}&8_{21}\\
\hline
X.95&.&.&.&-3&-9&3&-3&-9&3&3&-1&-5&-3&-3&-3&3&1&1&1&1&1&3&-3&1\\
X.96&.&-1&.&.&.&.&.&8&.&.&.&.&.&.&.&8&.&.&.&.&.&.&.&.\\
X.97&.&.&1&.&.&.&.&-8&.&.&.&.&.&.&-8&.&.&.&.&.&.&.&.&.\\
\end{array}
$$
} }

\newpage
Character table of $H_1(\Co_1) = \langle
j_1,k_1,h_1\rangle$(continued){
\renewcommand{\baselinestretch}{0.5}

{\tiny \setlength{\arraycolsep}{0,3mm}
$$
\begin{array}{r|rrrrrrrrrrrrrrrrrrrrrrrrrr}
2&10&10&10&10&10&10&8&8&9&9&9&9&9&9&9&9&9&9&9&9&8&8&8&8&8\\
3&.&.&.&.&.&.&1&1&.&.&.&.&.&.&.&.&.&.&.&.&.&.&.&.&.\\
\hline
&8_{22}&8_{23}&8_{24}&8_{25}&8_{26}&8_{27}&8_{28}&8_{29}&8_{30}
&8_{31}&8_{32}&8_{33}&8_{34}&8_{35}&8_{36}&8_{37}&8_{38}&8_{39}&8_{40}
&8_{41}&8_{42}&8_{43}&8_{44}&8_{45}&8_{46}\\
\hline
2P&4_{12}&4_{17}&4_{5}&4_{12}&4_{12}&4_{4}&4_{13}&4_{14}&4_{1}
&4_{4}&4_{3}&4_{17}&4_{18}&4_{4}&4_{18}&4_{18}&4_{12}&4_{17}&4_{12}
&4_{18}&4_{27}&4_{28}&4_{13}&4_{27}&4_{29}\\
3P&8_{22}&8_{23}&8_{24}&8_{25}&8_{26}&8_{27}&8_{28}&8_{29}&8_{30}
&8_{31}&8_{32}&8_{33}&8_{34}&8_{35}&8_{36}&8_{37}&8_{38}&8_{39}&8_{40}
&8_{41}&8_{42}&8_{43}&8_{44}&8_{45}&8_{46}\\
\hline
X.95&-3&1&3&1&1&1&3&3&1&1&1&1&1&1&3&-3&1&1&-1&-1&1&-1&1&-1&3\\
X.96&.&.&.&.&.&.&.&.&.&.&.&.&.&.&.&.&.&.&.&.&.&.&.&.&.\\
X.97&.&.&.&.&.&.&.&.&.&.&.&.&.&.&.&.&.&.&.&.&.&.&.&.&.\\
\end{array}
$$
} }

\newpage
Character table of $H_1(\Co_1) = \langle
j_1,k_1,h_1\rangle$(continued){
\renewcommand{\baselinestretch}{0.5}

{\tiny \setlength{\arraycolsep}{0,3mm}
$$
\begin{array}{r|rrrrrrrrrrrrrrrrrrrrrrrrrr}
2&8&8&8&8&8&8&8&8&8&8&7&7&7&7&7&7&7&7&7&7&10&10&9&8&8\\
3&.&.&.&.&.&.&.&.&.&.&.&.&.&.&.&.&.&.&.&.&2&2&2&2&2\\
\hline
&8_{47}&8_{48}&8_{49}&8_{50}&8_{51}&8_{52}&8_{53}&8_{54}&8_{55}
&8_{56}&8_{57}&8_{58}&8_{59}&8_{60}&8_{61}&8_{62}&8_{63}&8_{64}&8_{65}
&8_{66}&12_{1}&12_{2}&12_{3}&12_{4}&12_{5}\\
\hline
2P&4_{50}&4_{27}&4_{28}&4_{19}&4_{14}&4_{29}&4_{29}&4_{28}&4_{50}
&4_{14}&4_{56}&4_{57}&4_{58}&4_{58}&4_{56}&4_{57}&4_{49}&4_{51}&4_{6}
&4_{23}&6_{1}&6_{1}&6_{1}&6_{2}&6_{1}\\
3P&8_{47}&8_{48}&8_{49}&8_{50}&8_{51}&8_{52}&8_{53}&8_{54}&8_{55}
&8_{56}&8_{57}&8_{58}&8_{59}&8_{60}&8_{61}&8_{62}&8_{63}&8_{64}&8_{65}
&8_{66}&4_{2}&4_{1}&4_{3}&4_{9}&4_{10}\\
\hline
X.95&-1&1&1&1&1&1&1&-3&-1&1&-1&-1&1&1&1&1&1&-1&-1&-1&.&.&.&.&.\\
X.96&.&.&.&.&.&.&.&.&.&.&.&.&.&.&.&.&.&.&.&.&4&-4&-4&2&2\\
X.97&.&.&.&.&.&.&.&.&.&.&.&.&.&.&.&.&.&.&.&.&4&-4&-4&-2&-2\\
\end{array}%
$$
} }
\end{ctab}

\newpage

\begin{ctab}\label{Co_1ct_U}Character table of $UCo1 = \langle
y,j,k\rangle${
\renewcommand{\baselinestretch}{0.5}

{\tiny \setlength{\arraycolsep}{0,3mm}
$$
%
$$
} }

\newpage
Character table of $UCo1$ (continued){
\renewcommand{\baselinestretch}{0.5}

{\tiny \setlength{\arraycolsep}{0,3mm}
$$
%
$$
} }

\newpage
Character table of $UCo1$ (continued){
\renewcommand{\baselinestretch}{0.5}

{\tiny \setlength{\arraycolsep}{0,3mm}
$$
%
$$
} }

\newpage
Character table of $UCo1$ (continued){
\renewcommand{\baselinestretch}{0.5}

{\tiny \setlength{\arraycolsep}{0,3mm}
$$
%
$$
} }
\end{ctab}

\newpage

\begin{ctab}\label{Co_1ct_D}Character table of $D(\Co_1) = \langle x,y\rangle$
{ \setlength{\arraycolsep}{1mm}
\renewcommand{\baselinestretch}{0.5}

{\tiny
\setlength{\arraycolsep}{0,3mm}
$$
\begin{array}{r|rrrrrrrrrrrrrrrrrrrrrrr}
2&21&21&20&17&18&15&15&18&17&16&16&14&14&14&15&15&13&7&8&17&17&16\\
3&2&2&1&1&2&2&2&1&1&1&.&1&1&1&.&.&.&2&2&1&1&1\\
5&1&1&.&.&.&1&1&.&.&.&.&.&.&.&.&.&.&1&.&.&.&.\\
7&1&1&1&1&.&.&.&.&.&.&.&.&.&.&.&.&.&.&.&1&.&.\\
\hline
&1a&2a&2b&2c&2d&2e&2f&2g&2h&2i&2j&2k&2l&2m&2n&2o&2p&3a&3b&4_{1}
&4_{2}&4_{3}\\
\hline
2P&1a&1a&1a&1a&1a&1a&1a&1a&1a&1a&1a&1a&1a&1a&1a&1a&1a&3a&3b&2a&2a
&2b\\
3P&1a&2a&2b&2c&2d&2e&2f&2g&2h&2i&2j&2k&2l&2m&2n&2o&2p&1a&1a&4_{1}
&4_{2}&4_{3}\\
5P&1a&2a&2b&2c&2d&2e&2f&2g&2h&2i&2j&2k&2l&2m&2n&2o&2p&3a&3b&4_{1}
&4_{2}&4_{3}\\
\end{array}
$$
} }

\newpage
Character table of $D(\Co_1) =\langle x,y \rangle$ (continued){
\renewcommand{\baselinestretch}{0.5}

{\tiny
\setlength{\arraycolsep}{0,3mm}
$$
%
$$
} }

\newpage

Character table of $D(\Co_1)= \langle x,y\rangle$ (continued){
\renewcommand{\baselinestretch}{0.5}

{\tiny
\setlength{\arraycolsep}{0,3mm}
$$
%
$$
} }

\newpage
Character table of $D(\Co_1)= \langle x,y\rangle$ (continued){
\renewcommand{\baselinestretch}{0.5}

{\tiny
\setlength{\arraycolsep}{0,3mm}
$$
%
$$
} }

\newpage
Character table of $D(\Co_1)= \langle x,y\rangle$ (continued){
\renewcommand{\baselinestretch}{0.5}

{\tiny
\setlength{\arraycolsep}{0,3mm}
$$
%
$$
} }, where
$
A = -2\zeta(15)_3\zeta(15)_5^3 - 2\zeta(15)_3\zeta(15)_5^2 - \zeta(15)_3 -
\zeta(15)_5^3 - \zeta(15)_5^2 - 1,
B = -\zeta(7)^4 - \zeta(7)^2 - \zeta(7) - 1,
C = -2\zeta(7)^4 - 2\zeta(7)^2 - 2\zeta(7) - 2
$.
\end{ctab}

\newpage

\begin{ctab}\label{Co_1ct H}Character table of $H(\Co_1) = \langle x,y,h \rangle$
{ \setlength{\arraycolsep}{1mm}
\renewcommand{\baselinestretch}{0.5}

{\tiny
\setlength{\arraycolsep}{0,3mm}
$$
%
$$
} }

\newpage
Character table of $H(\Co_1)= \langle x,y,h \rangle$ (continued){
\renewcommand{\baselinestretch}{0.5}

{\tiny
\setlength{\arraycolsep}{0,3mm}
$$
%
$$
} }

\newpage
Character table of $H(\Co_1)= \langle x,y,h \rangle$ (continued){
\renewcommand{\baselinestretch}{0.5}

{\tiny
\setlength{\arraycolsep}{0,3mm}
$$
%
$$
} }

\newpage
Character table of $H(\Co_1)= \langle x,y,h \rangle$ (continued){
\renewcommand{\baselinestretch}{0.5}

{\tiny
\setlength{\arraycolsep}{0,3mm}
$$
%
$$
} }, where $
A = 4\zeta(15)_3\zeta(15)_5^3 + 4\zeta(15)_3\zeta(15)_5^2 + 2\zeta(15)_3 +
2\zeta(15)_5^3 + 2\zeta(15)_5^2 + 1,
B = 2\zeta(7)^4 + 2\zeta(7)^2 + 2\zeta(7) + 1
$
\end{ctab}


\newpage

\begin{thebibliography}{10}

\bibitem{magma}
John Cannon and Catherine Playoust.
\newblock {\em An Introduction to {\sc Magma}}.
\newblock School of Mathematics and Statistics, University of Sydney, 1993.

\bibitem{cannon}
John~J. Cannon, Derek~F. Holt.
\newblock Automorphism group computation and isomorphism
testing in finite groups.
\newblock {\em J. Symbolic Computat.}, 35:241-267, 2003.

\bibitem{carter}
R.W.\ Carter.
\newblock {\em Simple groups of Lie type}.
\newblock John Wiley and Sons, London, 1972.

\bibitem{conway1}
J.H.\ Conway.
\newblock A group of order 8,315,553,613,086,720,000.
\newblock {\em Bull. London Math. Soc.},1:79-88, 1969.

\bibitem{conway}
J.H.\ Conway.
\newblock Three lectures on exceptional groups, in M.B. Powel, G.
Higman (eds.){\em Finite simple groups}.
\newblock Academic Press, London, 1971, pp. 215
- 247.

\bibitem{atlas}
J.H.\ Conway, R.T.\ Curtis, S.P.\ Norton, R.A.\ Parker, and R.A.\
Wilson.
\newblock {\em Atlas of finite groups}.
\newblock Clarendon Press, Oxford, 1985.

\bibitem{holt2} Holt,D. F.,  B. Eick and E.A. O'Brien, \emph{Handbook of computational group theory},
Chapman and Hall/CRC, Boca Raton, 2005.

\bibitem{holt5} Holt, D. F., Cohomology and group extensions in
Magma, in W. Bosma, J. Cannon (eds), Discovering mathematics with
Magma, Springer, Berlin 2006, pp. 221--141.

\bibitem{james}
G.\ James.
\newblock The modular characters of the Mathieu groups.
\newblock {\em J. Algebra}, 27:57-111, 1973.

\bibitem{kim}
H. Kim, G. O. Michler.
\newblock Simultaneous constructions of the sporadic groups $\Co_2$ and $\Fi_{22}$.
\newblock in (L. -C. Kappe, A. Magidin, R. F. Morse, eds.)
\textit{Computational Group Theory and the Theory of Groups},
Contemporary Mathematics {\bf 470}, 141--234, Amer. Math. Soc,
Providence, RI., (2008). See also $\verb"arXiv:0906.0623v1"$.

\bibitem{kim1}
H. Kim.
\newblock Representation theoretic existence proof for Fischer's sporadic group $\Fi_{23}$.
\newblock Senior Thesis, Dept. Math. Cornell University, 2008. 
$\verb"arXiv:0904.0639v1"$.

\bibitem{kim2}
H. Kim, G. O. Michler.
\newblock Construction of Fischer's sporadic group $\Fi_{24}'$ inside $\GL_{8671}(13)$.
\newblock $\verb"arXiv:0906.1064v1"$.

\bibitem{michler}
G. O. Michler.
\newblock Theory of finite simple groups.
\newblock Cambridge University Press, Cambridge, 2006.

\bibitem{michler1}
G. O. Michler.
\newblock Constructing finite simple groups from irreducible subgroups of
$\GL_n(2)$,
\newblock in (L. -C. Kappe, A. Magidin, R. F. Morse, eds.)
\textit{Computational Group Theory and the Theory of Groups},
Contemporary Mathematics {\bf 470}, 235--262, Amer. Math. Soc, Providence,
RI., (2008).

\bibitem{michler2}
G. O. Michler.
\newblock Theory of finite simple groups II.
\newblock Cambridge University Press, Cambridge (to appear in 2009).


\bibitem{soicher1}
L. Soicher.
\newblock Presentations for Conway's group
$\Co_1$.
\newblock {\em Math. Proceedings Cambridge Phil. Soc.}, 102:1-3, 1987.
\end{thebibliography}
\end{document}